\documentclass[11pt]{article}

\usepackage{latexsym}
\usepackage{amsfonts}
\usepackage{amsmath}

\usepackage[dvips]{color}

\newtheorem{thrm}{Theorem}[section]
\newtheorem{lemma}[thrm]{Lemma}
\newtheorem{prop}[thrm]{Proposition}

\newtheorem{cor}[thrm]{Corollary}
\newtheorem{remark}[thrm]{Remark}

\numberwithin{equation}{section}

\usepackage{enumerate}

\usepackage{amssymb}
\usepackage{mathtools}
\mathtoolsset{showonlyrefs}
\usepackage{mathrsfs}
\usepackage{comment}
\usepackage{hyperref}
\usepackage{xcolor}

\def\E{\mathbb{E} }
\def\P{\mathbb{P} }
\def\Q{\mathbb{Q} }
\def\R{\mathbb{R} }
\def\N{\mathbb{N} }

\makeatletter
\begin{document}
\allowdisplaybreaks

\title{\Large \bf{Extremal process for  irreducible multitype branching Brownian motion}\footnote{The research of this project is supported
     by the National Key R\&D Program of China (No. 2020YFA0712900).}}
\author{ \bf  Haojie Hou \hspace{1mm}\hspace{1mm}
Yan-Xia Ren\footnote{The research of this author is supported by NSFC (Grant Nos. 11731009, 12071011 and 12231002) and LMEQF.
 } \hspace{1mm}\hspace{1mm} and \hspace{1mm}\hspace{1mm}
Renming Song\thanks{Research supported in part by a grant from the Simons
Foundation
(\#960480, Renming Song).}
\hspace{1mm} }
\date{}
\maketitle

\begin{abstract}
We first study the convergence of solutions of a system of F-KPP equations related to irreducible multitype branching Brownian motions with Heaviside-type initial conditions to traveling wave solutions. Then we apply this convergence result to prove that the extremal processes of irreducible multitype branching Brownian motions converges weakly to a cluster point process.
\end{abstract}

\medskip


\medskip


\section{Introduction and notation}

\subsection{Background}

A binary branching Brownian motion (BBM) is a continuous-time Markov process
which can be defined as follows. Initially, there is a particle at the origin and the
particle moves according to a standard Brownian motion.
After an exponential time with parameter $1$, this particle dies and splits into $2$ particles.
The offspring move
independently according to standard Brownian motion from the place they are born and obey the same branching mechanism as their parent.
We denote the law of this branching Brownian motion by $\P$.

The binary branching Brownian motion is related to the F-KPP equation.
Let $M_t$ be the right-most position among all the particles alive at time $t$.
McKean \cite{McKean} proved that the function
\[
u(t,x):= \P \left(M_t\leq x\right),\quad t\geq 0, x\in\R,
\]
solves the F-KPP equation
\begin{align}\label{F-KPP-1d}
	u_t = \frac{1}{2}u_{xx}+ u^2 -u
\end{align}
with the Heaviside initial condition $u(0,x)=1_{[0, \infty)}(x)$. Equation \eqref{F-KPP-1d} was first studied by Fisher \cite{Fisher} and Kolmogorov, Petrovskii and Piskounov \cite{KPP}. Later,
Bramson \cite[Theorems A, B and Example 2]{Bramson1} studied the
asymptotic behaviors of
solutions of \eqref{F-KPP-1d}
for a class of more general initial conditions.
Let $u$ be a solution of \eqref{F-KPP-1d} and $v= 1-u$.
Bramson proved that, under some conditions on $v(0,x)$,
\[
v\left(t,\sqrt{2}t - \frac{3}{2\sqrt{2}}\log t +x\right)\to 1-w(x),\quad \mbox{uniformly in }x\ \mbox{as}\ t\to\infty,
\]
here $w$ is the unique solution (up to a translation) of
\[
\frac{1}{2}w'' + \sqrt{2}w'+w^2 -w=0
\]
and $w$ is called a traveling wave solution.
In the Heaviside case, a probabilistic representation of the limit $w(x)$ was given by Lalley and Sellke \cite{L-S}.
For different proofs of this result,
see \cite{BDZ, Roberts}.

The extremal point process of branching Brownian motion has also been widely studied.
A\"{i}d\'ekon et. al \cite{ABBS} and Arguin et. al \cite{ABK}
studied this extremal point process
using different methods.
Suppose that the set of the positions of all particles alive at time $t$ is given by $\{X_u(t): u\in Z(t)\}$, where $Z(t)$ is the set of particles alive at time $t$.
It is known that
$$
D_{t}=\sum_{u \in Z_t}\left(\sqrt{2} t-{X}_u(t)\right) e^{\sqrt{2} {X}_u(t)-2 t},\quad  t\geq 0,
$$
is a martingale and has (non-negative) limit $D_\infty$ as $t\to\infty$.
As $t\to\infty$, the extremal point process
\[
\sum_{u\in Z(t)}\delta_{X_u(t)- \left(\sqrt{2}t-\frac{3}{2\sqrt{2}}\log t \right)}
\]
converges in distribution to a decorated Poisson point process DPPP($CD_\infty e^{-\sqrt{2}x}\mathrm{d}x, \mathcal{D}_0$),
in the sense of
vague topology,
where $\mathcal{D}_0$ is a point process.
More precisely, this limit has the following description:
Given  $D_\infty$, let $\mathcal{P}=\sum_{k} \delta_{p_k}$ be a Poisson point process
with intensity
$CD_\infty e^{-\sqrt{2}x}\mathrm{d}x$ and let
$\mathcal{D}^{(k)}=\sum_{n} \delta_{\Delta_n^{(k)}}$ be  iid copies of $\mathcal{D}_0$, then
\[
\sum_{u\in Z(t)} \delta_{X_u(t)- \left(\sqrt{2}t-\frac{3}{2\sqrt{2}}\log t \right)} \stackrel{\mathrm{d}}{\Longrightarrow} \sum_{k,n} \delta_{p_k + \Delta_n^{(k)}}
\]
in the sense of vague topology.
For the case of branching random walks,
see \cite{AE, HuShi, Madaule}.
For the case of $d$-dimensional branching Brownian motions, see \cite{BKLMZ}. For the case of super-Brownian motions, see \cite{RSZ, RYZ}.

In this paper, we consider (irreducible) multitype branching Brownian motions.
Let $S=\{1,...,d \}$ be the set of all types and $i\to \{p_{\mathbf{k}}(i): \mathbf{k}=(k_1,...,k_d)^T \in \mathbb{N}^d \}$ be the offspring distribution of type $i$ particles, here $\mathbb{N}=\{0,1,...\}$.
Let $a_i>0, i\in S$, be the branching rate of type $i$ particles.
A multitype branching Brownian motion can be defined as follows:  Initially, there is a particle of type $i$  at site $x$ and it moves according a standard Brownian motion.
After an exponential time with parameter $a_i$, it dies and
splits into $k_1$ offspring of type 1, $k_2$ offspring of type 2, \dots, $k_d$ offspring of type $d$ with probability $p_{\mathbf{k}}(i)$, where  $\mathbf{k}=(k_1,...,k_d)^T$.
The offspring evolve independently, each moves according to a standard Brownian motion and each type $j$ particle reproduces with law $\{p_{\mathbf{k}}(j):\ \mathbf{k}\in \mathbb{N}^d\}$  after an exponential distributed lifetime with parameter $a_j$.
This procedure goes on. We denote the law of this process by  $\P_{(x,i)}$. We use $\E_{(x,i)}$ to denote the expectation with respect to $\P_{(x,i)}$.
The multitype branching Brownian motion is related to the following system of  F-KPP equations:
\begin{align}\label{F-KPP}
	\mathbf{u}_t = \frac{1}{2}\mathbf{u}_{xx} + \Lambda \left(\mathbf{\psi}(\mathbf{u})- \mathbf{u}\right),
\end{align}
where
\begin{align}
	\mathbf{u}(t,x)&:= \left(u_1(t,x),..., u_d(t,x)\right)^T,\quad
	\Lambda := \mbox{diag}\left\{a_1,...,a_d \right\},\\
\mathbf{\psi}(\mathbf{u})&:= \left(\psi_1(\mathbf{u}),...,\psi_d(\mathbf{u}) \right)^T,
\quad
	\psi_i(\mathbf{u}):= \sum_{\mathbf{k}\in \N^d}p_{\mathbf{k}}(i)\prod_{j=1}^{d}
	u_j^{k_j}.
\end{align}
Let
\[
m_{i,j}:= \sum_{\mathbf{k}\in \N^d} p_{\mathbf{k}}(i)k_j < \infty,\quad i,j\in S,
\]
be the mean number of type $j$ offspring given birth by a type $i$ particle. Assume that the mean matrix $M= (m_{i,j})_{i,j\in S}$ is irreducible, i.e.,  there exists no permutation matrix $S$ such that $S^{-1}MS$ is block triangular.
We use $N_i(t)$ to denote the number of
type $i$ particles alive at time $t$.
Assume that $m_{i,j}(t)= \E_{(0,i)}\left(N_j(t)\right)<\infty$.
Then $M(t):=
\left(m_{i,j}(t)\right)_{i,j \in S}$ satisfies
(see the paragraph below \cite[(2)]{RY})
\begin{align}\label{Mean-Matrix}
	M(t)= \sum_{n=0}^\infty \frac{A^n}{n!} t^n,\quad \mbox{with}\ A:= \left(a_{i,j}\right)_{i,j\in S}\ \mbox{and}\ a_{i,j}=a_i\left(m_{i,j}-\delta_{i,j}\right).
\end{align}
For any $\mathbf{u}=(u_1,\cdots u_d)^T $ and $\mathbf{v}=(v_1,\cdots v_d)^T$, define $\langle\mathbf{u}, \mathbf{v} \rangle:=\sum^d_{i=1}u_iv_i.$
According to the Perron-Frobenius theorem,
the matrix $A$ admits a unique simple eigenvalue
$\lambda^*>0$, which is larger than the real part of any other eigenvalue, such that the associated left eigenvector $\mathbf{g}=(g_1,...,g_d)^T$  and right eigenvector $\mathbf{h}=(h_1,...,h_d)^T$
can be chosen to have all positive coordinates. We normalize $\mathbf{g}$  and $\mathbf{h}$
 so that $\langle \mathbf{g}, \mathbf{h} \rangle = \langle \mathbf{g}, \mathbf{1}\rangle = 1$,
 where $\mathbf{1}= (1,...,1)^T.$
We assume that $p_{\mathbf{0}}(i)=0$ for all $i\in S$, here $\mathbf{0}:=(0,...,0)^T$,
so the system survives with probability $1$.
We further assume that there exists $\alpha_0\in (0, 1]$ such that
\begin{align}\label{Second-Moment}
	 \sum_{\mathbf{k}\in \N^d} p_{\mathbf{k}}(i) k_j^{1+\alpha_0} < \infty,\quad \forall \ i,j \in S.
\end{align}
Define
\[
\mathbf{\varphi}(\mathbf{v}):=\mathbf{1} -\mathbf{\psi}(\mathbf{1} -\mathbf{v}).
\]
If  $\mathbf{u}$ is a solution of \eqref{F-KPP} and $\mathbf{v}:= \mathbf{1} -\mathbf{u}$, then $\mathbf{v}$ satisfies
\begin{align}\label{F-KPP2}
	\mathbf{v}_t = \frac{1}{2}\mathbf{v}_{xx} + \Lambda \left(\mathbf{1} -\mathbf{v} - \mathbf{\psi}(\mathbf{1}-\mathbf{v})\right)= \frac{1}{2}\mathbf{v}_{xx} + \Lambda \left(\mathbf{\varphi}(\mathbf{v})-\mathbf{v}\right).
\end{align}
Using the relationship between \eqref{F-KPP} and \eqref{F-KPP2}, by \cite[Lemma 5]{RY},  if $\mathbf{v}$ solves \eqref{F-KPP2}, then for all $i\in S, t>0, x\in \R$,
$v_i(t,x)$ has the following probabilistic representation
\begin{align}\label{Prob-Representation}
v_i(t,x)= 1- \E_{(x,i)}\left(\prod_{u \in Z(t)}\left(1- v_{I_u(t)}(0,X_u(t))\right) \right).
\end{align}
Here $Z(t)$ is the set of all the particles alive at time $t$ and for $u\in Z(t)$, $I_u(t)$ is the type of $u$ and $X_u(t)$ is the position of $u$.
In addition, we assume that
\begin{align}\label{Pure-Jump}
	m_{i,i} = 0,\quad \forall \ i\in S.
\end{align}
This assumption is not really necessary. When $m_{i,i}>0$ for some $i$, we only need to
 modify $a_i$ to $a_i+m_{i,i}.$

Define,
\begin{equation}\label{Additive-martingale}
	W_{\sqrt{2\lambda^*}}(s):= \sum_{u\in Z(s)} h_{I_u(s)}e^{-\sqrt{2\lambda^*}\left(X_u(s)+ \sqrt{2\lambda^*}s \right) },\quad s\geq 0,
 \end{equation}
 and
\begin{equation}\label{Derivative-martingale}	M_{\sqrt{2\lambda^*}}(s):= \sum_{u\in Z(s)} h_{I_u(s)}\left(X_u(s)+ \sqrt{2\lambda^*}s  \right)e^{-\sqrt{2\lambda^*}\left(X_u(s)+ \sqrt{2\lambda^*}s \right) },\quad s\geq 0.
\end{equation}
It is proved in
 \cite{RY} that  $\{W_{\sqrt{2\lambda^*}}(s), s\ge 0\}$ and $\{M_{\sqrt{2\lambda^*}}(s), s\ge 0\}$ are
 martingales, called the additive and derivative martingales of multitype branching Brownian motion, respectively.
 Note that the assumption  \eqref{Second-Moment} implies that $\sum_{\mathbf{k}\in \N^d} p_{\mathbf{k}}(i)k_j(\log_+k_j)^2 < \infty$. By \cite[Theorem 3]{RY},
\begin{equation}\label{addMlimit}
	\lim_{s\to\infty}W_{\sqrt{2\lambda^*}}(s)=0,\quad
\P_{(x,i)}\mbox{-a.s.}
\end{equation}
 According to \cite[Lemma 10, Theorem 5]{RY}, there is  a nonnegative and
 nondegenerate  random variable $ M_{\sqrt{2\lambda^*}}(\infty)$ such that
 \begin{equation}\label{derivMlimit}
 	\lim_{s\to\infty}M_{\sqrt{2\lambda^*}}(s)=M_{\sqrt{2\lambda^*}}(\infty),\quad \P_{(x,i)}
 	\mbox{-a.s.}
 \end{equation}

\subsection{Main results}

Our first main result is on the convergence of $\mathbf{v}$ to the traveling wave solution for a class of initial value conditions. Our second main result is about  the characterization of the extremal process of multitype branching Brownian motion.

For the initial value of $\mathbf{v}$, we assume that there exist $N_1 <N_2$ and $i_0 \in S$ such that
\begin{align}\label{Initial-Assumption}
v_i(0,x)\leq 1_{(-\infty, N_2)}(x),\quad \mbox{for all }i \in S\quad \mbox{and}\quad v_{i_0}(0,x)\geq 1_{(-\infty, N_1)}(x).
\end{align}
Let $m(t):= \sqrt{2\lambda^*}t - \frac{3}{2\sqrt{2\lambda^*}}\log t$ for $t>0$.
\begin{thrm}\label{thm1}
	Suppose that $\mathbf{v}$ solves \eqref{F-KPP2} with initial value satisfying \eqref{Initial-Assumption}, then it holds that for any $i\in S$ and $x\in \R$,
	\begin{align}
				\lim_{t\to\infty}(1- v_i \left(t, m(t)+x \right)) =
		\E_{(0,i)}\left(\exp\left\{ -  C_v(\infty) M_{\sqrt{2\lambda^*}}(\infty) e^{-\sqrt{2\lambda^*}x}\right\}  \right),
	\end{align}
where  $M_{\sqrt{2\lambda^*}}(\infty)$ is given in \eqref{derivMlimit} and $C_v(\infty)$ is defined by
	\begin{align}\label{Cv(infty)}
		C_v(\infty):= \lim_{r\to\infty} \sqrt{\frac{2}{\pi }} \int_{0}^\infty y e^{\sqrt{2\lambda^*}y}\ \left(\sum_{j=1}^d g_j
v_j(r,y +\sqrt{2\lambda^*}r)\right)  \mathrm{d} y\in (0,\infty).
	\end{align}
\end{thrm}

Let $v_i(0,x)=1_{(-\infty, 0)}(x)$, then by  \eqref{Prob-Representation} we have
\begin{align}
	v_i(t,x) =
	1- \E_{(x,i)}\left(\prod_{u \in Z(t)}\left(1- 1_{\{X_t(u)< 0 \}}\right) \right)= \P_{(x,i)}\left(\min_{u\in Z(t)} X_t(u)< 0 \right)= \P_{(0,i)}\left(M_t > x \right),
\end{align}
where $M_t:= \max_{u\in Z(t)} X_u(t)$ and we used the symmetry of Brownian motion in the last equality.
Using this, we get the following corollary of Theorem \ref{thm1}:
\begin{cor}\label{cor4}
	For any $i\in S$ and $x\in \R$,
	\begin{align}
	\lim_{t\to \infty}	\P_{(0,i)}\left(M_t \leq m(t)+ x \right)= \E_{(0,i)}\left(\exp\left\{ -  C_\infty M_{\sqrt{2\lambda^*}}(\infty) e^{-\sqrt{2\lambda^*}x}\right\}  \right),
	\end{align}
	where
	\begin{align}\label{C-Infinity}
		C_\infty= \lim_{r\to\infty} \sqrt{\frac{2}{\pi}} \int_{0}^\infty y e^{\sqrt{2\lambda^*}y}\ \left(\sum_{j=1}^d g_j \P_{(0,j)}\left(M_r> \sqrt{2\lambda^*}r +y\right)\right)
		 \mathrm{d} y.
	\end{align}
\end{cor}

For $j\in S$, define
 $$M_t^j:= \max_{u\in Z(t): I_u(t)=j} X_u(t).$$
Fix $i_1\in S$. Taking $v_{i_1}(0,x)=1_{(-\infty, 0)}(x)$,  $v_j(0,x)=0$ for $j\neq i_1$,
then by  \eqref{Prob-Representation} we have
\begin{align*}
	v_i(t,x) &=
 1- \E_{(x,i)}\left(\prod_{u \in Z(t), I_u(t)=i_1}\left(1- 1_{\{X_t(u)< 0 \}}\right) \right)\\
&= \P_{(x,i)}\left(\min_{u\in Z(t), I_u(t)=i_1} X_t(u)< 0 \right)= \P_{(0,i)}\left(M^{i_1}_t > x \right).
\end{align*}
Then we get the following corollary of Theorem \ref{thm1}:
\begin{cor}\label{cor3}
Fix $i_1\in S$.	For any $i\in S$ and $x\in \R$,
	\begin{align}
\lim_{t\to \infty}	\P_{(0,i)}\left(M_t^{i_1} \leq m(t)+ x \right)= \E_{(0,i)}\left(\exp\left\{ -  C_\infty^{(i_1)} M_{\sqrt{2\lambda^*}}(\infty) e^{-\sqrt{2\lambda^*}x}\right\}  \right),
	\end{align}
where
\begin{align}\label{C-Infinity-i-0}
C_\infty^{(i_1)}:= \lim_{r\to\infty} \sqrt{\frac{2}{\pi }} \int_{0}^\infty y e^{\sqrt{2\lambda^*}y}\ \left(\sum_{j=1}^d g_j \P_{(0,j)}\left(M_r^{i_1}> \sqrt{2\lambda^*}r +y\right)\right)
		\mathrm{d} y\in (0,\infty).
\end{align}
\end{cor}

\begin{thrm}\label{Conver-Gap-Process}
	Define
	\begin{align}
		\mathcal{D}_t:= \sum_{u\in Z(t)}\delta_{\left(X_u(t)-M_t, I_u(t) \right)},\quad t\ge0.
	\end{align}
Under $\mathbb{P}_{(0,i)}\left(\cdot \big| M_t>\sqrt{2\lambda^*}+z\right)$,
$\left(\mathcal{D}_t, M_t -\sqrt{2\lambda^*}-z\right)$ converges in distribution to
some $\left(\mathcal{D}, Y\right)$
as $t\to\infty$,
where $Y$ is an exponential random variable with parameter $\sqrt{2\lambda^*}$,
 $\mathcal{D}$ does not depend on  $i\in S$ and $z\in \R$,  and $\mathcal{D}$ and  $Y$ are independent.
\end{thrm}

Define
\[
\mathcal{E}_t:= \sum_{u\in Z(t)} \delta_{\left(X_u(t)-m(t), I_u(t) \right)}.
\]
Let $\mathcal{C}_c^+(\R \times S)$ be the set of all functions $\phi: \R \times S \to \R_+$ such that for any $j\in S$, $\phi(\cdot, j)$ is a non-negative
continuous function of compact support.

\begin{thrm}\label{thm2}
	Given $M_{\sqrt{2\lambda^*}}(\infty)$, let $\mathcal{P}\equiv \sum_{k\in \N} \delta_{p_k}$ be
	a  Poisson point process with intensity
	$C_\infty M_{\sqrt{2\lambda^*}}(\infty)\sqrt{2\lambda^*}e^{-\sqrt{2\lambda^*}x}\mathrm{d}x $, and let  $\{ \mathcal{D}^{(k)}:k\in \N\}$ be iid copies of
$\mathcal{D}$ defined in
Theorem \ref{Conver-Gap-Process}.
If $\mathcal{D}^{(k)}:= \sum_{n\in \N} \delta_{\left(\Delta_n^{(k)}, q_n^{(k)}\right)}$,
then for any $i\in S$, under $\mathbb{P}_{(0,i)}$, $\mathcal{E}_t$ converges in distribution to
	\begin{align}
		\mathcal{E}_\infty \stackrel{\mathrm{d}}{=} \sum_{k,n}
		\delta_{(p_k+\Delta_n^{(k)}, q_n^{(k)})}
\quad \mbox{as }t\to\infty.
	\end{align}
\end{thrm}

 For  $i_1\in S$, define
\[
 \mathcal{E}_t(i_1):= \sum_{u\in Z(t): I_u(t)= i_1}
\delta_{X_u(t)-m(t)}.
\]
As a consequence of Theorem \ref{thm2}, we have the following corollary:
\begin{cor}
$\mathcal{E}_t(i_1)$ converges
in distribution to
\begin{align}
 \mathcal{E}_\infty(i_1) \stackrel{\mathrm{d}}{=} \sum_{k}\sum_{n: q_n^{(k)}=i_1}\delta_{p_k+\Delta_n^{(k)}}
\quad \mbox{as }t\to\infty,
\end{align}
	where $p_k, \Delta_n^{(k)}$ and $q_n^{(k)}$ are defined in Theorem \ref{thm2}.
\end{cor}
\begin{remark}
The asymptotic behavior above for
irreducible multi-type branching Brownian motion is similar to the one obtained in  \cite{ABBS, ABK} for a single-type branching Brownian motion.
	Belloum and Mallein \cite{Belloum-Mallein} and Belloum \cite{Belloum} considered a 2-type reducible  branching Brownian motion and their results are quite different. In their model, particles of type 1 move as a Brownian motion with diffusion coefficient $\sigma^2$,
	reproduce with branching rate $\beta + \alpha$ and offspring distribution $\{p_{\mathbf{k}}(1)\}$ satisfying $p_{(2,0)}(1)= \frac{\beta}{\beta +\alpha},  p_{(1,1)}(1)= \frac{\alpha}{\alpha +\beta}$.
	Particles of type 2 evolve as a standard branching Brownian motion
	with branching rate 1 and binary branching.
In the spacial case when $\sigma^2=1$ and $\beta =1$,  $\lambda^*=1$, their
results show that the corresponding front $m(t)$ is $\sqrt{2}t -\frac{1}{2\sqrt{2}}\log t$,
which is quite different from the irreducible case where the corresponding front for $\lambda^*=1$
is $\sqrt{2}t -\frac{3}{2\sqrt{2}}\log t$.
\end{remark}

In the remainder of this paper, for a set $E$, the notation
\[
f(x)\lesssim g(x),\quad x\in E
\]
means
that there exists some constant $C$ independent of $x\in E$ such that $f(x)\leq C g(x)$ holds for all $x\in E$.
Also, the notation $f\asymp g, x\in E,$ means $f\lesssim g , x\in E$ and $g\lesssim f, x\in E$.

\section{Many-to-one formula and spine decomposition}\label{Many-to-one-spine}

Let $\mathbf{N}_t:= (N_1(t),...,N_d(t))^T$ and
let $\{\mathcal{F}_t\}$ be the natural filtration of  the multitype branching Brownian motion.
By \cite[Proposition 2]{Athreya},  under $\mathbb{P}_{(x,i)}$, $e^{-\lambda^* t}\langle \mathbf{N}_t, \mathbf{h}\rangle$ is a mean $h_i$ positive martingale with respect to $\{\mathcal{F}_t\}$.
Define $\widehat{\P}_{(x,i)}$ by
\begin{align}\label{Change-of-measure}
	\frac{\mathrm{d} \widehat{\P}_{(x,i)}}{\mathrm{d} \P_{(x,i)}}\bigg|_{\mathcal{F}_t}:= \frac{e^{-\lambda^* t }\langle \mathbf{N}_t, \mathbf{h}\rangle}{h_i}.
\end{align}
According to \cite[p. 224]{RY},
the multitype branching Brownian motion under $\widehat{\P}_{(x,i)}$ has the following \emph{spine decomposition}:

(i)
Initially there is a marked particle $\xi$, called the spine, of type $i$ at site $x$.

(ii)
After an exponential time $\zeta_\xi$ with parameter
$a_i+\lambda^*$,
 this marked particle dies and  produces $A_1:=A_1(i)$ offspring of type 1, \dots,
$A_d:=A_d(i)$ offspring of type $d$ with probability
$\widehat{p}_{\mathbf{A}}(i)
:=\frac{p_{\mathbf{A}}(i)\langle \mathbf{A}, \mathbf{h}\rangle}{(1+a_{i}^{-1}\lambda^*)h_i}$, where $\mathbf{A}=(A_1, \dots A_d)^T$.
Randomly choose one of these  $\langle \mathbf{A}, \mathbf{1}\rangle$ particle to continue as the spine, with each type $j$ particle being chosen with probability $h_j/ \langle \mathbf{A}, \mathbf{h}\rangle$.

(iii)
The $\langle \mathbf{A}, \mathbf{1}\rangle$ offspring particles evolve independently, with the
marked (spine) particle repeating step (ii) with law
$\widehat{\mathbb{P}}_{(X_\xi(\zeta_\xi), I_\xi(\zeta_\xi))}$
and each unmarked particle of type $j$, $j\in S$,
evolving as a multitype branching Brownian motion with law $\mathbb{P}_{(X_\xi(\zeta_\xi), j)}$.
The process then goes on.

If we only consider the spine process $(X_\xi(t), I_\xi(t))$, then, under $\widehat{\P}_{(x,i)}$,
$X_\xi$ is a standard Brownian motion starting from $x$, $I_\xi$ is an $S$-valued Markov chain with generator
\[
\widehat{G}= (\widehat{g}_{i,j})_{S\times S}\quad \mbox{with}\quad \widehat{g}_{i,j}:= \left(a_i +\lambda^*\right)\left(\frac{m_{i,j}h_j}{\left(1+\lambda^*/a_i\right)h_i} -\delta_{i,j}\right)= \frac{a_im_{i,j}h_j}{h_i}- (a_i +\lambda^*)\delta_{i,j}
\]
and $X_\xi$ is independent of $I_\xi(t)$.
 According to \cite[(12)]{RY}, we have
\begin{align}\label{Many-to-one-2}
	\widehat{\P}_{(x,i)}\left(\xi_t = u \big| \mathcal{F}_t\right)=\frac{h_{I_u(t)}}{\langle \mathbf{N}_t, \mathbf{h} \rangle}
	1_{\{u\in Z(t) \}}.
\end{align}
Using \eqref{Many-to-one-2}, we give a stronger version of the many-to-one formula in  \cite[Proposition 1]{RY}.
For the case of branching Brownian motions,
one can refer to \cite[Proposition 4.1]{MP} (In \cite{MP}, there is  also a change-of-measure for the spinal movement).
\begin{prop}\label{Many-to-one}
	For any $t>0$ and $u\in Z(t)$, let $H(u, t)$ be a  non-negative $\mathcal{F}_t$-measurable  random variable. Then
	\[
	\E_{(x,i)} \left(\sum_{u\in Z(t)} H(u,t) \right) =e^{\lambda^* t}\widehat{\E}_{(x,i)}\left(H(\xi_t, t)\frac{h_i}{h_{I_\xi(t)}}\right).
	\]
\end{prop}
\textbf{Proof: }  By \eqref{Change-of-measure} and \eqref{Many-to-one-2},
\begin{align*}
	&\E_{(x,i)} \left(\sum_{u\in Z(t)} H(u,t) \right)  = e^{\lambda^*t}\widehat{\E}_{(x,i)}\left( \frac{h_i}{\langle \mathbf{N}_t, \mathbf{h}\rangle}\sum_{u\in Z(t)} H(u,t)\right)\\
	& = e^{\lambda^*t}\widehat{\E}_{(x,i)}\left( \sum_{u\in Z(t)} \frac{h_{I_u(t)}}{\langle \mathbf{N}_t, \mathbf{h}\rangle}H(u,t) \frac{h_i}{h_{I_u(t)}}\right) =  e^{\lambda^*t}\widehat{\E}_{(x,i)}\left( \sum_{u\in Z(t)} 	\widehat{\P}_{(x,i)}\left(\xi_t = u \big| \mathcal{F}_t\right) H(u,t) \frac{h_i}{h_{I_u(t)}}\right) \\
	& = e^{\lambda^*t}\widehat{\E}_{(x,i)}\left( H(\xi_t,t)\frac{h_i}{h_{I_\xi(t)}}\sum_{u\in Z(t)} 	1_{\{\xi_t =u\}} \right) = e^{\lambda^*t}\widehat{\E}_{(x,i)}\left( H(\xi_t,t)\frac{h_i}{h_{I_\xi(t)}}\right).
\end{align*}

\hfill$\Box$

It is easy to deduce from  $A\mathbf{h}=\lambda^* \mathbf{h}$ and $\mathbf{g}^T A =\lambda^* \mathbf{g}^T$ that
\[
\sum_{j=1}^d m_{i,j} h_j =\frac{a_i +\lambda^*}{a_i}h_i,\quad \sum_{j=1}^d m_{j,i} a_j g_j = (a_i+ \lambda^*) g_i,
\quad i=1,2,...,d.
\]
Let $\mu_j:= h_j g_j$ and $\mathbf{\mu}= (\mu_1,...,\mu_d)^T$,
then we see that $\mu_j$ solves the equation
\[
\mu_j(a_j +\lambda^*)=\sum_{i=1}^d \mu_i \frac{a_i}{h_i}m_{i,j}h_j,\quad j=1,...,d,
\]
and $\langle \mathbf{1},\mathbf{\mu} \rangle  =1$,
which implies that $\mathbf{\mu}$ is an invariant measure for $I_\xi(t)$ under $\widehat{\P}_{(x,i)}$.

\section{Non-local Feynman-Kac formula}

Throughout this paper, $\left(X_t, t\ge 0; \mathbf{P}_x\right)$ is a standard Brownian motion starting from $x$.
Feynman-Kac formula plays an important role in the probabilistic treatment of the F-KPP equation \eqref{F-KPP-1d}.
The classical Feynman-Kac formula says that
a solution of the  linear equation
$$u_t=\frac{1}{2}u_{xx}+k(t,x)u$$
can be given by
\begin{equation}\label{F-K}
	u(t,x)=
	\mathbf{E}_x
	\left(e^{\int^t_0k(t-s, X_s)ds}u(0, X_t)\right).
\end{equation}
  If $u$ is a solution to  equation \eqref{F-KPP-1d}, the \eqref{F-K} holds with
$k(s, y)=\frac{(u^2-u)}{u}(s,y)=u(s,y)-1$.
For our multitype branching Brownian motion,
we will give similar representation for  a solution $\mathbf{v}$ of \eqref{F-KPP2} using  a non-local Feynman-Kac formula.

First note that
\eqref{F-KPP2} is equivalent to
\[
v_i(t,x)=\mathbf{E}_x \left(v_i(0, X_t)\right)+\mathbf{E}_x \left(\int_0^t a_i\left(\varphi_i\left(\mathbf{v}(t-s, X_s)\right)-v_i(t-s,X_s)  \right) \mathrm{d}s \right),\quad  i\in S.
\]
Let $n_i:= \sum_{j=1}^d m_{i,j}$ and $p_{i,j}:= m_{i,j}/n_i$.
Since $p_{\mathbf{0}}(i)=0$ and $m_{i,i}=0$ by assumption, we have $n_i\geq 1$ and $p_{i,i}=0$ for all $i\in S$.
Rewrite $A$ given in \eqref{Mean-Matrix} as
\begin{equation}\label{Mean-Matrix'}
A = \mbox{diag}\left\{ a_1 (n_1-1), a_2(n_2-1), \dots, a_d(n_d-1) \right\} + A^*,
\end{equation}
where $A^*:= (a_{i,j}^*)$ and $a_{i,j}^*= a_in_i(p_{i,j}-\delta_{i,j})$.
 Define
\begin{equation}\label{def-mathcal-A}
 \mathcal{A}:= \mbox{diag}\left\{\frac{1}{2}\frac{\mathrm{d}^2}{\mathrm{d}x^2},..., \frac{1}{2}\frac{\mathrm{d}^2}{\mathrm{d}x^2} \right\} +A^*.
\end{equation}
Let $I_t$ be a  continuous time Markov chain,  independent of $X$, with generator $A^*$.
We use $\mathbf{P}_{(x, i)}$ to denote the law of $(X_t, I_t)$ and use $\mathbf{E}_{(x, i)}$ to denote the corresponding expectation.
Then \eqref{F-KPP2} is equivalent to
\[
\mathbf{v}_t = \mathcal{A}\mathbf{v} + \Lambda \left(\mathbf{\varphi}(\mathbf{v})- \mathbf{v}\right)- A^*\mathbf{v},
\]
which in turn is equivalent to
\begin{align}\label{Integral-eq}
	v_i(t,x)&=
	\mathbf{E}_{(x,i)}\left(v_{I_t}\left(0, X_t\right)\right) +
	\mathbf{E}_{(x,i)}\left(\int_0^t a_{I_s}\left(\varphi_{I_s}\left(\mathbf{v}(t-s, X_s)\right)-v_{I_s}(t-s,X_s) \right)\mathrm{d}s\right)\nonumber\\
	&\quad -\mathbf{E}_{(x,i)}\left(\int_0^t \sum_{j=1}^d a_{I_s,j}^* v_j(t-s, X_s) \mathrm{d}s\right)\nonumber\\
	& =
		\mathbf{E}_{(x,i)}\left(v_{I_t}\left(0, X_t\right)\right)
	+ \mathbf{E}_{(x,i)}\left(\int_0^t a_{I_s}(n_{I_s}-1)v_{I_s}(t-s,X_s)\mathrm{d}s\right)\nonumber\\\
	&\quad + \mathbf{E}_{(x,i)} \left(\int_0^ta_{I_s}\left(\varphi_{I_s}\left(\mathbf{v}(t-s, X_s)\right) - n_{I_s}\sum_{j=1}^d p_{I_s,j}v_j(t-s,X_s) \right)\mathrm{d}s \right).
\end{align}

We will simplify the formula  above using the non-local Feynman-Kac formula introduced below.
Define a Feynman-Kac semigroup $P_t^{\mathcal{A} + \Lambda(N-I)}$ by
\[
P_t^{\mathcal{A} + \Lambda(N-I)}f (x,i) := \mathbf{E}_{(x,i)}\left(f(X_t,I_t)\exp\left\{\int_0^t a_{I_s}\left(n_{I_s}-1\right)\mathrm{d} s \right\} \right),
\]
 then (see, for instance, \cite[Lemma 2.1]{HHKM}),
since $m_{i,i}=0$ for all $i$, $P_t^{\mathcal{A} + \Lambda(N-I)}$ is the mean semigroup of a purely non-local branching Markov process with spatial motion $(X, I)$, branching rate function $\beta(x, i)=a_i$ and non-local probability distribution
$F\left((x,i), \cdot\right)$, on the space $\mathcal{M}(\mathbb{R}\times S)$ of finite measures on $\mathbb{R}\times S$, defined for all $(x,i)\in \mathbb{R}\times S$ by
\[
F\left((x,i) ,  \left\{ \delta_x \times \left(k_1 \delta_1+...+ k_d \delta_d\right)\right\} \right)= p_{\mathbf{k}}(i), \quad \mathbf{k}\in \mathbb{N}^d.
\]
Put
$$\mathbf{z} (t,x)= (z_1(t,x),\dots, z_d(t,x))^T\quad\mbox { with } z_i(t,x):=P_t^{\mathcal{A} + \Lambda(N-I)}f (x,i).$$
Then $\mathbf{z}$ solves the linear equation
\begin{align}
	\mathbf{z}_t= \mathcal{A} \mathbf{z} +\mbox{diag}\left\{ a_1 (n_1-1), a_2(n_2-1), \dots, a_d(n_d-1) \right\} \mathbf{z},
	\quad x\in\mathbb{R}, i\in S.
\end{align}
By \eqref{def-mathcal-A} and \eqref{Mean-Matrix'},
 we see that $\mathbf{z}$ solves the equation
\begin{align}\label{Mean-semigroup}
\mathbf{z}_t = \mbox{diag}\left\{\frac{1}{2}\frac{\mathrm{d}^2}{\mathrm{d}x^2},..., \frac{1}{2}\frac{\mathrm{d}^2}{\mathrm{d}x^2} \right\}  \mathbf{z} + A \mathbf{z}.
\end{align}

Let $H(s)=s$ and
\[
J((x,k), \mathrm{d}(y,\ell)):= \delta(x-y) a_kn_k\left(p_{k,\ell}-\delta_{k,\ell} \right) 1_{\{k\neq \ell\}} \mathrm{d}y\mathrm{d}\ell.
\]
Define
\[
    D_J:= \left\{t\big| \ (X_{t-}, I_{t-})\neq (X_t, I_t) \right\}= \left\{t \big| \ I_{t-}\neq I_t\right\}.
\]
It is easy to check that for any non-negative Borel function $f$ on $(\mathbb{R}\times S)^2 $ vanishing on the diagonal and any $x \in \mathbb{R}, i \in S$,
	\begin{align*}
		&\mathbf{E}_{(x,i)}\left(\sum_{s\in D_J, s \leq t} f\left((X_{s-}, I_{s-}), (X_{s}, I_s)\right)\right)\\
		&=\mathbf{E}_{(x,i)}\left(\int_{0}^{t} \int_{\mathbb{R}\times S} f\left((X_{s}, I_s), (y,\ell)\right) J\left((X_s, I_s), \mathrm{d}(y,\ell)\right) \mathrm{d} s\right),
	\end{align*}
and thus 	$(J, H)$ is a L\'{e}vy system for $(X,I)$.
Note that the non-linear term of \eqref{Integral-eq} is equal to
\begin{align*}
	&a_{I_s}\left(\varphi_{I_s}\left(\mathbf{v}(t-s, X_s)\right)-n_{I_s}\sum_{j=1}^d p_{I_s, j}v_j (t-s, X_s) \right) \\
	&= \left(\frac{\varphi_{I_s}\left(\mathbf{v}(t-s, X_s)\right) }{\sum_{j=1}^d m_{I_s, j}v_j (t-s, X_s)}-1\right) a_{I_s}n_{I_s} \sum_{j=1}^d p_{I_s, j}v_j (t-s, X_s).
\end{align*}
Applying \cite[Lemma A.1]{CRS}, similar to \cite[(4.8)]{CRS}, \eqref{Integral-eq} can be written as
\begin{align}\label{Feynman-Kac}
	v_i(t,x)=\mathbf{E}_{(x,i)}\left(\exp\left\{\sum_{s\in D_J, s\leq t} \log \left(\frac{ \varphi_{I_s}\left(\mathbf{v}(t-s, X_s)\right)}{\sum_{j=1}^d m_{I_s, j}v_j (t-s, X_s)}\right) + \int_0^t a_{I_s}(n_{I_s}-1)\mathrm{d}s\right\}v_{I_t}(0,X_t)\right).
\end{align}

For any bounded non-negative function $f$ and $\theta>0$, by \eqref{Prob-Representation},
\[
u_i(t,x; \theta):= \mathbb{E}_{(x,i)} \left(\prod_{u \in Z(t)} e^{-\theta f(X_u(t), I_u(t))}\right)
\]
solves equation \eqref{F-KPP}. By taking derivative with respect to $\theta$ and letting $\theta \downarrow 0$, it is easy to see that
\[
z_i(t,x)= \E_{(x,i)}\left(\sum_{u\in Z(t)} f(X_u(t), I_u(t))\right)
\]
also solves equation \eqref{Mean-semigroup}.
Therefore, $P_t^{\mathcal{A} + \Lambda(N-I)}$ is also the mean-semigroup of the multitype branching Brownian motion, i.e., for every bounded measurable function $f$,
\begin{align}
	P_t^{\mathcal{A} + \Lambda(N-I)}f (x,i) = \E_{(x,i)}\left(\sum_{u\in Z(t)} f(X_u(t), I_u(t))\right).
\end{align}
It follows from Proposition \ref{Many-to-one} that for $\phi(x,i)=  h_i,$
\begin{equation}\label{e:hasevector}
P_t^{\mathcal{A} + \Lambda(N-I)}\phi (x,i) =e^{\lambda^* t}\phi(x,i).
\end{equation}
Using the definition of $P_t^{\mathcal{A} + \Lambda(N-I)}$ and \eqref{e:hasevector}, we can easily see that
$$
e^{-\lambda^* t}\exp\left\{\int_0^t a_{I_s}\left(n_{I_s}-1\right)\mathrm{d} s \right\}
\frac{h_{I_t}}{h_{I_0}}
$$
is a non-negative martingale of mean 1 under $\mathbf{P}_{(x,i)}$. Now we define
\begin{equation}\label{def-P^h}
\frac{\mathrm{d}\mathbf{P}_{(x,i)}^h}{\mathrm{d}\mathbf{P}_{(x,i)}}\bigg|_{\sigma(X_s, I_s, s\leq t)}:= e^{-\lambda^* t}\exp\left\{\int_0^t a_{I_s}\left(n_{I_s}-1\right)\mathrm{d} s \right\}
\frac{h_{I_t}}{h_{I_0}}.
\end{equation}
Then by the definition of $P_t^{\mathcal{A} + \Lambda(N-I)}$,  we have
\[
P_t^{\mathcal{A} + \Lambda(N-I)}f (x,i) = e^{\lambda^* t}\mathbf{E}_{(x,i)}^h\left(f(X_t,I_t)\frac{h_i}{h_{I_t}} \right).
\]
Combining this with Proposition \ref{Many-to-one},  we get that
$$\left((X,I), \mathbf{P}_{(x,i)}^h\right)\stackrel{\mathrm{d}}{=} \left((X_\xi, I_\xi), \widehat{\P}_{(x,i)}\right).$$
Now \eqref{Feynman-Kac} can be rewritten as
\begin{align}\label{Feynman-Kac-2}
	v_i(t,x)=e^{\lambda^* t}h_i\mathbf{E}_{(x,i)}^h\left(\exp\left\{\sum_{s\in D_J, s\leq t} \log \left(\frac{ \varphi_{I_s}\left(\mathbf{v}(t-s, X_s)\right)}{\sum_{j=1}^d m_{I_s, j}v_j (t-s, X_s)}\right) \right\}\frac{v_{I_t}(0,X_t)}{h_{I_t}}\right).
\end{align}
Define for
$0\le r<t\leq t_1$,
\begin{align}\label{Non-Linear-R}
	R_{t_1}\left((r,t]; v \right):= \exp\left\{\sum_{s\in D_J, r< s\leq t} \log \left(\frac{ \varphi_{I_s}\left(\mathbf{v}(t_1-s, X_s)\right)}{\sum_{j=1}^d m_{I_s, j}v_j (t_1-s, X_s)}\right) \right\}
\end{align}
and
\begin{align}\label{Non-Linear-R-2}
R\left((r,t]; v \right):= R_{t}\left((r,t]; v \right), \quad R\left(t; v \right):= R\left((0,t]; v \right).
\end{align}
Now for $0< r<t$, by
the Markov property, we get from \eqref{Feynman-Kac-2} that
\begin{align}\label{Feynman-Kac-3}
	v_i(t,x)& = e^{\lambda^* (t-r)}h_i\mathbf{E}_{(x,i)}^h\left(
	R_t(t-r;v)
	\frac{v_{I_{t-r}}(r,X_{t-r})}{h_{I_{t-r}}}\right)\nonumber \\
	& =e^{\lambda^* (t-r)}h_i\int_\R  \frac{e^{-\frac{(x-y)^2}{2(t-r)}}}{\sqrt{2\pi (t-r)}}  \mathbf{E}_{(x,i)}^h\left(
	R_t(t-r;v)
	\frac{v_{I_{t-r}}(r,y)}{h_{I_{t-r}}}\bigg| X_{t-r}=y \right)\mathrm{d} y.
\end{align}
The above representation \eqref{Feynman-Kac-3} of $v_i$ will play an important role in this paper.
\bigskip

For any $i\in S$, by Bernoulli's inequality,
\[
\varphi_i (\mathbf{v})= 1-\psi_i(\mathbf{1}-\mathbf{v}) = 1- \sum_{\mathbf{k}\in \N^d} p_{\mathbf{k}}(i)\prod_{j =1}^d (1-v_j)^{k_j}\leq 1-\sum_{\mathbf{k}\in \N^d} p_{\mathbf{k}}(i)(1-\sum_{j=1}^d k_j v_j)= \sum_{j=1}^d m_{i,j}v_j,
\]
Thus, $\mathbf{P}_{(x,i)}^h$-a.s., for any $0\leq r<t$ and $t_1\geq t$, $R_{t_1}\left((r,t]; v \right) \leq 0.$

The assumption \eqref{Second-Moment} implies the following estimate on $\varphi_i (\mathbf{v})$:

\begin{lemma}
It holds uniformly for all $i\in S$ that
\begin{align}\label{Assumption1}
	\frac{\varphi_i (\mathbf{v})}{\sum_{j=1}^d m_{i,j}v_j} = 1- O\left(\Vert \mathbf{v} \Vert^{\alpha_0} \right),\quad \mbox{for all } \mathbf{v}\in [0,1]^d.
\end{align}
\end{lemma}
\textbf{Proof: }
For any  $\mathbf{v}\in [0,1]^d$, let
$F(r):= \varphi_i(r\mathbf{v} )$ for $r\in [0,1]$, then there exists $\theta \in [0,1]$ such that
\[
\varphi_i(\mathbf{v})  =  F(1)- F(0)= F'(\theta )= \triangledown \varphi_i(\theta\mathbf{v})\cdot \mathbf{v} .
\]
Let $J(i):= \left\{ j\in S: m_{i,j}>0 \right\}$, then for any $j\notin J(i)$, $p_{\mathbf{k}}(i)=0$ for any $\mathbf{k}\in \N^d$ with $k_j>0$. Therefore, by the trivial
inequalities
$$
1-\prod_{j=1}^d x_j\leq \sum_{j=1}^d (1-x_j), \quad x_j\in [0, 1], j=1, \dots, d
$$
and
$$
1-(1-x)^{k}\leq k^{\alpha_0}x^{\alpha_0}, \quad k\ge 1, x\in [0, 1],
$$
we have
\begin{align*}
	&\left|\varphi_i(\theta \mathbf{v}) - \triangledown \varphi_i(0)\cdot \mathbf{v} \right|\leq \sum_{\ell=1}^d v_\ell \left|\sum_{\mathbf{k}\in \N^d} k_\ell p_{\mathbf{k}}(i)\prod_{j =1}^d (1-\theta v_j)^{k_j -\delta_{j,\ell}} - m_{i,\ell}\right|\\
	& =\sum_{\ell\in J(i)} v_\ell \left(\sum_{\mathbf{k}\in \N^d} k_\ell p_{\mathbf{k}}(i) - \sum_{\mathbf{k}\in \N^d} k_\ell p_{\mathbf{k}}(i)\prod_{j =1}^d (1-\theta v_j)^{k_j -\delta_{j,\ell}} \right)\\
	& \leq \sum_{\ell\in J(i)} v_\ell \sum_{\mathbf{k}\in \N^d} k_\ell p_{\mathbf{k}}(i) \left(1- \prod_{j =1}^d (1-\theta v_j)^{k_j } \right)\leq \sum_{\ell\in J(i)} v_\ell \sum_{\mathbf{k}\in \N^d} k_\ell p_{\mathbf{k}}(i) \sum_{j=1}^d \left(1- (1-\theta v_j)^{k_j } \right)\\
	& \leq \sum_{\ell\in J(i)} v_\ell \sum_{\mathbf{k}\in \N^d} k_\ell p_{\mathbf{k}}(i) \sum_{j=1}^d k_j^{\alpha_0} v_j^{\alpha_0}\leq \frac{1}{\min_{\ell \in J(i)}m_{i,\ell} } \sum_{\ell\in J(i)} m_{i,\ell} v_\ell \sum_{\mathbf{k}\in \N^d} k_\ell p_{\mathbf{k}}(i) \sum_{j=1}^d k_j^{\alpha_0} \cdot \Vert \mathbf{v}\Vert^{\alpha_0}\\
	& =: \Gamma(i) \left(\triangledown \varphi_i(0)\cdot \mathbf{v}\right) \Vert \mathbf{v} \Vert^{\alpha_0},
\end{align*}
where we used \eqref{Second-Moment}  at the end of the display above. Thus \eqref{Assumption1} is valid.
\hfill$\Box$

\bigskip

In the remainder of this paper,
when we consider the spine process $(X_\xi, I_\xi)$  only under $\widehat{\P}_{(x,i)}$,
we sometimes use $\left((X, I), \mathbf{P}_{(x,i)}^h\right)$ to denote the law of the spine process for simplicity.

 \section{Estimates in the case of Heaviside initial conditions}

In this section, we consider two kinds of initial conditions. The first kind is
\begin{align}\label{Heav1}
    v_i(0,x)= 1_{(-\infty,0)}(x),\quad \mbox{for all }i \in S.
\end{align}
Fix $i' \in S$.
The second kind of initial condition is
\begin{align}\label{Heav2}
 v_{i'}(0,x)=1_{(-\infty, 0)}(x),\quad v_j(0,x)=0, \quad\mbox{ for }i\neq i'.
\end{align}
Note that if $\mathbf{v}$ solves  \eqref{F-KPP2} with  initial condition \eqref{Heav1}, then
$$v_i(t,x) = \P_{(x,i)}\left(\min_{u\in Z(t)} X_t(u)< 0 \right)= \P_{(0,i)}\left(M_t > x \right);$$
and that if $\mathbf{v}$ solves \eqref{F-KPP2} with initial condition \eqref{Heav2}, then
$$v_i(t,x)=\P_{(x,i)}\left(\min_{u\in Z(t), I_u(t)=i'} X_t(u)< 0 \right)= \P_{(0,i)}\left(M^{i'}_t > x \right).$$

The purpose of this section is to get estimates on solutions
 $\mathbf{v}(t,x)$ of  \eqref{F-KPP2} with Heaviside initial conditions \eqref{Heav1} or \eqref{Heav2}, and with $x=m(t)+y$,
that is to say, we want to  get some upper and lower bounds for
$\P_{(0,i)} \left(M_t \geq m(t) +y\right)$
    and $\P_{(0,i)} \left(M_t^{i'} \geq m(t) +y\right)$ with $y>0$.
See Proposition \ref{Upper-Bound-M-T}  below for the upper bound and Proposition \ref{Lower-Bound-M-T} for the lower bound.
Then we use Propositions \ref{Upper-Bound-M-T} and \ref{Lower-Bound-M-T} to
prove that for any $i\in S$, $\left(M_t - m(t),\ t\geq1; \P_{(0,i)}\right) $
	is tight, and that,  for any $i,i' \in S$, $\left(M_t^{i'} - m(t),\ t\geq1; \P_{(0,i)}\right) $ is tight.

We first prove an estimate on the path of Brownian motion.

\begin{lemma}\label{lemma1}
Let $K>0$, $\alpha < 1/2$ and $t\ge 1$. For any function $f$ satisfying
	\[
	\sup_{s\leq t}\left(\frac{|f(s)|}{s^\alpha} + \frac{|f(t)-f(s)|}{(t-s)^\alpha}\right)< K,
	\]
	there exists a constant $\Gamma_1$
	depending only on $K$ and $\alpha$ such that
	\begin{align}
		\mathbf{P}_0 \left(B_s\geq -y +f(s), s\leq t,\ B_t +y -f(t) \in [z,z+1]\right)\leq \Gamma_1 \frac{(y\land \sqrt{t})(z\land \sqrt{t})}{ t^{3/2}},\quad
		y,z \geq 1,
	\end{align}
where $\left(B_t, t\geq 0; \mathbf{P}_x\right)$ is a standard Brownian motion starting from $x$.
\end{lemma}
\textbf{Proof: }
Using $\left(\inf_{s\leq t} B_s,\ \mathbf{P}_0\right)\stackrel{\mathrm{d}}{=} \left(-|B_t|,\ \mathbf{P}_0 \right)$, we can easily get
\begin{align}\label{Ballot-Theorem1}
	\mathbf{P}_0 \left( B_s \geq -y,\ s\leq t \right)  = \mathbf{P}_0\left(|B_t| \leq y\right)\lesssim \frac{y\land \sqrt{t}}{\sqrt{t}},\quad y,t >0.
\end{align}
Next, we prove that, for any $K_1\in \R$ and $\alpha \in (0, 1/2)$, it holds that
\begin{align}\label{Ballot-Theorem2}
	\mathbf{P}_0 \left( B_s \geq -y + K_1 s^\alpha,\ s\leq t \right)\lesssim \frac{y\land \sqrt{t}}{\sqrt{t}},\quad y,t \geq 1.
\end{align}
If $K_1>0$, then by \eqref{Ballot-Theorem1},
\begin{align}
	\mathbf{P}_0 \left( B_s \geq -y + K_1 s^\alpha,\ s\leq t \right) \leq \mathbf{P}_0 \left( B_s \geq -y ,\ s\leq t \right) \lesssim \frac{y\land \sqrt{t}}{\sqrt{t}},\quad y,t \geq 1.
\end{align}
If $K_1<0$, then by \cite[Lemma 3.6]{BM2}, when $y\leq \sqrt{t}$, it holds that
\begin{align}
	\mathbf{P}_0 \left( B_s \geq -y + K_1 s^\alpha,\ s\leq t \right) \leq \mathbf{P}_0 \left( B_j \geq -y - |K_1| j^\alpha ,\ j= 1,2,... [t] \right) \lesssim \frac{y}{\sqrt{[t]}} \lesssim \frac{y}{\sqrt{t}} , \quad t,y\geq 1.
\end{align}
When $y\geq \sqrt{t}$, we use the trivial upper bound $1$.

Now we prove the desired result. When $t\leq 3$,
we use the trivial upper-bound $1$.
When $t>3$, by the Markov property at time $t/3$ and \eqref{Ballot-Theorem2},
\begin{align}
	&\mathbf{P}_0 \left(B_s\geq -y +f(s), s\leq t,\ B_t +y -f(t) \in [z,z+1]\right) \\
	&\leq \mathbf{P}_0\left(B_s\geq -y -K s^\alpha, s\leq\frac{t}{3} \right)\cdot \sup_{x\in \R} \mathbf{P}_x \left(B_s\geq -y +f\left(s+ \frac{t}{3}\right), s\leq \frac{2t}{3},\ B_{2t/3} +y -f(t) \in [z,z+1]\right)\\
	& \lesssim \frac{y\land \sqrt{t}}{\sqrt{t}}\cdot \sup_{x\in \R} \mathbf{P}_x \left(B_s\geq -y +f\left(s+ \frac{t}{3}\right), s\leq \frac{2t}{3},\ B_{2t/3} +y -f(t) \in [z,z+1]\right) ,\quad
	y,z \geq 1.
\end{align}
For any $x\in \R$,
\begin{align}
	&\mathbf{P}_x \left(B_s\geq -y +f\left(s+ \frac{t}{3}\right), s\leq \frac{2t}{3},\ B_{2t/3} +y -f(t) \in [z,z+1]\right) \\
	& \leq \mathbf{P}_x \left(B_s- B_{2t/3}\geq -(z+1)+f\left(s+ \frac{t}{3}\right)-f(t),\frac{t}{3}\leq  s\leq \frac{2t}{3},\ B_0 - B_{2t/3} \in [h-1,h]\right) \\
	& = \mathbf{P}_0 \left(\widetilde{B}_s \geq -(z+1)+f\left(t-s\right)-f(t),  s\leq \frac{t}{3},\ \widetilde{B}_{2t/3} \in [h-1,h]\right)
\end{align}
with $h=x+y-z-f(t)$ and $\widetilde{B}_s:= B_{2t/3-s} - B_{2t/3}$ being still a Brownian motion
starting from $0$. By the Markov property of $\widetilde{B}$ at time $t/3$ and \eqref{Ballot-Theorem2}, we get that
\begin{align}
	&\mathbf{P}_0 \left(\widetilde{B}_s \geq -(z+1)+f\left(t-s\right)-f(t),  s\leq \frac{t}{3},\ \widetilde{B}_{2t/3} \in [h-1,h]\right)\\
	& \leq \mathbf{P}_0 \left(\widetilde{B}_s \geq -(z+1)-K s^\alpha,  s\leq \frac{t}{3}\right)\cdot \sup_{x'\in \R} \mathbf{P}_{x'}\left(\widetilde{B}_{t/3}\in [h-1, h] \right)\\
	& \lesssim \frac{(z+1)\land \sqrt{t}}{\sqrt{t}}\cdot \sup_{x'\in \R} \int_{h-1}^h \frac{e^{(y'-x')^2 /(t/3)}}{\sqrt{2\pi t/3}}\mathrm{d}y'\\
	& \leq \frac{(2z )\land \sqrt{t}}{\sqrt{t}}\cdot \frac{1}{\sqrt{2\pi t/3}}\lesssim \frac{z\land \sqrt{t}}{t}.
\end{align}
Therefore, the desired result is valid.
\hfill$\Box$
\bigskip

Fix $y,t>0$, we define for $s\in [0, t]$,
\[
h_s^{t,y}:=\frac{3}{2\sqrt{2\lambda^*}}\log \left(\frac{t+1}{t-s+1}\right)-y, \quad
f_s^{t,y}:= \sqrt{2\lambda^*}s - h_s^{t,y}.
\]
Then $f_t^{t,y}= \sqrt{2\lambda^*}t -\frac{3}{2\sqrt{2\lambda^*}}\log \left(t+1\right)+y\leq m(t)+y$ and $f_0^{t,y} =y > 0$.

The following result gives the upper bound.

\begin{prop}\label{Upper-Bound-M-T}
There exists a positive constant $C_0$  such that for any $y,t \geq 1$ and $i, i'\in S$,
\[
\P_{(0,i)} \left(M_t^{i'} \geq m(t) +y\right)\leq \P_{(0,i)} \left(M_t \geq m(t) +y\right)\leq C_0 (y\land \sqrt{t}) e^{-\sqrt{2\lambda^*} y}.
\]
\end{prop}
\textbf{Proof: }
The first inequality is trivial since
$M_t^{i'}\leq M_t$. Now we prove the second inequality.
Let $[x]$ be the largest integer less than or equal to $x$.  Then
\begin{align}\label{step_37}
	&\P_{(0,i)} \left(M_t \geq m(t) +y\right)\leq \P_{(0,i)} \left(M_t \geq f_{t}^{t,y}\right)\nonumber  \\ &\leq \sum_{k=0}^{[t]}\E_{(0,i)}\left(\sum_{u\in Z((k+1)\land t)} 1_{\left\{\sup_{s \in [k, (k+1)\land t]} X_u(s) \geq f_k^{t,y} \right\}}  1_{\left\{ X_u(s) \leq f_s^{t,y}, s\leq k \right\}} \right)\nonumber\\
	& = \sum_{k=1}^{[t]}\E_{(0,i)}\left(\sum_{u\in Z((k+1)\land t)} 1_{\left\{\sup_{s \in [k, (k+1)\land t]} X_u(s) \geq f_k^{t,y} \right\}}  1_{\left\{ X_u(s) \leq f_s^{t,y}, s\leq k \right\}} \right)=: \sum_{k=1}^{[t]} D_k .
\end{align}
Since all components of $\mathbf{h}$ are positive, we have by Proposition \ref{Many-to-one} that
\begin{align}\label{step_38}
	&D_k \lesssim  e^{\lambda^*(k+1)} \widehat{\P}_{(0,i)}\left(\sup_{s \in [k, (k+1)\land t]} X_\xi(s) \geq f_k^{t,y} ,\   X_\xi(s) \leq f_s^{t,y}, s\leq k  \right).
\end{align}
Note that, under $\widehat{\P}_{(0,i)}$, $X_\xi(t)$ is a standard Brownian motion. Thus,
\begin{align}\label{step_1}
	&\widehat{\P}_{(0,i)}\left(\sup_{s \in [k, (k+1)\land t]} X_\xi(s) \geq f_k^{t,y} ,\   X_\xi(s) \leq f_s^{t,y}, s\leq k  \right)\\
	&= \mathbf{E}_0 \left(\mathbf{P}_0\left(B_k\geq f_k^{t,y} -x , B_s \leq f_s^{t,y}, s\leq k\right)\big|_{x=\sup_{s\in [k, (k+1)\land t]}B_s - B_k } 	\right).
\end{align}
For any $\lambda \in \R$, define
\begin{align}\label{Change-of-measure2}
	\frac{\mathrm{d}\mathbf{P}_0^{\lambda}}{\mathrm{d}\mathbf{P}_0}\bigg|_{\sigma(B_s, s\leq t)}:= e^{\lambda B_t - \frac{1}{2}\lambda^2 t},
\end{align}
then under $\mathbf{P}_0^{\lambda}$, $B_t$ is a Brownian motion with drift $\lambda$. Using this change of measure, we get that
\begin{align}\label{step_2}
	&\mathbf{P}_0\left(B_k\geq f_k^{t,y} -x , B_s \leq f_s^{t,y}, s\leq k\right) = \mathbf{E}_0^{\sqrt{2\lambda^*}}\left(e^{-\sqrt{2\lambda^*}B_k +\lambda^* k}1_{\left\{B_k\geq f_k^{t,y} -x , B_s \leq f_s^{t,y}, s\leq k \right\}}
	\right)\nonumber\\ & \leq e^{-\sqrt{2\lambda^*}\left(f_k^{t,y}-x\right) +\lambda^* k}\mathbf{P}_0 \left(B_k +\sqrt{2\lambda^*}k \geq f_k^{t,y}-x, B_s +\sqrt{2\lambda^*}s \leq f_s^{t,y}, s\leq k\right)\nonumber\\
	& = e^{-\sqrt{2\lambda^*}\left(f_k^{t,y}-x\right) +\lambda^* k}\mathbf{P}_0 \left(B_k  \leq x +h_k^{t,y}, B_s  \geq h_s^{t,y}, s\leq k\right)\nonumber\\
	& \leq e^{-\sqrt{2\lambda^*}\left(f_k^{t,y}-x\right) +\lambda^* k}\mathbf{P}_0 \left( B_s  \geq h_s^{t,y+1}, s\leq k,\ B_k- h_{k}^{t,y+1} \in [1, x+1]\right)\nonumber\\
	& = : e^{-\sqrt{2\lambda^*}\left(f_k^{t,y}-x\right) +\lambda^* k} F_k^{t,y}\left([1,x+1]\right).
\end{align}
Let
\[
f(s):=\frac{3}{2\sqrt{2\lambda^*}}\log\left(\frac{t+1}{t-s+1}\right).
\]
Note  that, since  $\log (1+x)\lesssim  x^{1/4}, x>0$,
we have for all $s\leq k$ and all $1\leq k\leq t$,
\begin{align}
	& \frac{|f(s)|}{s^{1/4}}+ \frac{|f(k)- f(s)|}{(k-s)^{1/4}} \lesssim   \frac{1 }{s^{1/4}}\log\left(1+\frac{s}{t-s+1}\right) + \frac{1}{(k-s)^{1/4}} \log \left(1+ \frac{k-s}{t-k+1}\right)\\
	& \leq \frac{\log\left(1+s\right)}{s^{1/4}}+ \frac{\log\left(1+k-s\right)}{(k-s)^{1/4}} \lesssim 1.
\end{align}
Therefore,  applying  Lemma \ref{lemma1}  to the function $f$ above with $\alpha=1/4$, $y$ replaced by $y+1$
and $z = 1,.., [x]+1$, we get
 \begin{align}\label{step_67}
	&F_k^{t,y}\left([1,x+1]\right) \leq
	\sum_{z=1}^{[x]+1} F_k^{t,y}\left([z,z+1]\right) \nonumber\\
&\lesssim \frac{(y+1)\land \sqrt{k}}{k^{3/2}}\left(\sum_{z=1}^{[x]+1} z\land \sqrt{k}  \right)\lesssim \frac{y\land \sqrt{t}}{k^{3/2}} (x+1)^2,\quad y\geq 1, x\geq 0,\ 1 	\leq k\leq [t].
\end{align}
Plugging this upper bound into \eqref{step_2}, we get that
\begin{align}\label{step_3}
	\mathbf{P}_0\left(B_k\geq f_k^{t,y} -x , B_s \leq f_s^{t,y}, s\leq k\right) \lesssim e^{-\sqrt{2\lambda^*}\left(f_k^{t,y}-x\right) +\lambda^* k}\frac{y\land \sqrt{t}}{k^{3/2}}
	(x+1)^2.
\end{align}
Note that $0\leq \sup_{s \in [k, (k+1)\land t]} B_s - B_k \leq \sup_{s \in [k, k+1]} B_s - B_k$ which is equal in law to $W:= \sup_{s \in [0, 1]} B_s$ under $\mathbf{P}_0$. Combining this with \eqref{step_1} and \eqref{step_3},
we get that for all $y\geq 1, 1\leq k \leq \sqrt{t}$,
\begin{align}\label{step_39}
	& \widehat{\P}_{(0,i)}\left(\sup_{s \in [k, (k+1)\land t]} X_\xi(s) \geq f_k^{t,y} ,\   X_\xi(s) \leq f_s^{t,y}, s\leq k  \right)\nonumber\\
& \lesssim   e^{-\sqrt{2\lambda^*}f_k^{t,y} +\lambda^* k}\left(\frac{y\land \sqrt{t}}{k^{3/2}}\right) \mathbf{E}_0\left( e^{\sqrt{2\lambda^*}W}(W+1)^2\right)\asymp e^{-\sqrt{2\lambda^*}f_k^{t,y} +\lambda^* k}\left(\frac{y\land \sqrt{t}}{k^{3/2}}\right).
\end{align}
Combining \eqref{step_37}, \eqref{step_38} and \eqref{step_39}, we finally get that
\begin{align}
	 \P_{(0,i)} \left(M_t \geq m(t) +y\right)&\lesssim    \sum_{k=1}^{[t]}e^{\lambda^*(k+1)}e^{-\sqrt{2\lambda^*}f_k^{t,y} +\lambda^* k}\left(\frac{y\land \sqrt{t}}{k^{3/2}}\right)\\
	& \asymp   (y\land \sqrt{t})e^{-\sqrt{2\lambda^*}y}   \sum_{k=1}^{[t]}\left(\frac{t+1}{t-k+1}\right)^{3/2}\frac{1}{k^{3/2}},\quad t,y\geq 1.
\end{align}
Note that for all $t\geq 1$,
\begin{align}
	& \sum_{k=1}^{[t]}\left(\frac{t+1}{t-k+1}\right)^{3/2}\frac{1}{k^{3/2}} \leq \sum_{k=1}^{[t]}\left(\frac{[t]+2}{[t]-k+1}\right)^{3/2}\frac{1}{k^{3/2}} \\
	& \leq 2 \sum_{k=1}^{[[t]/2]+1}\left(\frac{3[t]}{[t]-k+1}\right)^{3/2}\frac{1}{k^{3/2}} \leq 2\left(\frac{3[t]}{[t]-[[t]/2]}\right)^{3/2}\sum_{k=1}^{[[t]/2]+1}\frac{1}{k^{3/2}}\\
	&\leq 2\left(\frac{3[t]}{[t]-[t]/2}\right)^{3/2}\sum_{k=1}^{\infty}\frac{1}{k^{3/2}} = 2\times 6^{3/2} \sum_{k=1}^{\infty}\frac{1}{k^{3/2}}\lesssim 1.
\end{align}
Therefore, for all $y\geq 1, t\geq 1$ and $i\in S$,
\begin{align}
	\P_{(0,i)} \left(M_t \geq m(t) +y\right) \lesssim (y\land \sqrt{t})e^{-\sqrt{2\lambda^*}y},
\end{align}
which is the desired result.
\hfill$\Box$
\bigskip

Next, we are going to get a lower bound for $\P_{(0,i)} \left(M_t \geq m(t) +y\right)$.
For $i' \in S$, let
\[
\mathcal{A}^{t,y}(i'):=  \# \left\{u\in Z(t): \  I_u(t)= i'\
 \mbox{and} \ \forall s\leq t,\ X_u(s)\leq f_s^{t,y},\ X_u(t)\geq f_t^{t,y}-4\right\}.
\]
Since $\lim_{t\to\infty}\widehat{\P}_{(0,i)}\left(I_\xi(t)=i' \right)= g_{i'}h_{i'}>0$ for all $i,i' \in S$, we have $\inf_{i,i' \in S} \inf_{t>1} \widehat{\P}_{(0,i)}\left(I_\xi(t)=i' \right)>0$.
Therefore, for $t>1$, by  Proposition \ref{Many-to-one}, and the independence of $I_\xi $ and $X_\xi$, we get that
\begin{align}\label{step_25}
	&\E_{(0,i)}\left( \mathcal{A}^{t,y}(i') \right)=e^{\lambda^* t}\widehat{\E}_{(0,i)}\left(\frac{h_i}{h_{I_\xi(t)}}1_{\left\{X_\xi(s)\leq f_s^{t,y}, s\leq t,\ X_\xi(t)\geq f_t^{t,y}-4 \right\}} 1_{\{I_\xi(t)=i' \}}\right)\nonumber\\
	&\gtrsim  e^{\lambda^* t}\mathbf{P}_0\left(B_s\leq f_s^{t,y}, s\leq t,\ B_t\geq f_t^{t,y}-4 \right).
\end{align}

We first show that for all
$ 4\sqrt{t}-3\geq y\geq 1, t\geq1$,
\begin{align}\label{step_22}
	e^{\lambda^* t}\mathbf{P}_0\left(B_s \leq y +\frac{s}{t}m(t),\ s\leq t,\quad B_t \geq y +m(t) -1 \right)\gtrsim y e^{-\sqrt{2\lambda^*}y}.
\end{align}
Let $q_t:= m(t)/t$. Taking $\lambda=q_t$ in \eqref{Change-of-measure2}, we get that for all $t, y\geq 1$,
\begin{align}
	&e^{\lambda^* t}\mathbf{P}_0\left(B_s \leq y +\frac{s}{t}m(t),\ s\leq t,\quad B_t \geq y +m(t) -1 \right)\\
	&= e^{\lambda^* t}\mathbf{E}_0^{q_t}\left(e^{-q_t B_t + q_t^2 t/2};B_s \leq y +q_t s,\ s\leq t,\quad B_t \geq y + q_t t -1 \right)\\
	& \geq e^{\lambda^* t} e^{-q_t(y+q_t t)+ q_t^2 t/2}\mathbf{P}_0\left(B_s \geq -y,\ s\leq t, \quad B_t \leq -y+1 \right) \\
	& \gtrsim e^{-\sqrt{2\lambda^*}y}t^{3/2}\mathbf{P}_0\left(B_s \geq -y\ s\leq t; \ B_t \leq -y+1 \right)=e^{-\sqrt{2\lambda^*}y}t^{3/2}
	\mathbb{E}_{\Q_y}
	\left(\frac{y}{R_t} ;R_t\leq 1\right),
\end{align}
where $(R_t, \Q_{y})$ is a Bessel-3 process starting from $y$, and in the last equation we use the following well-known change-of-measure
\[
\frac{\mathrm{d} \Q_y}{\mathrm{d} \mathbf{P}_0}\bigg|_{\sigma(B_s, s\leq t) } = \frac{B_t +y}{y}1_{\{ B_s+y \geq 0, \ s\leq t\}},
\]
here $(B_t+y, \ \Q_y)$ is equal in law to $(R_t, \Q_y)$.
The density of $R_t$ under $\Q_{y}$ is given by
\[
\frac{x}{y\sqrt{2\pi t}}
e^{-(x-y)^2/(2t)}
\left(1- e^{-2xy/t}\right)1_{\{x >0 \} }.
\]
Note that $1-e^{-x}\gtrsim x, 0\leq x \leq 8$, and that $0\leq 2xy/t\leq 8\sqrt{t}/t\leq 8$
and $(x-y)^2 / (2t) \leq (4\sqrt{t})^2/(2t)\lesssim 1$
for all $0\leq x \leq 1, 1\leq y\leq 4\sqrt{t}-3<4\sqrt{t}$ and $t\geq 1$.
Thus,
\begin{align}
	& e^{-\sqrt{2\lambda^*}y}t^{3/2}
		\mathbb{E}_{\Q_y}
	\left(\frac{y}{R_t} ;R_t\leq 1\right)= y e^{-\sqrt{2\lambda^*}y} t^{3/2}\int_0^1 \frac{1}{y \sqrt{2\pi t}}
	e^{-(x-y)^2/(2t)}
	 (1- e^{-2xy/t})\mathrm{d} x\\
	& \gtrsim y e^{-\sqrt{2\lambda^*}y} t^{3/2}\int_0^1 \frac{1}{y \sqrt{2\pi t}}
\frac{2xy}{t}\mathrm{d} x\\
	& \gtrsim y e^{-\sqrt{2\lambda^*}y}\int_0^1 x\mathrm{d} x \gtrsim y e^{-\sqrt{2\lambda^*}y},
	\quad 1\leq y\leq  4\sqrt{t}-3,\ t\geq 1.
\end{align}
Thus \eqref{step_22} is true.

Since, for fixed $t$,
\begin{align}
	\frac{s}{t}m(t) +y - \left(f_s^{t,y}+3\right)=\frac{3}{2\sqrt{2\lambda^*}} \left(\log \left(\frac{t+1}{t-s+1} \right) - \frac{s}{t}\log t\right)-3=:G(s),
\end{align}
we see that
\[
G'(s)= \frac{3}{2\sqrt{2\lambda^*}}\left(\frac{1}{t-s+1}-\frac{\log t}{t}\right)= 0\ \Longleftrightarrow\ s= t+1-\frac{t}{\log t}.
\]
Since $G(0)=0$, $G'(0)<0$ and $G'(t)>0$, $t \geq e$, we have for $t\geq e$ and $s\le t$,
\begin{align}
	G(s)\leq
	G(t)=
	-3 \frac{3}{2\sqrt{2\lambda^*}}\log \left(1+\frac{1}{t}\right)\leq 0.
\end{align}
Thus
\begin{equation}\label{e:eqa}
\left\{ B_s \leq y +\frac{s}{t}m(t),\ s\leq t \right\} \subset \left\{ B_s \leq f_s^{t,y}+3,\quad s\leq t  \right\}.
\end{equation}
On the other hand, since
\begin{align}
	y +m(t) -1 - \left(f_t^{t,y} +3 -4\right)= \frac{3}{2\sqrt{2\lambda^*}}\log \left(1+\frac{1}{t}\right)\geq 0,
\end{align}
we have that
\begin{equation}\label{e:eq2}
\left\{ B_t \geq y +m(t) -1 \right\} \subset \left\{ B_t \geq f_t^{t,y} +3 -4 \right\}.
\end{equation}
Combining \eqref{e:eqa} and \eqref{e:eq2},  and noting that $f_{s}^{t,y}+3 = f_s^{t, y+3}$, we get that
\begin{equation}\label{ineq-m-f}
\mathbf{P}_0\left(B_s \leq y +\frac{s}{t}m(t),\ s\leq t,\quad B_t \geq y +m(t) -1 \right)\le \mathbf{P}_0\left(B_s\leq f_s^{t,y+3}, s\leq t,\ B_t\geq f_t^{t,y+3}-4 \right).
\end{equation}
By \eqref{step_25}, \eqref{step_22}, and \eqref{ineq-m-f},
\begin{align}\label{step_24}
	 \E_{(0,i)}\left( \mathcal{A}^{t,y}(i') \right) \gtrsim (y-3) e^{-\sqrt{2\lambda^*}(y-3)} \gtrsim y e^{-\sqrt{2\lambda^*}y},\quad i,i'\in S,
	\quad t\geq 1,\quad  4 \leq y \leq 4\sqrt{t}.
\end{align}

Now we state our result for the lower bound.

\begin{prop}\label{Lower-Bound-M-T}
      Let $\beta_0:= 4+ \frac{3}{2\sqrt{2\lambda^*}}\log 2$.
      There exists a positive constant $K_0$ such that for any
      $t\geq (\beta_0/3)^2, y\in [1, \sqrt{t}]$ and $i,i'\in S$,
	\[
	\P_{(0,i)} \left(M_t\geq m(t) +y\right)\geq  \P_{(0,i)} \left(M_t^{i'} \geq m(t) +y\right)
	\geq K_0 y e^{-\sqrt{2\lambda^*} y}.
	\]
\end{prop}
\textbf{Proof: }
The first inequality is trivial, so we only need to  prove the second inequality.
It is easy to see that for all
$t\geq 1, y\geq 1$,
\begin{align}
	\P_{(0,i)}\left(M_t^{i'} \geq m(t) +y -\beta_0 \right)\geq \P_{(0,i)}\left(\mathcal{A}^{t,y}(i')\geq 1\right).
\end{align}
If we can prove that for all
$t\geq 1, y\in [4, 4\sqrt{t}]$
and $i,i'\in S$,
\begin{equation}\label{toprove}
 \P_{(0,i)}\left(\mathcal{A}^{t,y}(i')\geq 1\right)\gtrsim  y
  e^{-\sqrt{2\lambda^*}y},
\end{equation}
then for any $t\geq (\beta_0/3)^2$ (which is equivalent to $3\sqrt{t}\geq \beta_0$)
and $1\leq y\leq \sqrt{t}$, we have $4\leq y+\beta_0 \leq 4\sqrt{t}$,  and  thus
	\begin{align*}
	& \P_{(0,i)}\left(M_t^{i'} \geq m(t) +y  \right)
= \P_{(0,i)}\left(M_t^{i'} \geq m(t) +y+\beta_0 - \beta_0 \right)\geq \P_{(0,i)}\left(\mathcal{A}^{t,y+ \beta_0 }(i')\geq 1\right)\\
&	 \gtrsim \left(y+ \beta_0 \right)	e^{-\sqrt{2\lambda^*}(y+ \beta_0 )}	\gtrsim y e^{-\sqrt{2\lambda^*}y},
\end{align*}
which completes the proof. To prove \eqref{toprove},
we use the trivial inequality $\mathbb{E}\left(|Y|^{1+\alpha_0}\right)\mathbb{E}(1_A)^{\alpha_0}\geq \mathbb{E}\left(|Y|1_A\right)^{1+\alpha_0}$ to
get that for all
$t\geq 1, 4\leq y\leq 4\sqrt{t}$ and
$i,i'\in S$,
\begin{align}
	\P_{(0,i)}\left(\mathcal{A}^{t,y}(i')\geq 1\right)\geq \left\{\frac{\left(\E_{(0,i)}\left( \mathcal{A}^{t,y}(i')\right)\right)^{1+\alpha_0}}{\E_{(0,i)}\left(\left(\mathcal{A}^{t,y}(i')\right)^{1+\alpha_0} \right)}\right\}^{1/\alpha_0}.
\end{align}
If we can prove that for all $y,t \geq 1$ and
$i,i'\in S$,
\begin{equation}\label{toprove2}
\E_{(0,i)}\left( (\mathcal{A}^{t,y}(i'))^{1+\alpha_0} \right)
\lesssim (y\land \sqrt{t})e^{-\sqrt{2\lambda^*}y},
\end{equation}
then using \eqref{step_24}, we get \eqref{toprove}.
Now we  prove \eqref{toprove2}.

{\bf Step 1}\quad
For $u\in Z(t),$ define
\[
\Upsilon(u):= \left\{X_u(s)\leq f_s^{t,y},\ s\leq t, X_u(t)\geq f_t^{t,y}-4 \right\}.
\]
We first estimate $\E_{(0,i)}\left( (\mathcal{A}^{t,y}(i'))^{1+\alpha_0} \right)$ from above.
Define $\mathcal{A}^{t,y}:= \sum_{j=1}^d \mathcal{A}^{t,y}(j)$. It follows from Proposition \ref{Many-to-one} that
\begin{align}\label{step_48}
	&\E_{(0,i)}\left( (\mathcal{A}^{t,y}(i'))^{1+\alpha_0} \right)
	\leq \E_{(0,i)}\left( (\mathcal{A}^{t,y})^{1+\alpha_0} \right)
	= \E_{(0,i)}\left( \sum_{u\in Z(t)} (\mathcal{A}^{t,y})^{\alpha_0} 1_{\Upsilon(u)}  \right)\nonumber\\
	& = h_i e^{\lambda^* t}\widehat{\E}_{(0,i)}\left( \left(\mathcal{A}^{t,y} \right)^{\alpha_0}\frac{1}{h_{I_\xi(t)}} 1_{\Upsilon(\xi)}\right) \lesssim e^{\lambda^* t} \widehat{\E}_{(0,i)}\left( \left(\mathcal{A}^{t,y}\right)^{\alpha_0}  1_{\Upsilon(\xi)}	\right).
\end{align}

In the following substeps (i) and (ii), we estimate $ \widehat{\E}_{(0,i)}\left( \left(\mathcal{A}^{t,y}\right)^{\alpha_0}  1_{\Upsilon(\xi)}	\right)$ from above.

{\bf Substep (i)}\quad
Let $\mathcal{G}_t:= \sigma\left(X_\xi(s), I_\xi(s), s\leq t\right)$. Let $\tau_k$ be the $k$-th time that $I_\xi(s-)\neq I_\xi(s)$ and let $\mathbf{A}\left(I_\xi(\tau_{\ell}-)\right)= \left(A_1\left(I_\xi(\tau_{\ell}-)\right),...,A_d\left(I_\xi(\tau_{\ell}-)\right)\right)^T$.
Let $K_t:= \sup\left\{k : \tau_k\leq t\right\}$. By the trivial inequality
$\left(\sum_{j=1}^k x_j\right)^{\alpha_0}\leq \sum_{j=1}^k x_j^{\alpha_0}$,
on the event $\Upsilon(\xi)$,
\begin{align}\label{step_5}
	&\widehat{\E}_{(0,i)}\left(\left( \mathcal{A}^{t,y}\right)^{\alpha_0} \big| \mathcal{G}_t\cap\left(\mathbf{A}\left(I_\xi(\tau_{\ell}-)\right), \ell\leq K_t \right) \right)\leq 1+ \sum_{\ell=1}^{K_t}\sum_{j=1}^d \left(A_{j}\left(I_\xi(\tau_{\ell}-)\right)- \delta_{j, I_\xi(\tau_\ell)}\right)^{\alpha_0} \\ &\quad \times \left\{\E_{(X_\xi(\tau_{\ell}),j)} \left(\# \left\{u\in Z(t-\tau_{\ell}):\ \forall s\leq t-\tau_{\ell}, X_u(s)\leq f_{s+\tau_{\ell}}^{t,y}, X_u(t-\tau_{\ell})\geq f_t^{t,y}-4 \right\}\right)\right\}^{\alpha_0}.
\end{align}
Let $z= X_\xi(\tau_{\ell}) \leq f_{\tau_\ell}^{t,y}, z':= f_{\tau_\ell}^{t,y}-z \geq 0$ and $r = t-\tau_\ell$, then
\begin{align}\label{step_4}
	&\E_{(z,j)} \left(\# \left\{u\in Z(r):\ \forall s\leq r, X_u(s)\leq f_{s+\tau_{\ell}}^{t,y}, X_u(r)\geq f_t^{t,y}-4 \right\}\right)\nonumber \\
	& \lesssim e^{\lambda^* r}\mathbf{P}_0 \left(B_s \leq f_{s+\tau_{\ell}}^{t,y}-z, s\leq r,\ B_r \geq f_t^{t,y}-4-z\right)\nonumber\\
	& =  e^{\lambda^* r}\mathbf{E}_0^{\sqrt{2\lambda^*}} \left(e^{-\sqrt{2\lambda^*}B_r + \lambda^* r}1_{\{B_s \leq f_{s+\tau_{\ell}}^{t,y}-z, s\leq r,\ B_r \geq f_t^{t,y}-4-z\}}\right)\nonumber\\
	& \lesssim e^{-\sqrt{2\lambda^*}\left(f_t^{t,y}-4-z\right) +2\lambda^* r} \mathbf{P}_0\left(B_s \geq h_s^{r,z'}, s\leq r,\ B_r \leq 4+h_r^{r,z'} \right).
\end{align}
Recall that $h_s^{t,y}= \frac{3}{2\sqrt{2\lambda^*}}\log \left(\frac{t+1}{t-s+1}\right)-y$.
It follows from Lemma \ref{lemma1} that for  $r\geq 1$,
\begin{align}
	&\mathbf{P}_0\left(B_s \geq h_s^{r,z'}, s\leq r,\ B_r \leq 4+h_r^{r,z'} \right)\leq \mathbf{P}_0\left(B_s \geq h_s^{r,z'+1}, s\leq r,\ B_r \leq 3+h_r^{r,z'+1} \right) \lesssim \frac{ (z'+1)}{r^{3/2}}.
\end{align}
For $r <1$, we use the trivial bound $1$.
 Plugging these
into \eqref{step_4}, together with \eqref{step_5}, we conclude that on the event $\Upsilon(\xi)$,
\begin{align}\label{step_7}
	&\widehat{\E}_{(0,i)}\left(\left( \mathcal{A}^{t,y} \right)^{\alpha_0}\big| \mathcal{G}_t\cap\left(\mathbf{A}\left(I_\xi(\tau_{\ell}-)\right), \ell\leq K_t\right) \right)\lesssim 1+ \sum_{\ell=1}^{K_t}\sum_{j=1}^d (A_{j}\left(I_\xi(\tau_{\ell}-)\right))^{\alpha_0} \\
	& \quad \times e^{-\alpha_0 \sqrt{2\lambda^*}\left(f_t^{t,y}-2-X_\xi(\tau_\ell)\right) +2\alpha_0\lambda^* (t-\tau_\ell)}  \left(\frac{ (f_{\tau_\ell}^{t,y}-X_\xi(\tau_\ell)+1)}{(t-\tau_\ell)^{3/2}} 1_{\left\{t-\tau_\ell\geq 1 \right\}} + 1_{\left\{t-\tau_\ell\leq 1\right\}} \right)^{\alpha_0}.
\end{align}
Note that the distribution of the number of offspring $\mathbf{A}(j)$ of a spine particle of type  $j$ is given by
\[
\widehat{\P}_{(0,i)}\left(\mathbf{A}(j) = \mathbf{k}\right)=\frac{p_{\mathbf{k}}(j) \langle \mathbf{k}, \mathbf{h} \rangle}{(1+\lambda^*/ a_j)h_j}=: \widehat{p}_{\mathbf{k}}(j).
\]
Thus, given $\mathcal{G}_t$,
    the law of  $\mathbf{A}\left(I_\xi(\tau_{\ell}-)\right)$ is equal to  $\widehat{\P}_{(0,i)}\left(\cdot \Big| A_{I_\xi(\tau_{\ell})}\left(I_\xi(\tau_{\ell}-)\right)\geq 1,  I_\xi(\tau_{\ell}-), I_\xi(\tau_{\ell}) \right)$ since
 there must be at least one particle of   type  $I_\xi(\tau_{\ell})$ among the $\mathbf{A}\left(I_\xi(\tau_{\ell}-)\right)$ offspring.
    So for any $\mathbf{k}\in \mathbb{N}^d$,
\begin{align}
	\widehat{\P}_{(0,i)}\left(\mathbf{A}\left(I_\xi(\tau_{\ell}-)\right) = \mathbf{k} \Big| \mathcal{G}_t \right)
&= \frac{\widehat{\P}_{(0,i)}\left(\mathbf{A}\left(I_\xi(\tau_{\ell}-)\right) = \mathbf{k} , A_{I_\xi(\tau_{\ell})}\left(I_\xi(\tau_{\ell}-)\right)\geq 1 \Big| I_\xi(\tau_{\ell}-), I_\xi(\tau_{\ell}) \right)}{\widehat{\P}_{(0,i)}\left(A_{I_\xi(\tau_{\ell})}\left(I_\xi(\tau_{\ell}-)\right)\geq 1 \Big| I_\xi(\tau_{\ell}-), I_\xi(\tau_{\ell}) \right)}\\
	& =\frac{1}{\sum_{\mathbf{k}\in \N^d: k_{I_\xi(\tau_{\ell})}\geq 1 } \widehat{p}_{\mathbf{k}}(I_\xi (\tau_{\ell}-)) } \widehat{p}_{\mathbf{k}}(I_\xi (\tau_{\ell}-)) 1_{\{k_{I_\xi(\tau_{\ell})}\geq 1 \}}.
\end{align}
Define
\[
\mathcal{S}:= \left\{(j_1, j_2)\in S \times S: \sum_{\mathbf{k}\in \N^d: k_{j_2}\geq 1 } \widehat{p}_{\mathbf{k}}(j_1)>0 \right\}.
\]
There exists a constant $c_1\in (0, 1]$ such that for all $(j_1, j_2)\in \mathcal{S}$, $\sum_{\mathbf{k}\in \N^d: k_{j_2}\geq 1 } \widehat{p}_{\mathbf{k}}(j_1)\geq c_1$.
Note that $\sum_{\mathbf{k}\in \N^d: k_{I_\xi(\tau_{\ell})}\geq 1 } \widehat{p}_{\mathbf{k}}(I_\xi (\tau_{\ell}-))>0$.
Therefore, $\left(I_\xi (\tau_{\ell}-), I_\xi (\tau_{\ell}) \right)\in \mathcal{S}$ and
\begin{align}\label{step_8}
&\widehat{\E}_{(0,i)}\left( \sum_{j=1}^d \left(A_j (I_\xi(\tau_\ell)-)\right)^{\alpha_0} \Big| \mathcal{G}_t \right)= \frac{1}{\sum_{\mathbf{k}\in \N^d: k_{I_\xi(\tau_{\ell})}\geq 1 }\widehat{p}_{\mathbf{k}}(I_\xi (\tau_{\ell}-))  } \sum_{j=1}^d \sum_{\mathbf{k}\in \N^d: k_{I_\xi(\tau_{\ell})}\geq 1} \widehat{p}_{\mathbf{k}}(I_\xi(\tau_\ell)-)  k_j^{\alpha_0} \nonumber\\
	& \leq \frac{1}{c_1}\sum_{j=1}^d \sum_{\mathbf{k}\in \N^d} \widehat{p}_{\mathbf{k}}(I_\xi(\tau_\ell)-)  k_j^{\alpha_0} \lesssim \sup_{\ell \in S}\sum_{j=1}^d \sum_{\mathbf{k} \in \N^d} \widehat{p}_{\mathbf{k}}(\ell)k_j^{\alpha_0}\nonumber
	\\& =  \sup_{\ell \in S} \sum_{j=1}^d \sum_{\mathbf{k} \in \N^d} \frac{p_{\mathbf{k}}(\ell) \langle \mathbf{k}, \mathbf{h} \rangle}{(1+\lambda^*/ a_\ell)h_\ell} k_j^{\alpha_0}\lesssim \sup_{\ell , j, q \in S}\sum_{\mathbf{k} \in \N^d} p_{\mathbf{k}}(\ell) k_q k_j^{\alpha_0} \nonumber\\
	&\leq \sup_{\ell , j, q \in S}\sum_{\mathbf{k} \in \N^d} p_{\mathbf{k}}(\ell) \left(\frac{k_q^{1+\alpha_0}}{1+\alpha_0} + \frac{\alpha_0 k_j^{1+\alpha_0}}{1+\alpha_0} \right) \lesssim 1,
\end{align}
where in the last inequality we used the assumption \eqref{Second-Moment}.
By \eqref{step_7} and \eqref{step_8},  on the event $\Upsilon(\xi)$,
\begin{align}\label{step_9}
	&\widehat{\E}_{(0,i)}\left( \left(\mathcal{A}^{t,y}\right)^{\alpha_0} \big| \mathcal{G}_t \right)\nonumber \lesssim 1+  \sum_{\ell=1}^{K_t} e^{-\alpha_0\sqrt{2\lambda^*}\left(f_t^{t,y}-2-X_\xi(\tau_\ell)\right) +2\alpha_0\lambda^* (t-\tau_{\ell})} \\
	& \quad \quad \times  \left(\frac{ (f_{\tau_\ell}^{t,y}-X_\xi(\tau_\ell)+1)}{(t-\tau_\ell)^{3/2}} 1_{\left\{t-\tau_\ell\geq 1 \right\}} + 1_{\left\{t-\tau_\ell\leq 1\right\}} \right)^{\alpha_0}.
\end{align}

{\bf Substep (ii)}\quad
Note that $\tau_\ell$ is measurable with respect to $\sigma(I_\xi(s): s\geq 0)$, which means that $\tau_\ell$ is independent of $X_\xi$. For the Brownian motion $B$, define
\[
\Upsilon:= \left\{B_s\leq f_s^{t,y},\ s\leq t, B_t\geq f_t^{t,y}-4 \right\}.
\]
Then
\begin{align}\label{step_59}
	& \widehat{\E}_{(0,i)}\left( \left(\mathcal{A}^{t,y}\right)^{\alpha_0}1_{\Upsilon(\xi)} \big| \tau_{\ell}: \ell \geq 1 \right) \lesssim \mathbf{P}_0 \left(B_s\leq f_s^{t,y},\ s\leq t, B_t\geq f_t^{t,y}-4 \right)\nonumber\\&\quad + \sum_{\ell=1}^{K_t} \mathbf{E}_0\left(e^{-\alpha_0\sqrt{2\lambda^*}\left(f_t^{t,y}-B_{\tau_\ell}\right) +2\alpha_0\lambda^* (t-\tau_\ell)}  \left(\frac{ f_{\tau_\ell}^{t,y}-B_{\tau_\ell}+1}{(t-\tau_\ell )^{3/2}} 1_{\{t-\tau_\ell \geq 1\}}+ 1_{\{t-\tau_\ell \leq 1\}}\right)^{\alpha_0}1_{\Upsilon} \right)\nonumber\\ &=: L_1 + \sum_{\ell =1}^{K_t}L_2(\ell).
\end{align}

For $L_1$, note that the argument leading to \eqref{step_3} also works when $k$ is not an integer. Letting $k=t$ and $x=4$ in \eqref{step_3}, we get
\begin{align} \label{ineq-L_1}
	& L_1 \lesssim  e^{-\sqrt{2\lambda^*}\left(f_t^{t,y}-4\right) +\lambda^* t}\frac{(y+1)\land \sqrt{t}}{t^{3/2}}\times 25 \lesssim  e^{-\lambda^* t}
(y\land \sqrt{t})e^{-\sqrt{2\lambda^*}y}.
\end{align}

Note that $\tau_\ell$ can be regarded as a constant with respect to $\mathbf{E}_0$.
For $L_2(\ell)$, we deal with two cases separately. We first deal with the case $t-\tau_\ell \geq 1$. Set $r:= \tau_\ell$. By the Markov property,
\begin{align}\label{step_52}
	L_2(\ell)&= \mathbf{E}_0\left(e^{-\alpha_0\sqrt{2\lambda^*}\left(f_t^{t,y}-B_r\right) +2\alpha_0\lambda^* (t-r)}  \left(\frac{ f_{r}^{t,y}-B_r+1}{(t-r)^{3/2}} \right)^{\alpha_0}1_{\Upsilon} \right)\nonumber\\
	& =  \mathbf{E}_0\bigg( e^{-\alpha_0\sqrt{2\lambda^*}\left(f_t^{t,y}-B_r\right) +2\alpha_0\lambda^* (t-r)} \left(\frac{f_{r}^{t,y}-B_r+1}{(t-r)^{3/2}} \right)^{\alpha_0}1_{\left\{B_s \leq f_s^{t,y},\ s\leq r \right\}} \nonumber\\
	& \quad \times \mathbf{P}_{B_r}\left(B_s \leq f_{s+r}^{t,y},\ s\leq t-r,\ B_{t-r}\geq f_t^{t,y}-4\right) \bigg).
\end{align}
For  $z\leq f_r^{t,y}$, set $z':= f_r^{t,y}-z$. Using \eqref{Change-of-measure2},
the fact that $f_{s+r}^{t,y}-z = \sqrt{2\lambda^*}s - h_s^{t-r,z'}$
and \eqref{step_67} (which is still valid when $k=t$ is not an integer)
with
$F_{t-r}^{t-r,z'}([1,5])$ defined in \eqref{step_2}, we get that
\begin{align}\label{step_53}
	&\mathbf{P}_z\left(B_s \leq f_{s+r}^{t,y},\ s\leq t-r,\ B_{t-r}\geq f_t^{t,y}-4\right)\nonumber\\
	&= \mathbf{E}_0^{\sqrt{2\lambda^*}}\left(e^{-\sqrt{2\lambda^*} B_{t-r}+\lambda^*(t-r)} 1_{\left\{B_s \leq f_{s+r}^{t,y}-z,\ s\leq t-r,\ B_{t-r}\geq f_t^{t,y}-4-z \right\}}\right)\nonumber \\
	& \leq e^{-\sqrt{2\lambda^*} \left(f_t^{t,y}-z\right)+\lambda^*(t-r)}\mathbf{P}_0 \left( B_s \leq -h_s^{t-r,z'},\ s\leq t-r,\ B_{t-r}\geq-4-h_{t-r}^{t-r,z'} \right)\nonumber \\
	&=e^{-\sqrt{2\lambda^*} \left(f_t^{t,y}-z\right)+\lambda^*(t-r)}\mathbf{P}_0 \left( B_s \geq h_s^{t-r,z'},\ s\leq t-r,\ B_{t-r}\leq 4+h_{t-r}^{t-r,z'} \right)\nonumber \\
&=e^{-\sqrt{2\lambda^*} \left(f_t^{t,y}-z\right)+\lambda^*(t-r)}\mathbf{P}_0 \left( B_s \geq 1+h_s^{t-r,z'+1},\ s\leq t-r,\ B_{t-r}\leq 5+h_{t-r}^{t-r,z'+1} \right)\nonumber \\
	& \lesssim e^{-\sqrt{2\lambda^*} \left(f_t^{t,y}-z\right)+\lambda^*(t-r)}\left(\frac{z'+1}{(t-r)^{3/2}}\right).
\end{align}
Combining \eqref{step_52} and \eqref{step_53}, we get
that
\begin{align}\label{step_54}
	L_2(\ell) \lesssim & \mathbf{E}_0\left( e^{-(1+\alpha_0)\sqrt{2\lambda^*}\left(f_t^{t,y}-B_{\tau_\ell}\right) +(1+2\alpha_0)\lambda^* (t-\tau_{\ell})} \left(\frac{f_{\tau_\ell}^{t,y}-B_{\tau_\ell}+1}{(t-\tau_\ell)^{3/2}} \right)^{1+\alpha_0}1_{\left\{B_s \leq f_s^{t,y},\ s\leq \tau_{\ell} \right\}} \right)\nonumber \\
	=:& \frac{1}{(t-\tau_\ell)^{3(1+\alpha_0)/2}} e^{-(1+\alpha_0)\sqrt{2\lambda^*}f_t^{t,y}+(1+2\alpha_0)\lambda^*(t-\tau_\ell) } \mathbf{E}_0\left(G1_{\left\{B_s \leq f_s^{t,y},\ s\leq \tau_{\ell} \right\}}\right),
\end{align}
where $G:= e^{(1+\alpha_0)\sqrt{2\lambda^*}B_{\tau_\ell}}(f_{\tau_{\ell}}^{t,y}-B_{\tau_{\ell}} +1)^{1+\alpha_0}.$
Next,
using \eqref{Change-of-measure2}  with $\lambda = \sqrt{2\lambda^*}$, and noticing that
 $(B_s, \mathbf{P}_0^{\sqrt{2\lambda^*}})\stackrel{\mathrm{d}}{=}(B_s+\sqrt{2\lambda^*}s, \mathbf{P}_0)$
and that $f_s^{t,y}=\sqrt{2\lambda^*}s-h_s^{t,y}$, we get that
\begin{align}
	& \mathbf{E}_0\left(G1_{\left\{B_s \leq f_s^{t,y},\ s\leq \tau_{\ell} \right\}}\right)=\sum_{k=0}^\infty \mathbf{E}_0\left(G1_{\left\{B_s \leq f_s^{t,y},\ s\leq \tau_{\ell} \right\}} 1_{\{B_{\tau_{\ell}}- f_{\tau_{\ell}}^{t,y}\in [-k-1, -k] \}}\right)\\
	&\leq \sum_{k=0}^\infty e^{\alpha_0\sqrt{2\lambda^*}\left(f_{\tau_{\ell}}^{t,y}-k\right)}(k+2)^{1+\alpha_0} \mathbf{E}_0\left(e^{\sqrt{2\lambda^*}B_{\tau_{\ell}}} 1_{\left\{B_s \leq f_s^{t,y},\ s\leq \tau_{\ell} \right\}} 1_{\{B_{\tau_{\ell}}- f_{\tau_{\ell}}^{t,y}\in [-k-1, -k] \}} \right) \\
	&=e^{\alpha_0 \sqrt{2\lambda^*}f_{\tau_{\ell}}^{t,y} +\lambda^* \tau_{\ell}} \sum_{k=0}^\infty e^{-\alpha_0\sqrt{2\lambda^*}k}(k+2)^{1+\alpha_0}
 \mathbf{P}_0 \left(B_s\leq -h_s^{t,y},\ s\leq \tau_{\ell} , B_{\tau_{\ell}}+ h_{\tau_{\ell}}^{t,y}\in [-(k+1), -k] \right).
\end{align}
When $\tau_\ell\geq 1$, by Lemma \ref{lemma1}, we have for all $k\geq 0$,
\begin{align}\label{step_56}
	&  \mathbf{P}_0 \left(B_s\leq -h_s^{t,y},\ s\leq \tau_{\ell} , B_{\tau_{\ell}}+ h_{\tau_{\ell}}^{t,y}\in [-(k+1), -k] \right)\nonumber\\
	&= \mathbf{P}_0 \left(B_s\geq h_s^{t,y},\ s\leq \tau_{\ell} , B_{\tau_{\ell}}- h_{\tau_{\ell}}^{t,y}\in [k,(k+1)] \right)\nonumber\\
	& \leq \mathbf{P}_0\left(B_s \geq  h_s^{t,y+1},\ s\leq \tau_\ell,\ B_r-h_{\tau_{\ell}}^{t,y+1} \in [k+1, k+2] \right)\nonumber\\
	& \lesssim  \frac{((y+1)\land \sqrt{\tau_{\ell}})(k+1)}{\tau_{\ell}^{3/2}} \lesssim  \frac{(y\land \sqrt{t})(k+1)}{\tau_{\ell}^{3/2}}.
\end{align}
When $\tau_\ell \leq 1$,
we use the trivial upper bound $1$. Therefore, using the fact that $\sum_{k=0}^\infty e^{-\alpha_0\sqrt{2\lambda^*}k}(k+2)^{1+\alpha_0} (k+1)<\infty$,
 we conclude that
\begin{align}\label{step_55}
	&\mathbf{E}_0\left(G1_{\left\{B_s \leq f_s^{t,y},\ s\leq \tau_{\ell} \right\}}\right)\nonumber\\
	& \lesssim  e^{\alpha_0 \sqrt{2\lambda^*}f_{\tau_{\ell}}^{t,y} +\lambda^* \tau_{\ell}} \sum_{k=0}^\infty e^{-\alpha_0\sqrt{2\lambda^*}k}(k+2)^{1+\alpha_0} \left(\frac{(y\land \sqrt{t})(k+1)}{\tau_{\ell}^{3/2}}1_{\{\tau_\ell \geq 1\}}+ 1_{\{\tau_\ell \leq 1\}} \right)\nonumber \\
	& \lesssim \left(y\land \sqrt{t}\right)e^{\alpha_0 \sqrt{2\lambda^*}f_{\tau_{\ell}}^{t,y} +\lambda^* \tau_{\ell}}  \left(\frac{1}{\tau_{\ell}^{3/2}}1_{\{\tau_\ell \geq 1\}}+ 1_{\{\tau_\ell \leq 1\}} \right).
\end{align}
Combining \eqref{step_54} and \eqref{step_55}, we get that in the case when $t-\tau_\ell \geq 1$,
\begin{align}\label{step_57}
	L_2(\ell)&\lesssim \left(y\land \sqrt{t}\right)\frac{e^{-(1+\alpha_0)\sqrt{2\lambda^*}f_t^{t,y}+(1+2\alpha_0)\lambda^*(t-\tau_\ell) }e^{\alpha_0 \sqrt{2\lambda^*}f_{\tau_{\ell}}^{t,y} +\lambda^* \tau_{\ell}}  }{(t-\tau_\ell)^{3(1+\alpha_0)/2}}  \left(\frac{1}{\tau_{\ell}^{3/2}}1_{\{\tau_\ell \geq 1\}}+ 1_{\{\tau_\ell \leq 1\}} \right) \nonumber \\
	& =\left(y\land \sqrt{t}\right)e^{-\sqrt{2\lambda^*}y} \frac{(t+1)^{3/2}(t-\tau_\ell+1)^{3\alpha_0/2}}{(t-\tau_\ell)^{3(1+\alpha_0)/2}} e^{-\lambda^* t}\left(\frac{1}{\tau_{\ell}^{3/2}}1_{\{\tau_\ell \geq 1\}}+ 1_{\{\tau_\ell \leq 1\}} \right)\nonumber \\
	&\lesssim \left(y\land \sqrt{t}\right)e^{-\sqrt{2\lambda^*}y} \frac{(t+1)^{3/2}}{(t-\tau_\ell)^{3/2}} e^{-\lambda^* t}\left(\frac{1}{\tau_{\ell}^{3/2}}1_{\{\tau_\ell \geq 1\}}+ 1_{\{\tau_\ell \leq 1\}} \right)\nonumber\\
	&  \lesssim  \left(y\land \sqrt{t}\right)e^{-\sqrt{2\lambda^*}y}  e^{-\lambda^* t}\left(  \frac{1}{\tau_\ell^{3/2}} 1_{\{t/2\geq \tau_\ell\geq 1\}} + \frac{1}{(t-\tau_\ell)^{3/2}}1_{\left\{t/2 < \tau_\ell \leq t-1 \right\}} + 1_{\{\tau_\ell\leq 1\}}\right).
\end{align}

Now we deal with the case $t-\tau_\ell < 1$. For $z = B_{\tau_\ell},$
\begin{align}
	&\mathbf{P}_{z}\left( B_s\leq f_{s+\tau_\ell}^{t,y}, s\leq t-\tau_\ell, B_{t-\tau_\ell} \geq f_t^{t,y}-4\right) \leq  \mathbf{P}_0\left( B_{t-\tau_\ell} \geq f_t^{t,y}-4-z\right) \\
	& = \mathbf{E}_0^{\sqrt{2\lambda^*}}\left(e^{-\sqrt{2\lambda^*}B_{t-\tau_\ell}+\lambda^*(t-\tau_\ell)}1_{\{ B_{t-\tau_\ell} \geq f_t^{t,y}-4-z \}}\right) \leq e^{-\sqrt{2\lambda^*}(f_t^{t,y}-4-z)+\lambda^*(t-\tau_\ell)},
\end{align}
which implies that for all $y,t \geq 1$ and $\ell $ with $t-\tau_\ell <1$,
\begin{align}
	L_2(\ell) &= \mathbf{E}_0\left(e^{-\alpha_0\sqrt{2\lambda^*}\left(f_t^{t,y}-B_{\tau_\ell}\right) +2\alpha_0\lambda^* (t-\tau_\ell)} 1_{\left\{B_s \leq f_s^{t,y},\ s\leq t, B_t \geq f_t^{t,y}-4 \right\}} \right) \\
	& \lesssim   \mathbf{E}_0\left(e^{-(1+\alpha_0)\sqrt{2\lambda^*}\left(f_t^{t,y}-B_{\tau_\ell}\right) +(1+2\alpha_0)\lambda^*(t-\tau_\ell) } 1_{\left\{B_s \leq f_s^{t,y},\ s\leq \tau_\ell \right\}} \right)\\
	&= e^{-(1+\alpha_0)\sqrt{2\lambda^*}f_t^{t,y} +(1+2\alpha_0)\lambda^*(t-\tau_\ell) } \sum_{k=0}^\infty \mathbf{E}_0\left(e^{(1+\alpha_0)\sqrt{2\lambda^*}B_{\tau_\ell}} 1_{\left\{B_s \leq f_s^{t,y},\ s\leq \tau_\ell \right\}} 1_{\left\{ B_{\tau_{\ell}}-f_{\tau_{\ell}}^{t,y}\in[-(k+1), -k] \right\}}\right)   \\
	& \lesssim  e^{-(1+\alpha_0)\sqrt{2\lambda^*}f_t^{t,y} +(1+2\alpha_0)\lambda^*(t-\tau_\ell) }e^{\alpha_0\sqrt{2\lambda^*} f_{\tau_{\ell}}^{t,y}} e^{\lambda^*\tau_\ell}\\
	&\quad \times  \sum_{k=0}^\infty e^{-\alpha_0\sqrt{2\lambda^*}k}\mathbf{P}_0^{\sqrt{2\lambda^*}}\left( B_s \leq f_s^{t,y},\ s\leq \tau_\ell ,B_{\tau_{\ell}}-f_{\tau_{\ell}}^{t,y}\in[-(k+1), -k] \right) \\
	& \lesssim e^{-(1+\alpha_0)\sqrt{2\lambda^*}f_t^{t,y} +(1+2\alpha_0)\lambda^*(t-\tau_\ell) }e^{\alpha_0\sqrt{2\lambda^*} f_{\tau_{\ell}}^{t,y}} e^{\lambda^*\tau_\ell} \sum_{k=0}^\infty e^{-\alpha_0\sqrt{2\lambda^*}k} \frac{(y\land \sqrt{t})(k+1)}{\tau_{\ell}^{3/2}},
\end{align}
where in the last  inequality we used \eqref{step_56}. Therefore, when $t-\tau_{\ell}<1$,
\begin{align}\label{step_58}
	L_2(\ell)&\lesssim e^{-(1+\alpha_0)\sqrt{2\lambda^*}f_t^{t,y} +(1+2\alpha_0)\lambda^*(t-\tau_\ell) }e^{\alpha_0\sqrt{2\lambda^*} f_{\tau_{\ell}}^{t,y}} e^{\lambda^*\tau_\ell} \frac{(y\land \sqrt{t})}{\tau_{\ell}^{3/2}}\nonumber\\
	& \lesssim \left(y\land \sqrt{t}\right)e^{-\lambda^* t}e^{-\sqrt{2\lambda^*}y}.
\end{align}
Using \eqref{step_59}, \eqref{ineq-L_1}, \eqref{step_57} and \eqref{step_58},
taking expectation with respect to $\widehat\P_{(0,i)}$,
we get
\begin{align}
&\widehat{\E}_{(0,i)}\left( \left(\mathcal{A}^{t,y}\right)^{\alpha_0}1_{\Upsilon(\xi)} \right) \lesssim  \widehat{\E}_{(0,i)}\left(L_1 + \sum_{\ell=1}^{K_t} L_2(\ell)\right)\nonumber\\
 &\lesssim (y\land \sqrt{t})e^{-\lambda^* t}e^{-\sqrt{2\lambda^*}y}\widehat{\E}_{(0,i)}\left( 1   +  \sum_{\ell=1}^{K_t} \left(  \frac{1_{\{t/2\geq \tau_\ell\geq 1\}}}{\tau_\ell^{3/2}}  + \frac{1_{\left\{t/2 < \tau_\ell \leq t-1 \right\}} }{(t-\tau_\ell)^{3/2}} + 1_{\{\tau_\ell\leq 1\}}+ 1_{\{\tau_\ell \geq  t-1\}}\right)  \right).\end{align}
We have finished the upper estimate  on $ \widehat{\E}_{(0,i)}\left( \left(\mathcal{A}^{t,y}\right)^{\alpha_0}  1_{\Upsilon(\xi)}	\right)$.

Therefore, by \eqref{step_48}, we have
\begin{align}\label{step_40}
	& \E_{(0,i)}\left( (\mathcal{A}^{t,y}(i'))^{1+\alpha_0} \right) \nonumber\\
	& \lesssim  \frac{(y\land \sqrt{t})}{e^{\sqrt{2\lambda^*}y} }
 \widehat{\E}_{(0,i)}\left( 1   +  \sum_{\ell=1}^{K_t} \left(  \frac{1_{\{t/2\geq \tau_\ell\geq 1\}}}{\tau_\ell^{3/2}}  + \frac{1_{\left\{t/2 < \tau_\ell \leq t-1 \right\}} }{(t-\tau_\ell)^{3/2}} + 1_{\{\tau_\ell\leq 1\}}+ 1_{\{\tau_\ell \geq  t-1\}}\right)  \right).
\end{align}

{\bf Step 2}\quad
Let $0< \underline{a}\le \overline{a}$ be such that
\[
0< \underline{a} \leq \min_{i\in S}(a_i +\lambda^*) \leq \max_{i\in S}(a_i +\lambda^*) \leq \overline{a}.
\]
Recall that
$D_J:= \{t: I_\xi(t-)\neq I_\xi(t) \}.$
We can define two processes
$I_t^{\underline{a}}$ and $I_t^{\overline{a}}$ with the same jumping probability as $I_\xi(t)$ and with constant jump rates
$\underline{a}$ and $\overline{a}$ respectively.
Similarly,  we define $D^{\underline{a}}$ and $D^{\overline{a}}$
to be the jumping times of $I_t^{\underline{a}}$ and $I_t^{\overline{a}}$.
We can construct a coupling of $(I_\xi(t), I_t^{\underline{a}}, I_t^{\overline{a}})$ such that
the embedded chain of the three processes are the same and
the jump times
\begin{align}\label{Coupling4}
\!\!\!\!\!D_J= \{t_n: 0< t_1< t_2<... \},\,\,  D^{\underline{a}}= \{t_n^{\underline{a}}: 0< t_1^{\underline{a}} < t_2^{\underline{a}} <... \},\,\, D^{\overline{a}} =\{t_n^{\overline{a}}: 0< t_1^{\overline{a}} < t_2^{\overline{a}}<... \}
\end{align}
satisfy $t_n^{\overline{a}}\leq t_n \leq t_n^{\underline{a}}$ for every $n$.
More precisely, let $\{Y_n: n=0, 1,...\}$ be the embedded chain of $I_\xi(t)$ with $Y_0 = I_
\xi(0)$. Let $T_0:= t_1$ and $T_n:= t_{n+1}-t_{n}$ for $n\geq 1$. Then by the strong Markov property, given $Y_j, 0\leq j\leq n$, $T_n$ is an exponential distribution with parameter $a_{Y_n}+\lambda^*$. Let $T_n^{\overline{a}}:= (a_{Y_n}+\lambda^*)T_n/\overline{a} \leq T_n \leq (a_{Y_n}+\lambda^*)T_n/\underline{a}=: T_n^{\underline{a}}$. Then we see that given $Y_j, 0\leq j\leq n$, $T_n^{\overline{a}}$ and $T_n^{\underline{a}}$ are exponential distribution with parameter $\overline{a}$ and $\underline{a}$ respectively. Now for $n\geq 1$, define $t_n^{\overline{a}}:= \sum_{j=0}^{n-1}T_j^{\overline{a}}$ and $t_n^{\underline{a}}:=  \sum_{j=0}^{n-1} T_j^{\underline{a}}$.
Define $K_t^{\overline{a}}:= \sup\{n:  t_{n}^{\overline{a}}\leq  t\}$, $K_t^{\underline{a}}:= \sup\{n: t_{n}^{\underline{a}}\leq t\}$,
$I_t^{\overline{a}}:= Y_{K_t^{\overline{a}}}$ and $I_t^{\underline{a}}:= Y_{K_t^{\underline{a}}}$. Then $(I_\xi(t), I_t^{\underline{a}}, I_t^{\overline{a}})$ is the desired coupling.
Therefore, for any non-negative  and non-increasing  function $f$,
\begin{align}\label{Coupling1}
	\sum_{s\in D_J: s\leq t} f(s)= \sum_{n=1}^{K_t} f(t_n)\leq \sum_{n=1}^{K_t} f(t_n^{\overline{a}})\leq \sum_{i=1}^{K_t^{\overline{a}}} f(t_n^{\overline{a}}) = \sum_{s\in D^{\overline{a}}: s\leq t} f(s),
\end{align}
here $K_t\leq \# \left\{ D^{\overline{a}}\cap [0,t] \right\}= K_t^{\overline{a}}$ by the coupling.
Applying \eqref{Coupling1} to
$f(s)= 1_{\{s\leq 1\}} + s^{-3/2}1_{\{s\in [1, t/2]\}}$ and
$f(s)= 1_{\{s\leq 1\}}$, by the Markov property at time $t-1$,
 we get that
\begin{align}\label{step_41}
	& \widehat{\E}_{(0,i)}\left( 1   +  \sum_{\ell=1}^{K_t} \left(  \frac{1_{\{t/2\geq \tau_\ell\geq 1\}}}{\tau_\ell^{3/2}} +1_{\left\{\tau_\ell \leq 1 \right\}} + 1_{\{\tau_\ell \geq  t-1\}}\right)  \right)\nonumber\\
	&\leq \widehat{\E}_{(0,i)}\left( 1   +  \sum_{\ell=1}^{K_t} \left(  \frac{1_{\{t/2\geq \tau_\ell\geq 1\}}}{\tau_\ell^{3/2}} +1_{\left\{\tau_\ell \leq 1 \right\}} \right) \right) +\sup_{j\in S} \widehat{\E}_{(0,j)}\left(\sum_{\ell=1}^{K_1} 1\right)\nonumber \\
	& \leq 1+ \E \left(\sum_{\ell=1}^{K_t^{\overline{a}}} \frac{1}{(\tau_\ell^{\overline{a}})^{3/2}} 1_{\{t/2\geq  \tau_\ell^{\overline{a}}\geq 1\}} + 1_{\{\tau_\ell^{\overline{a}}\leq 1\}} \right) + \mathbb{E}\left(\sum_{\ell=1}^{K_1^{\overline{a}}} 1\right)\nonumber \\
	& = 1+ \overline{a} \int_0^t\left(1_{\{s\leq 1\}} +\frac{1}{s^{3/2}} 1_{\{s\in (1, t-1]\}}  \right)\mathrm{d} s + \overline{a} \int_0^1 1 \mathrm{d}s \lesssim 1.
\end{align}
For $\tau_\ell \in (t/2, t-1]$,  note that
\[
\sum_{\ell=1}^{K_t}  \frac{1_{\{t-1 \geq  \tau_\ell >t/2\}}}{(t-\tau_\ell)^{3/2}} \leq \sum_{k=[t/2]}^{\lceil t-1\rceil-1}\frac{1_{\{k< \tau_\ell \leq (k+1)\}}}{(t-k-1)^{3/2}},
\]
where $\lceil t-1\rceil$ is the smallest integer larger than or equal to $t-1$. For each $k$, using  the Markov property at time $k$, and using \eqref{Coupling1}, we have
\begin{align}\label{step_51}
	&\widehat{\E}_{(0,i)}\left(\sum_{\ell=1}^{K_t}  \frac{1_{\{t-1 \geq  \tau_\ell >t/2\}}}{(t-\tau_\ell)^{3/2}} \right)  \leq  \sum_{k=[t/2]}^{\lceil t-1\rceil-1} \frac{1}{(t-k-1)^{3/2}} \sup_{j\in S}\widehat{\E}_{(0,j)} \left(\sum_{\ell=1}^{K_1}1 \right)\nonumber\\
	& \leq \sum_{k=[t/2]}^{\lceil t-1\rceil-1} \frac{1}{(t-k-1)^{3/2}} \E \left(\sum_{\ell=1}^{K_1^{\overline{a}}}1 \right) \lesssim  \sum_{k=[t/2]}^{\lceil t-1\rceil-1} \frac{1}{(t-k-1)^{3/2}} \lesssim 1.
\end{align}
Combining \eqref{step_40}, \eqref{step_41} and \eqref{step_51},  we get \eqref{toprove2}.
The proof is now complete.
\hfill$\Box$

\begin{remark}
	As a consequence of \eqref{Coupling1}, we have the following useful inequality:
	for any $r<t$ and any decreasing non-negative function $f$ on $[r,t]$,
	\begin{align}\label{Coupling2}
		&\widehat{\E}_{(0,i)}\left(\exp\left\{ -\sum_{s\in D_J: r< s\leq t} f(s)\right\} \right)\geq \inf_{\ell \in S}\widehat{\E}_{(0,\ell)}\left(\exp\left\{ -\sum_{s\in D_J:  s\leq t-r} f(s+r)\right\} \right)\nonumber \\
		& \geq \inf_{\ell \in S}\widehat{\E}_{(0,\ell)}\left(\exp\left\{ -\sum_{s\in D^{\overline{a}}:  s\leq t-r} f(s+r)\right\} \right) = \exp\left\{ -\overline{a}\int_0^{t-r}\left(1- e^{-f(s+r)}\right)\mathrm{d}s \right\}.
	\end{align}
Using the fact that
$\# \left\{ D_J \cap [0,t]\right \} \geq \# \left\{ D^{\underline{a}} \cap [0,t]\right\} $, we also have that for any $\theta >0$,
	\begin{align}\label{Coupling3}
		&\widehat{\E}_{(0,i)}\left(\exp\left\{ -\theta \sum_{s\in D_J: r< s\leq t} 1\right\} \right) \leq  \sup_{\ell \in S}\widehat{\E}_{(0,\ell)}\left(\exp\left\{ -\theta\sum_{s\in D^{\underline{a}}:  s\leq t-r} 1\right\} \right)
		\\ &= \exp\left\{- \underline{a}\int_0^{t-r} \left(1- e^{-\theta}\right)\mathrm{d}s\right\}= e^{-\underline{a}(t-r)(1-e^{-\theta})}.
	\end{align}
\end{remark}

Using Propositions \ref{Upper-Bound-M-T} and \ref{Lower-Bound-M-T}, we can get the following result:
\begin{thrm}\label{Tightness-2}
	For any $i\in S$, $\left(M_t - m(t),\ t\geq1, \P_{(0,i)}\right) $
	is tight. Also, $\left(M_t^{i'} - m(t),\ t\geq1, \P_{(0,i)}\right) $ is tight for any $i,i' \in S$.
\end{thrm}
\textbf{Proof: }
Fix $i \in S$. For any $\varepsilon>0$,
choose $y>4$ and $\delta$ small so that $K_0 ye^{-\sqrt{2\lambda^*}y}\geq \delta$, where $K_0$ is the constant in Proposition \ref{Lower-Bound-M-T}.
Now choose $L$ large so that $\E_{(0,i)}\left((1-\delta)^{\langle \mathbf{N}_L , \mathbf{1}\rangle} \right)< \varepsilon/2.$
Indeed, we can find a large $n$ such that $(1-\delta)^n < \varepsilon/4$, therefore,
\[
\E_{(0,i)}\left((1-\delta)^{\langle \mathbf{N}_L , \mathbf{1}\rangle } \right)<\frac{ \varepsilon}{4}+ \P_{(0,i)}\left(\langle \mathbf{N}_L , \mathbf{1}\rangle < n\right),
\]
which is less than $\varepsilon/2$ for large $L$ since $\langle \mathbf{N}_L , \mathbf{1}\rangle \to +\infty\ \P_{(0,i)}$-a.s.
Let $b>0$ be a constant such that
\[
\P_{(0,i)}\left(\min_{u\in Z(L)}X_u(L) < -b \right)<
\frac{\varepsilon}{2}.
\]
By Proposition \ref{Lower-Bound-M-T}, for $t$ large enough so that
$t-L > \max\{\beta_0^2, y^2\}$ where $\beta_0$ is the constant in Proposition \ref{Lower-Bound-M-T},
\[
\inf_{i, i'\in S} \P_{(0,i)} \left(M_{t-L}^{i'}> m(t-L) +y\right)
\geq K_0 y e^{-\sqrt{2\lambda^*}y}\geq \delta.
\]
Thus,
\begin{align}
	&\P_{(0,i)}\left(M_t< m(t-L)-b +y\right)\leq \P_{(0,i)}\left(M_t^{i'}< m(t-L)-b +y\right)\\
	& \leq \P_{(0,i)}\left(\min_{u\in Z(L)} X_u(L)< -b\right) + \P_{(0,i)}\left(\min_{u\in Z(L)}X_u(L) \geq -b,  \max_{u\in Z(L)} M_{t-L}^{i'}(u)< m(t-L)-b +y \right)\\
	& \leq
	\frac{\varepsilon}{2}
	+ \E_{(0,i)}\left(\prod_{u\in Z(L)}
	 \P_{(0,I_u(L))}\left(M_{t-L}^{i'}< m(t-L)+y\right) \right)
	\leq \frac{\varepsilon}{2}
	+ \E_{(0,i)}\left((1-\delta)^{\langle \mathbf{N}_L , \mathbf{1}\rangle} \right)< \varepsilon.
\end{align}
From this one can easily see that there exists $\widetilde{y}>0$ such that for all $t$ large
$$
\P_{(0,i)}\left(M_t< m(t)-\widetilde{y}\right)\leq \P_{(0,i)}\left(M_t^{i'}< m(t)-
\widetilde{y}\right)< \varepsilon.
$$
Also, by  Proposition \ref{Upper-Bound-M-T},
there exists $y^*$ large enough such that for all $t\ge 1$
\begin{align}\label{e:newthm4.5}
 \P_{(0,i)}\left(M_t^{i'}> m(t) +y^*\right)
 \leq \P_{(0,i)}\left(M_t> m(t) +y^*\right)\leq C_0 y^* e^{-\sqrt{2\lambda^*}y^*}< \varepsilon.
\end{align}
Thus there exists $T>1$ such that $\left(M_t - m(t),\ t\geq T, \P_{(0,i)}\right) $
and $\left(M_t^{i'} - m(t),\ t\geq T, \P_{(0,i)}\right) $ are tight.
Since $ \min_{1\leq t\leq T} M_t^{i'}$ is finite $\P_{(0,i)}$-a.s., for $y^*$ large enough we
also have for $1\leq t\leq T$,
\begin{align}
	&\P_{(0,i)}\left(M_t> m(t) -y^*\right)\geq	\P_{(0,i)}\left(M_t^{i'}> m(t) -y^*\right)\\
	&\geq  \P_{(0,i)}\left(\min_{1\leq t\leq T}	\left( M_t^{i'}-m(t)\right)>  -y^*\right)
	> 1-\varepsilon.
\end{align}
Combining this with \eqref{e:newthm4.5}, we get
\begin{align*}
\mathbb{P}_{(0,i)}\left( -y^* < M_t -m(t) < y^* \right)= 	\mathbb{P}_{(0,i)}\left(  M_t -m(t) > -y^* \right)- 	\mathbb{P}_{(0,i)}\left(  M_t -m(t) > y^* \right)> 1-2\varepsilon
\end{align*}
and
\begin{align*}
\mathbb{P}_{(0,i)}\left( -y^* < M^{i'}_t -m(t) < y^* \right)= 	\mathbb{P}_{(0,i)}\left(  M^{i'}_t -m(t) > -y^* \right)- 	\mathbb{P}_{(0,i)}\left(  M^{i'}_t -m(t) > y^* \right)> 1-2\varepsilon.
\end{align*}
This completes the proof.
\hfill$\Box$

\section{Proof of Theorem \ref{thm1}}

\subsection{Upper bound for $\mathbf{v}$}\label{ss:5.1}

In this subsection,  we first give some estimates involving $R(t;v)$ defined  in \eqref{Non-Linear-R-2} and then use these estimates to get some upper bound for solutions of  \eqref{F-KPP2} with initial condition satisfying \eqref{Initial-Assumption}.
We roughly follow the arguments  of \cite[Sections 6-8]{Bramson1}.  However, some of the arguments in \cite{Bramson1} do not work in the multi-type case. We will explain these later in this section.

Let $\delta\in (0,1/2)$ and $K>2+ \sqrt{2\lambda^*}$ be fixed constants.
If  $L$  is a function on $[0,t]$, for $t>4r>0$, as in \cite[(6.11), p. 88]{Bramson1},   we define
\begin{align}
	\theta_{r,t}\circ L(s) :=\left\{\begin{array}{ll} L\left(s+s^\delta \land (t-s)^\delta \right)+ Ks^\delta \land (t-s)^\delta,\quad &r\leq s \leq t-2r,\\
		\displaystyle L(s),\quad &\mbox{otherwise}.
	\end{array}\right.
\end{align}
We  define $\theta_{r,t}^{-1}$  to be the inverse of  $\theta_{r,t}$.
Similar to \cite[(7.6), p. 99, (6.13), p. 88, and  (6.14), p. 89]{Bramson1}, we define
\begin{align}
	L_{r,t}(s)&:= m(s) -\frac{s}{t}m(t)-\frac{t-s}{t}\alpha(r), \\
	\underline{L}_{r,t} (s)&:= \theta_{r,t}^{-1}\circ L_{r,t}(s) ,\\
	\overline{L}_{r,t}(s)&:= \left(\theta_{r,t} \circ L_{r,t}(s)\right)\vee \underline{L}_{r,t} (s) \vee L_{r,t} (s),
\end{align}
and similar to \cite[(7.44), p. 111]{Bramson1} or \cite[(2.11)]{RSZ}, we define
\begin{align}
	\underline{\mathcal{M}}_{r,t}'(s):= 	\left\{\begin{array}{ll} \displaystyle\underline{L}_{r,t}(s) + \frac{s}{t}m(t) +\frac{t-s}{t}\alpha(r),\quad
		 &r+r^\delta \leq s \leq t-2r,\\
		\displaystyle -\infty ,\quad
		&\mbox{otherwise},
	\end{array}\right.
\end{align}
where the function $\alpha(r)$ is either taken to be $-\log r$ or taken to be identically 0.
Let
\[
    \mathcal{B}_{up}:= \left\{X_s> \underline{\mathcal{M}}_{r,t}'(t-s) \mbox{ for all } s\in [0,t-r]   \right\}.
\]
Note that
when $s\in [0, 2r)\cup [t-r-r^\delta, t]$, $\underline{\mathcal{M}}_{r,t}'(t-s) = -\infty,$ therefore,
\[
(\mathcal{B}_{up})^c=
    \left\{ \exists\ s\in [2r, t-r-r^\delta] \mbox{ such that } X_s \leq \underline{\mathcal{M}}_{r,t}'(t-s) \right\}.
\]
Similar to \cite[(7.19)--(7.20), p.102--103]{Bramson1}, we define
\begin{align}
	S^1(r,t)&:= \sup \left\{s: 2r \leq s\leq t/2,\ X_s \leq \underline{\mathcal{M}}_{r,t}'(t-s) \right\},\\
	S^2(r,t)&:= \inf \left\{s: t/2 \leq s\leq t-r,\ X_s \leq \underline{\mathcal{M}}_{r,t}'(t-s) \right\},\\
	S(r,t)&:= S^1(r,t)1_{\left\{ S^1(r,t)+S^2(r,t)> t\right\}}+S^2(r,t)1_{\left\{ S^1(r,t)+S^2(r,t)\leq  t\right\}}.
\end{align}
We use the convention that $S^1(r,t) =0$ if  $X_s> \underline{\mathcal{M}}_{r,t}'(t-s)$ for all $s\in [2r, t/2]$ and that
$S^2(r,t)=t$ if $X_s> \underline{\mathcal{M}}_{r,t}'(t-s)$ for all $s\in [t/2, t-r]$.
Next, similar to \cite[(7.21), p.103]{Bramson1}, we define for $r_1 \in [r, t/2]$,
\[
(\mathcal{B}_{up}^{r_1})^c:= \left\{r_1\leq S(r,t)\leq t-r_1 \right\}
 =\left\{\exists\ s\in[r_1\vee (2r), t-r_1] \mbox{ such that }  X_s \leq \underline{\mathcal{M}}_{r,t}'(t-s) \right\}.
\]
Let $j_1$ be the integer such that $j_1< t/2 \leq j_1 +1$. Define $G_j = [j, j+1)\cup (t-j-1, t-j], j=0,\cdots, j_1-1$
and $G_{j_1}= [j_1, t-j_1]$.
Similar to \cite[(7.23)--(7.24), p. 103]{Bramson1},
we define for $j=0,\cdots,j_1$,
\begin{align}
	A_j(r,t)&:= \left\{S(r,t)\in G_j \right\},
	\\
	A_j^1(r,t)&:= A_j(r,t) \cap \left\{X_s> -(s\land (t-s))+ \frac{s}{t}y +\frac{t-s}{t}x
	\mbox{ for all }s\in G_j \right\}
\end{align}
and $A_j^2(r,t)= A_j(r,t)\setminus A_j^1(r,t).$
The following result is \cite[Lemma 7.1, p. 104]{Bramson1}.
\begin{lemma}\label{lemma4}
	For large $r$, $t>4r$ and $r_1\in [r, t/2]$, for any $y\geq -\log r$ and $x\geq m(t)$,
	\[
	\mathbf{P}_x\left(\mathcal{B}_{up}^{r_1} \big| X_t=y \right)\lesssim \frac{r_1}{r}\mathbf{P}_x\left(\mathcal{B}_{up}\big| X_t=y\right).
	\]
\end{lemma}

Recall that by the definitions \eqref{Non-Linear-R}--\eqref{Non-Linear-R-2},
\begin{align}
	\exp\left\{\sum_{s\in D_J, s\leq t-r} \log \left(\frac{ \varphi_{I_s}\left(\mathbf{v}(t-s, X_s)\right)}{\sum_{j=1}^d m_{I_s, j}v_j (t-s, X_s)}\right) \right\} = R_t(t-r;v).
\end{align}
If $\mathbf{v}$ solves \eqref{F-KPP2}  with  initial value satisfying \eqref{Initial-Assumption}, then by \eqref{Prob-Representation},
for any $i\in S$,
\begin{align}\label{Lower-Bound-M-T2}
	v_i(t,x)&\geq 1- \E_{(x,i)}\left(\prod_{u \in Z(t): I_u(t)=i_0}\left(1- 1_{(-\infty, N_1)}(X_u(t))\right) \right)\nonumber \\
&= 1- \E_{(0,i)}\left(\prod_{u \in Z(t): I_u(t)=i_0}1_{\{X_u(t)+x \geq N_1\}} \right) = 1- \E_{(0,i)}\left(\prod_{u \in Z(t): I_u(t)=i_0}1_{\{X_u(t) \leq x- N_1\}} \right) \nonumber\\ & = \P_{(0,i)}\left(M_t^{i_0} > x-N_1\right).
\end{align}
By induction, for any $0\leq x_k ,y_k\leq 1$ with $x_k +y_k\leq 1$ for all $1\leq k\leq n$, it holds that
\begin{align}\label{step_60}
1- \prod_{k=1}^n (1-x_k-y_k)\leq \left\{1- \prod_{k=1}^n (1-x_k) \right\} +\left\{1- \prod_{k=1}^n (1-y_k) \right\}.
\end{align}
Indeed, it is easy to see that \eqref{step_60} holds for $n=1$.
 If \eqref{step_60}  holds for $n$ and $0\leq x_k, y_k, x_k +y_k\leq 1$ with $1\leq k\leq n$, then
 \begin{align*}
 	&1- \prod_{k=1}^{n+1} (1-x_k-y_k) = (1-x_{n+1}-y_{n+1})\left(1-\prod_{k=1}^n (1-x_k-y_k) \right) + x_{n+1}+ y_{n+1}\\
 	& \leq (1-x_{n+1}-y_{n+1})\left(\left\{1- \prod_{k=1}^n (1-x_k) \right\} +\left\{1- \prod_{k=1}^n (1-y_k) \right \} \right)+ x_{n+1}+y_{n+1}\\
 	& \leq (1-x_{n+1})\left\{1- \prod_{k=1}^n (1-x_k) \right\}  + (1-y_{n+1})\left\{1- \prod_{k=1}^n (1-y_k) \right\} +x_{n+1 } + y_{n+1}  \\
 	& = \left\{1- \prod_{k=1}^{n+1} (1-x_k) \right\} +\left\{1- \prod_{k=1}^{n+1} (1-y_k) \right\},
 \end{align*}
which implies \eqref{step_60}.

Suppose that  $\mathbf{v}$ satisfies \eqref{F-KPP2} and \eqref{Initial-Assumption}. For $r> 1$ and $t\geq r$, define $v_j^*(0,x):=v_j(r,x)1_{\{x > -\log r\}}$ for all $j\in S$ and
\begin{align}\label{V-Star}
	v^*_i(t-r,x) = 1- \mathbb{E}_{(x, i)}\left(\prod_{u\in Z(t-r)} \left(1- v_{I_u(t-r)}^*(0, X_u(t-r))\right)\right).
\end{align}

The next lemma is slightly  different from \cite[Proposition 8.3 (b), p.  136]{Bramson1}.
In \cite{Bramson1} (see the argument \cite[from (8.44) to (8.46), p.137]{Bramson1}), Bramson used the
Feynman-Kac formula to get
\begin{align*}
	v(t, m(t)+ x) = \int_\mathbb{R} v(r,y)\frac{e^{-\frac{(m(t)+x-y)^2}{2(t-r)}}}{\sqrt{2\pi (t-r)}} \mathbf{E}_{m(t)+x}\left(\exp\left\{\int_0^{t-r}k(v(t-r-s, X_s))\mathrm{d}s\right\}\Big| X_{t-r}=y\right)\mathrm{d}y,
\end{align*}
where $k(x)=f(x)/x$.
Then he separated the integral into $\int_{-\infty}^{-\log r} + \int_{-\log r}^\infty$. For the $\int_{-\infty}^{-\log r}$ part (see  \cite[(8.45), p. 137]{Bramson1} ), he used  the fact that
 if $f:[0,1] \to [0,\infty)$ is a concave function with $f(0)=0$
 (which is obviously the case for $f(x) = 1- \sum_{k=0}^\infty p_k (1-x)^k$ with $\{p_k: k\ge 1\}$
a distribution with finite mean), then
$k(x)=f(x)/x$
 is decreasing in $(0,1]$.
But when $f$ is a multivariable function like $\varphi_i$, it no longer holds that
\[
\frac{f(\mathbf{v})}{\triangledown f(0) \cdot \mathbf{v}}\leq \frac{f(\widetilde{\mathbf{v}})}{\triangledown f(0) \cdot \widetilde{\mathbf{v}}},\quad\quad \mbox{if}\ 1\geq v_j\geq \widetilde{v}_j\geq 0\ \mbox{for all}\ j\in S.
\]
To avoid this difficulty, we deal with the part $\int_{-\infty}^{-\log r}$ by first using probabilistic representation \eqref{Prob-Representation},  Propositions \ref{Upper-Bound-M-T} and \ref{Lower-Bound-M-T},
and then using the Feynman-Kac formula.
This is accomplished in the following lemma.

\begin{lemma}\label{lemma2}
	 Suppose that $\mathbf{v}$ satisfies \eqref{F-KPP2} and \eqref{Initial-Assumption}.
	 If $r$ is large enough,
         then  there exists a positive function $C(r)$
	 with $\lim_{r\to\infty}C(r)=1$ such that for all $r< x \leq  \sqrt{t}$,
	\begin{align}
		&v_i (t,m(t)+ x)
		\leq C(r) v_i^* \left(t-r, m(t)+x\right)\\
		& = C(r)e^{\lambda^* (t-r)}h_i\int_{-\log r} ^\infty \frac{e^{-\frac{(m(t)+x-y)^2}{2(t-r)}}}{\sqrt{2\pi (t-r)}}  \mathbf{E}_{(m(t)+x,i)}^h\left(R(t-r;v^*) \frac{v_{I_{t-r}}(r,y)}{h_{I_{t-r}}}\bigg| X_{t-r}=y \right)
		\mathrm{d} y.
	\end{align}
\end{lemma}
\textbf{Proof: }  The equality in the lemma follows from
the Feynman-Kac formula \eqref{Feynman-Kac-3}.
Note that $1\geq v_j^*(0,x) + 1_{\left\{x \leq -\log r \right\}} \geq v_j(r,x)$, by \eqref{Prob-Representation} and \eqref{step_60},
we have
\begin{align}\label{step_64}
	&v_i(t,m(t)+x)\leq v_i^*(t-r, m(t)+ x) + 1- \mathbb{E}_{(m(t)+x, i)}\left(\prod_{u\in Z(t-r)} \left(1- 1_{\{X_u(t-r)\leq -\log r \}}\right)\right)\nonumber\\
	& = v_i^*(t-r, m(t)+ x)+\P_{(0,i)} \left(M_{t-r}> m(t)+ x +\log r\right).
\end{align}
By \eqref{Lower-Bound-M-T2},  Propositions \ref{Upper-Bound-M-T} and \ref{Lower-Bound-M-T},
for all  $r<x \leq  \sqrt{t}$, if $r$ is large enough so that $ r-N_1 \geq \frac{r}{2}$ and $\sqrt{2\lambda^*}r -1 \leq m(t)-m(t-r)\leq \sqrt{2\lambda^*}r$,  then
\begin{align}
	&\frac{ \P_{(0,i)} \left(M_{t-r}> m(t)+x +\log r\right)}{v_i(t,m(t)+x)}\leq \frac{ \P_{(0,i)} \left(M_{t-r}> m(t)+x +\log r\right)}{\P_{(0,i)}\left(M_{t}^{i_0}>m(t)+ x- N_1\right)}
	\\ & \lesssim \frac{ (m(t)-m(t-r)+x+ \log r) e^{-\sqrt{2\lambda^*}(m(t)-m(t-r)+x+\log r)}}{ \left(x-N_1 \right) e^{-\sqrt{2\lambda^*}\left(x-N_1 \right) } } \\ & \lesssim \frac{ (\sqrt{2\lambda^*} r+x+ \log r) e^{-\sqrt{2\lambda^*}(\sqrt{2\lambda^*}r+\log r)}}{ \left(x-N_1 \right)  } \lesssim e^{-\sqrt{2\lambda^*}(\sqrt{2\lambda^*}r+\log r)},
\end{align}
 where in the  last inequality we used the fact that for $x>r$ and $r-N_1 \geq r/2$,
\begin{align}
	\frac{ (\sqrt{2\lambda^*} r+x+ \log r) }{ \left(x-N_1 \right)  } \leq \frac{(\sqrt{2\lambda^*}+1+ 1)x}{x/2}\lesssim 1.
\end{align}
Therefore,  if $\Gamma$ is a constant such that  for large $r$ and $r<x\leq  \sqrt{t}$,
\[
\frac{ \P_{(0,i)} \left(M_{t-r}> m(t)+x +\log r\right)}{v_i(t,m(t)+x)} \leq \Gamma e^{-\sqrt{2\lambda^*}(\sqrt{2\lambda^*}r+\log r)}<1,
\]
we can choose $C(r)$ to be
\begin{align}\label{C-R}
         C(r):=
	\left(1-  \Gamma  e^{-\sqrt{2\lambda^*}(\sqrt{2\lambda^*}r+\log r)}\right)^{-1}.
\end{align}
\hfill$\Box$

Next, we will prove a lemma (Lemma \ref{lemma3}) similar to \cite[Lemma 7.2, p. 105]{Bramson1}.
A key  step in the proof of Lemma \ref{lemma3} is the inequality \eqref{step_26}.
In Bramson's argument for the analog of \eqref{step_26} (see \cite[(7.32) on p.  107, Proposition 7.1 on p. 97]{Bramson1}),  the Kolmogorov-Petrovsky-Piscounov theorem (see
\cite[p. 34]{Bramson1})  was used, see \cite[Proposition 3.4 on p. 47 and (3.71) on p. 49]{Bramson1}.
In the multi-type case, the analog of the Kolmogorov-Petrovsky-Piscounov theorem has not been proved yet. So we have to overcome this difficulty.
Lemma \ref{lemma5} below is the key to \eqref{step_26},  which is different from
\cite[(7.32)]{Bramson1}.
Roughly speaking, since $v_i^*(0,y)=v_i(r,y)1_{\{ y> -\log r\}}$ is very close to $1$ when $|y|< \log r$,
by representation \eqref{Prob-Representation}
for $v_i^*(t,y)$, it suffices to show that, under  $\mathbb{P}_{(y,i)}$, the probability of the event that there is at least one particle locating in $[-\log r, \log r]$ is close to $1$ when $r$ is large enough.
This is easy to prove
since we know the  behavior of the maximal position $M_t$ very well by Theorem \ref{Tightness-2}.
Using this, we can get Lemma \ref{lemma5} below.

Recall the definition \eqref{V-Star} of $\mathbf{v}^*$,  where $\mathbf{v}$ is a solution to \eqref{F-KPP2} with initial value satisfying \eqref{Initial-Assumption}. Define
\begin{align}
			m_+(t):=\max\left\{ \sqrt{2\lambda^*} t -\frac{3}{2\sqrt{2\lambda^*}} \log_+ t, 0\right\}.
\end{align}

\begin{lemma}\label{lemma5}
         For any $\varepsilon>0$, 	there exists $N=N(\varepsilon)$ such that when $r> N$,
	\begin{align}
		v_i^*(t, y)\geq 1-\varepsilon,\quad \mbox{for all}\ t\geq 0,\ i\in S,\ y\in [0, m_+(t)].
	\end{align}
\end{lemma}
\textbf{Proof: }  We first prove that
for any $\varepsilon_1>0$, there exists $N^*= N^*(\varepsilon_1)$ such that when $r\geq N^*$,
\begin{align}\label{step_61}
	\mathbb{P}_{(0,i)}\left(\# \left\{u\in Z(t):\ |X_u(t)| \leq r \right\} = 0\right)<\varepsilon_1,\quad i\in S, t\geq 0.
\end{align}
Using Theorem \ref{Tightness-2} and symmetry, we get that, for any $\varepsilon_1>0$,  there exists $N_1$ such that
\begin{align}\label{step_63}
\sup_{t\geq 0} \sup_{i\in S} \mathbb{P}_{(0,i)}\left(|M_t - m_+(t)| > N_1\right) =\sup_{t\geq 0} \sup_{i\in S} \mathbb{P}_{(0,i)}\left(|M_t^-+m_+(t)| > N_1\right)  < \frac{\varepsilon_1}{2}.
\end{align}
Here $M_t^-:= \inf_{u\in Z(t)} X_u(t)$ is the leftmost position among all the particles. By the Markov property and branching property at time $t/2$, we have
\begin{align*}
	&\mathbb{P}_{(0,i)}\left(\# \left\{u\in Z(t):\ |X_u(t)| \leq 2 N_1\right\} = 0\right)\\
	& \leq \frac{\varepsilon_1}{2}+ \mathbb{P}_{(0,i)} \left(|M_{t/2}-m_+(t/2)| \leq N_1, \# \left\{u\in Z(t):\ |X_u(t)| \leq 2 N_1\right\} =0\right)\\
	& \leq  \frac{\varepsilon_1}{2} + \sup_{z: | z-m_+(t/2)| \leq N_1} \sup_{j\in S} \mathbb{P}_{(z,j)}\left(|M_{t/2}^-| > 2N_1\right)\\
	& \leq \frac{\varepsilon_1}{2} + \sup_{z: | z-m_+(t/2)| \leq N_1} \sup_{j\in S} \mathbb{P}_{(z,j)}\left(|M_{t/2}^-+ m_+(t/2)| > 2N_1- | z-m_+(t/2)| \right)\\
	& \leq \frac{\varepsilon_1}{2} +\sup_{j\in S} \mathbb{P}_{(0,j)}\left(|M_{t/2}^-+ m_+(t/2)| > N_1 \right)<\varepsilon_1.
\end{align*}
Therefore, \eqref{step_61} holds with $N^*= 2N_1$.

Next, we prove that, for any $\varepsilon_2>0$, there exists $N'= N'(\varepsilon_2) $ such that when $r\geq N'$,
\begin{align}\label{step_62}
	\mathbb{P}_{(y,i)}\left(\# \left\{u\in Z(t):\ |X_u(t)| \leq r \right\} = 0\right)<\varepsilon_2,\quad i\in S, t\geq 0, y\in [0, m_+(t)].
\end{align}
Let $t_0>1$ be a constant such that $m_+(t)= m(t)$ for all $t\geq t_0$. When $t\leq t_0$, we use
the trivial upper-bound
\[
\mathbb{P}_{(y,i)}\left(\# \left\{u\in Z(t):\ |X_u(t)| \leq r \right\} = 0\right) \leq \mathbf{P}_y \left(|B_t| > r\right)
\]
and the tail probability of normal random variables; when
$y\leq m(t_0)$, we use
the bound
\[
\mathbb{P}_{(y,i)}\left(\# \left\{u\in Z(t):\ |X_u(t)| \leq r \right\} = 0\right) \leq  \mathbb{P}_{(0,i)}\left(\# \left\{u\in Z(t):\ |X_u(t)| \leq r -m(t_0)\right\} = 0\right)
\]
and \eqref{step_61}. So we only deal with the case when $t\geq t_0$ and $y\in [m(t_0), m_+(t)]= [m(t_0), m(t)]$.  Suppose that $y= m(s)$ for some
$t_0\leq s \leq t.$
 Let $\varepsilon_1 = \varepsilon_2/2$,
 Using \eqref{step_61}, \eqref{step_63} and the fact that $\mathbb{P}_{(0,i)}\left(|M_s^-+m_+(s)| > N_1\right) = \mathbb{P}_{(y,i)}\left(|M_s^-| > N_1\right) $,
\begin{align}
	&\mathbb{P}_{(y,i)}\left(\# \left\{u\in Z(t):\ |X_u(t)| \leq N_1+ N^*\right\} = 0\right) \\
	& \leq \frac{\varepsilon_1}{2} +
	\mathbb{P}_{(m(s),i)}\left(|M_s^-| \leq N_1, \# \left\{u\in Z(t-s):\ |X_u(t-s)| \leq N_1+ N^*\right\} = 0\right) \\
	& \leq \frac{\varepsilon_2}{4} +
	 \sup_{z: |z|\leq N_1} \sup_{j\in S} \mathbb{P}_{(z,j)}\left(\# \left\{u\in Z(t-s):\ |X_u(t-s)| \leq N_1+ N^*\right\} = 0\right)\\
	& \leq \frac{\varepsilon_2}{4} +
	\sup_{j\in S} \mathbb{P}_{(0,j)}\left(\# \left\{u\in Z(t-s):\ |X_u(t-s)| \leq N^*\right\} = 0\right)
	< \frac{\varepsilon_2}{4} + \frac{\varepsilon_2}{2}< \varepsilon_2,
\end{align}
which implies \eqref{step_62}.

For any $\varepsilon>0$,  by \eqref{Lower-Bound-M-T2} and Theorem \ref{Tightness-2}, when $r$ is large enough,  we have for all $x\in [-\log r, \log r]$ and any $i\in S$,
\[
v_i(r, x)\geq \P_{(0,i)}\left(M_r^{i_0} > x-N_1\right)\geq \P_{(0,i)}\left(M_r^{i_0}-m(r)> \log r - m(r)-N_1\right)> 1-\frac{\varepsilon}{2}.
\]
Taking $\varepsilon_2 =\frac{\varepsilon}{2}$ in \eqref{step_62}, we get that for any $i\in S, t\geq 0$ and $y\in [0, m_+(t)]$, as long as $\log r> N' \Leftrightarrow r > e^{N'}=:N$,
\begin{align*}
	v_i^*(t,y)& = 1- \mathbb{E}_{(y, i)}\left(\prod_{u\in Z(t)} \left(1- v_{I_u(t)}^*(0, X_u(t))\right)\right)\\
	& \geq 1- \mathbb{E}_{(y, i)}\left(\prod_{u\in Z(t)} \left(1- v_{I_u(t)}(r, X_u(t))1_{\{|X_u(t)| \leq \log r\}}\right)\right)\\
	& \geq 1- \mathbb{E}_{(y, i)}\left(\prod_{u\in Z(t)} \left(1- \left(1-\frac{\varepsilon}{2}\right)1_{\{|X_u(t)| \leq \log r\}}\right)\right)\\
	& \geq 1-\frac{\varepsilon}{2} - \mathbb{P}_{(y,i)}\left(\# \left\{u\in Z(t): |X_u(t)| \leq \log r\right\}=0 \right)> 1-\varepsilon.
\end{align*}
This completes the proof.

\hfill$\Box$

 For single-type BBM, the assumption $p_0=0, p_1\neq 1$ implies that the offspring mean is strictly larger than $1$.
This fact is used  in  the inequality
$k(v(t-s, \mathfrak{Z}_{x,y}(s)))\leq 1/2$
 above (7.34) on   page 107  of \cite{Bramson1}
to prove the exponential decay in \cite[Lemma 7.2]{Bramson1},
where $\mathfrak{Z}_{x,y}$ is the Brownian bridge starting at $x$ and ending at $y$.
 But for multitype BBM, the assumptions $p_{\mathbf{0}}(i)=0$ for all $i\in S$ and $\lambda^*>0$ do not imply $n_i=\sum_{j=1}^d m_{i,j}>1$ for all $i$. We can only get that there exists $j_0\in S$ such that $n_{j_0}>1$.
For multitype BBM, the fact $n_{j_0}>1$ will play a role in getting the exponential decay in Lemma \ref{lemma3} below.
In Lemma \ref{lemma8} below, we give an estimate which will replace the role of  the inequality
$k(v(t-s, \mathfrak{Z}_{x,y}(s)))\leq 1/2$   of \cite{Bramson1} in the multi-type case.

\begin{lemma}\label{lemma8}
	Suppose that $n_{j_0}>1$. Then for any $\theta >0$, there exist $C_\theta >0$ and $\varepsilon= \varepsilon_\theta >0$ such that for all $t>0$,
	\begin{align}
		\sup_{x\in\mathbb{R}, i\in S} \mathbf{E}_{(x,i)}^h\left( \exp\left\{ -\theta \sum_{s\in D_J: s\leq t} 1_{\{I_s= j_0\}}\right\}\right) \leq C_\theta e^{-\varepsilon t}.
	\end{align}
\end{lemma}
\textbf{Proof: }  Since $I$ and $X$ are independent,  we only consider the case  $x=0.$ Let $\{Y_n: n=0,1,...\}$ be the embedded chain of $I_t$ under $\mathbf{P}_{(0,i)}$.
We first prove that there exist $C_1, \delta_1>0$ such that
\begin{align}\label{Ineq:Emb-chain}
	\sup_{i\in S} \mathbf{P}_{(0,i)}^h \left(Y_0 \neq j_0,\cdots, Y_n\neq j_0\right) \leq C_1 e^{-\delta_1 n}.
\end{align}
Since $\{Y_n\}$ is irreducible, for each $i\in S$, there exists $L_i\in \mathbb{N}$ such that $\mathbf{P}_{(0,i)}^h (Y_{L_i}\neq j_0)<1$.
Let $L:= \max_{i\in S} L_i$. Then
\[
\sup_{i\in S} \mathbf{P}_{(0,i)}^h \left(Y_0 \neq j_0,\cdots, Y_L\neq j_0\right) \leq  \sup_{i\in S}\mathbf{P}_{(0,i)}^h (Y_{L_i}\neq j_0)=: e^{-\varepsilon_1}.
\]
Therefore, for $n\geq L$, we have
\[
	\sup_{i\in S} \mathbf{P}_{(0,i)}^h \left(Y_0 \neq j_0,\cdots, Y_n\neq j_0\right) \leq e^{-\varepsilon_1}\sup_{i\in S} \mathbf{P}_{(0,i)}^h \left(Y_0 \neq j_0,..., Y_{n-L}\neq j_0\right),
\]
which implies \eqref{Ineq:Emb-chain} with $C_1:= e^{\varepsilon_1}, \delta_1:= \varepsilon_1/L.$

Next, define $U_{j_0}:= \inf\{t\in D_J, I_t = j_0 \}$. We prove that there exists $\delta_2>0$ such that
\begin{align}\label{Exp-moment-hitting-time}
	\sup_{i\in S} \mathbf{E}_{(0,i)}^h\left(e^{\delta_2 U_{j_0}} \right)=: C_2<\infty.
\end{align}
To this end,
it suffices to show that there exist constants $C_3, \varepsilon_2>0$ such that for $t$ large enough,
\begin{align}\label{Tail-Pro-Hitting-time}
        \sup_{i\in S} \mathbf{P}_{(0,i)}^h \left(U_{j_0}> t \right)\leq C_3 e^{-\varepsilon_2 t}.
\end{align}
Recall that in the paragraph containing \eqref{Coupling4},
we defined a coupling $(I_t, I_t^{\underline{a}})$ so that the embedded chain of $I$ and $I^{\underline{a}}$ are the same, and the jump times $D_J=\{t_n: 0<t_1<t_2<\dots\}$ of $I$ and the jumps times $D^{\underline{a}}=\{0<t_1^{\underline{a}}<t_2^{\underline{a}}<\dots\} $ of $I^{\underline{a}}$
satisfy $t_n\le t_n^{\underline{a}}$ for all $n\ge 1$.
Let $U_{j_0}^{\underline{a}}:= \inf\{t\in D^{\underline{a}}, I_t^{\underline{a}} = j_0 \}$, then $U_{j_0} \leq U_{j_0}^{\underline{a}}$.
For $n\in \mathbb{N}$, on the event that the first hitting time of $j_0$ by the embedded chain is larger than $n$, by \eqref{Ineq:Emb-chain} we can bound $\mathbf{P}_{(0,i)}^h \left(U_{j_0}> t \right)$ from above by $C_1e^{-\delta_1 n}$.
On the event that the first hitting time of $j_0$ by the embedded chain is less than or equal to $n$, we bound $\mathbf{P}_{(0,i)}^h \left(U_{j_0}> t \right)$ from above by
\[
\sup_{i\in S} \sum_{m=1}^n \mathbf{P}_{(0,i)}^h\left(Y_1^{\underline{a}}+...+Y_m^{\underline{a}} >t\right) \leq n e^{-\underline{a}t/2} \sup_{i\in S} \left(\mathbf{E}_{(0,i)}^h\left( e^{\underline{a}Y_1^{\underline{a}}/2}\right)\right)^n = ne^{-\underline{a}t/2} 2^n.
\]
Taking $n= [\varepsilon_2 t]$ for $0<\varepsilon_2 \log 2 < \underline{a}/2$, we get
\begin{align}
	\sup_{i\in S} \mathbf{P}_{(0,i)}^h \left(U_{j_0}> t \right) \leq C_1 e^{-\delta_1[\varepsilon_2 t]} + [\varepsilon_2 t]e^{-\underline{a}t/2} e^{\log 2 [\varepsilon_2 t]},
\end{align}
which implies \eqref{Exp-moment-hitting-time}.

Define $V_{j_0}^1:= \inf\left\{t\in D_J, I_t =j_0\right\}$ and $V_{j_0}^n:= \inf\left\{t\in D_J: t> V_{j_0}^{n-1}: I_t= j_0 \right\}$ for $n\geq 2$. Set $U_{j_0}^1= V_{j_0}^1$ and $U_{j_0}^n:= V_{j_0}^n - V_{j_0}^{n-1}$. By the strong Markov property,  $\{U_{j_0}^n:n\ge 1\}$ are independent.
Define $S_t:= \sum_{s\in D_J: s\leq t} 1_{\{I_s= j_0\}}
= \sup\left\{ n:\ \sum_{m=1}^n U_{j_0}^m \leq t\right\}$, then for any $n$, $\left\{S_t = n\right\}\subset \left\{ \sum_{m=1}^{n+1} U_{j_0}^m >t\right\}$. Thus, by \eqref{Exp-moment-hitting-time} and the strong Markov property,
\begin{align}
	&\sup_{x\in\mathbb{R}, i\in S} \mathbf{E}_{(x,i)}^h\left(e^{-\theta S_t}\right)\leq e^{-\theta n}+ \sup_{x\in\mathbb{R}, i\in S} \sum_{\ell =1}^n \mathbf{P}_{(x,i)}^h\left(S_t =\ell \right)\nonumber\\
	& \leq e^{-\theta n}+ \sup_{x\in\mathbb{R}, i\in S} \sum_{\ell=1}^n \mathbf{P}_{(x,i)}^h\left( \sum_{m=1}^{\ell+1} U_{j_0}^m >t \right) \leq e^{-\theta n}+ n \sup_{x\in\mathbb{R}, i\in S} \mathbf{P}_{(x,i)}^h\left( \sum_{m=1}^{n+1} U_{j_0}^m >t \right) \nonumber\\
	& \leq e^{-\theta n}+ n e^{-\delta_2 t } \sup_{i\in S} \mathbf{E}_{(0,i)}^h \left(e^{\delta_2 \sum_{m=1}^{n+1} U_{j_0}^m}\right) \leq e^{-\theta n}+ n e^{-\delta_2 t } \prod_{m=1}^{n+1} 	\sup_{i\in S} \mathbf{E}_{(0,i)}^h\left(e^{\delta_2 U_{j_0}} \right)\nonumber\\
	& = e^{-\theta n} + n e^{(n+1)\delta_2 -\delta_2 t}.
\end{align}
Taking $n= [t/2]$, we get the conclusion of the lemma.

\hfill$\Box$

Our goal is to  get the upper bound for $v_i(t, m(t)+x)$ for large $r$ and $r<x\leq \sqrt{t}$ in Proposition \ref{prop3}.
Lemma \ref{lemma2} implies that the upper bound is related to $\mathbf{E}_{(x,i)}^h\left(R(t-r;v^*) \frac{v_{I_{t-r}}(r,y)}{h_{I_{t-r}}}\bigg| X_{t-r}=y \right)$ for $m(t)+ r< x\leq m(t)+ \sqrt{t}$ and $y>-\log r$.
In Lemma \ref{lemma3} and Proposition \ref{prop1} below,
we will estimate $\mathbf{E}_{(x,i)}^h\left(R(t-r;v^*) \frac{v_{I_{t-r}}(r,y)}{h_{I_{t-r}}}; (\mathcal{B}_{up})^c\bigg| X_{t-r}=y \right)$ under the condition $x \geq m(t)$ and $y> -\log r$.
We first prove Lemma \ref{lemma3},
which is an analog of \cite[Lemma 7.2, p.105]{Bramson1}.
 Recall that $j_1$ is the integer such that $j_1<t/2 \leq j_1+1$.

\begin{lemma}\label{lemma3}
	Let $\mathbf{v}$ be a solution to \eqref{F-KPP2} with initial value satisfying
	\eqref{Initial-Assumption} and let $\mathbf{v}_*$ be given by \eqref{V-Star}.
	Then for $r$ large enough, $t>4r$ and	$j_1\ge j\geq [r+r^{\delta}]$, it holds that
	\[
	\mathbf{E}_{(x,i)}^h\left(R((r, t-r];v_*) \frac{v_{I_{t-r}}(r,y)}{h_{I_{t-r}}}; A_j^1(r,t) \big| X_t =y \right)
	\lesssim e^{-j^\delta/C}\mathbf{E}_{(x,i)}^h \left(\frac{v_{I_{t-r}}(r,y)}{h_{I_{t-r}}} \right){\mathbf{E}_{(x,i)}^h \left(A_j^1(r,t)\big| X_t =y \right)}
	\]
	for all $y > -\log r, x \geq m(t)$ and some constant $C$.
	\end{lemma}
\textbf{Proof: }
First note that $j_1\ge j\geq [r+r^{\delta}]$ implies that $r\le j-j^\delta/2<j\le t-r$ and $r\le t-j<t-j+j^\delta/2 \leq t-r$.
When $r$ is large enough, we have
\begin{align}
	& \mathbf{E}_{(x,i)}^h\left(
	R((r, t-r];v_*) \frac{v_{I_{t-r}}(r,y)}{h_{I_{t-r}}}; A_j^1(r,t) \big| X_t =y \right)\\
	& \leq \mathbf{E}_{(x,i)}^h\left(R_{t-r}((j-j^\delta/2, j];v_*) \frac{v_{I_{t-r}}(r,y)}{h_{I_{t-r}}}; A_j^1(r,t), S(r,t)=S^1(r,t) \big| X_t =y \right) \\
	&\quad + \mathbf{E}_{(x,i)}^h\left(R_{t-r}((t-j, t-j+j^\delta/2];v_*) \frac{v_{I_{t-r}}(r,y)}{h_{I_{t-r}}}; A_j^1(r,t), S(r,t)=S^2(r,t) \big| X_t =y \right).
\end{align}
For the first term, when $s\in [2r, t/2]$,  by \cite[(7.30), p. 106]{Bramson1},
\begin{align}\label{e:newRS}
	\underline{\mathcal{M}}_{r,t}'(t-s)&= m(s_1)- \left(K +\frac{\alpha(r) - m(t)}{t}\right)s^\delta +o_1(1),
\end{align}
where $s_1= t-s-s^{\delta}+o_2(1)$ and $o_1(1), o_2(1)\to 0$ as $r\to\infty$. Since $K> 2+ \sqrt{2\lambda^*}$, for $r$ large enough, we have for any $s\in [2r, t/2], t>4r$ and $s'\in [2r, t/2]\cap [s, s+1)$,
\begin{align}
	\underline{\mathcal{M}}_{r,t}'(t-s') \leq m(s_1)-2s^{\delta}.
\end{align}
Now we prove that, for any $\varepsilon>0$,  there exists $N(\varepsilon)$ such that for all $r>N(\varepsilon)$,
$i\in S, t>4r$,
$s\in [2r, t/2]$ and
$0\leq  y' \leq \underline{\mathcal{M}}_{r,t}'(t-s') + 2s^\delta\leq m(s_1)$
with some $ s'\in [2r, t/2]\cap [s, s+1)$,
\begin{align}\label{step_26}
	v_i^* (t-r-s, y')\geq 1- \varepsilon.
\end{align}
When $y' \in [0, m_+(t-r-s)]$, by Lemma \ref{lemma5}, we can find  $N_1(\varepsilon)$ such that for $r> N_1(\varepsilon)$,
\begin{align}\label{step_65}
	v_i^* (t-r-s, y') \geq  1-\varepsilon.
\end{align}
Combining the above with \eqref{e:newRS} we get \eqref{step_26}
when  $m_+(t-r-s)> m(s_1)$. If $m_+(t-r-s)\leq m(s_1)$, then  for $y'\in [m_+(t-r-s), m(s_1)]$, by \eqref{Lower-Bound-M-T2} and \eqref{step_64},
\begin{align}\label{step_66}
	 &v_i^* (t-r-s, y')\geq v_i (t-s, y')- \P_{(0,i)} \left(M_{t-r-s}> y'+\log r\right)\nonumber\\
	 & \geq \P_{(0,i )}\left(M_{t-s}^{i_0} > y'-N_1\right)
	 - \P_{(0,i)} \left(M_{t-r-s}> y'+\log r\right)\nonumber\\
	  & \geq \P_{(0,i )}\left(M_{t-s}^{i_0} > m(s_1)-N_1\right)
	 - \P_{(0,i)} \left(M_{t-r-s}> m_+(t-r-s)+ \log r\right),
\end{align}
where $i_0\in S$ is the type fixed in \eqref{Initial-Assumption}.
Note that $ t-s -s_1= s^{\delta}- o_2(1)$ and $m'(s)\geq \sqrt{\lambda^*}$ for large $s$, when $r$ is large enough,
\begin{align}
	m(t-s)- m(s_1)+ N_1 \geq  \sqrt{\lambda^*}(t-s- s_1)+N_1 \geq \frac{\sqrt{\lambda^*}}{2} s^\delta +N_1 \to +\infty.
\end{align}
Therefore,  by Theorem \ref{Tightness-2}, there exists $N_2(\varepsilon)$ such that for $r\geq N_2(\varepsilon)$,
\begin{align*}
	 & \P_{(0,i )}\left(M_{t-s}^{i_0} - m(t-s)> m(s_1)-m(t-s)-N_1\right)
	\geq 1-\frac{\varepsilon}{2},\\
	& \P_{(0,i)} \left(M_{t-r-s}- m_+(t-r-s)> \log r\right)
	\leq \frac{\varepsilon}{2}.
\end{align*}
Putting these inequalities back to \eqref{step_66}, we have that $v_i^* (t-r-s, y') \geq 1-\varepsilon$ when $r>N_2(\varepsilon)$ and $y'\in [m_+(t-r-s), m(s_1)]$, which, together with \eqref{step_65}, implies \eqref{step_26}.

On the event $A_j^1(r,t)\cap \{S(r,t)= S^1(r,t)\in [j, j+1) \}$, set
\[
E_j := \left\{X_s - X_{S^1(r,t)}\leq 2j^\delta,\quad
\forall \
s\in [j-j^\delta/2, S^1(r,t)]
\right\}.
\]
Then on $A_j^1(r,t)\cap \{S(r,t)= S^1(r,t)\in [j, j+1) \}\cap E_j$, it holds that
\begin{align}
	X_s \leq  2j^{\delta} + X_{S^1(r,t)}= 2j^{\delta} + \underline{\mathcal{M}}_{r,t}'\left(t-S^1(r,t)\right) .
\end{align}
By \eqref{step_26}, uniformly for $i \in S$,  on $A_j^1(r,t)\cap \{S(r,t)=S^1(r,t) \} \cap E_j$,
\[
v_i^*(t-s, X_s)\geq 1-\varepsilon,\quad \mbox{for }s\in [j-j^{\delta}/2, j].
\]
This implies that on $A_j^1(r,t)\cap \{S(r,t)=S^1(r,t) \} \cap E_j$,
for $s\in [j-j^\delta /2, j]$,
\[
\frac{ \varphi_{I_s}\left(\mathbf{v}^*(t-s, X_s)\right)}{\sum_{j=1}^d m_{I_s, j}v_j^*(t-s, X_s)}\leq \frac{1}{(1-\varepsilon)n_{I_j}}1_{\{I_j= j_0\}} + 1_{\{I_{j}\neq j_0\}} .
\]
We have $n_{j_0}>1$ by assumption. Choose an $\varepsilon >0$ sufficient small and an appropriate $\eta <1$ so that
\[
\frac{1}{(1-\varepsilon)n_{j_0}} \leq \eta <1.
\]
By  \cite[(7.36)]{Bramson1}, under the assumption $y>-\log r$ and $x \geq m(t)$, for $r$ large enough, we have
\begin{align}
	\mathbf{P}_x \left(E_j^c \big| A_j^1(r,t), S(r,t)=S^1(r,t), X_t =y\right)\leq e^{-j^\delta /4}
\end{align}
Using the independence of $X$ and $I$, we get
\begin{align}
	&\mathbf{E}_{(x,i)}^h\left(R_{t-r}((j-j^\delta/2, j];v_*)  \frac{v_{I_{t-r}}(r,y)}{h_{I_{t-r}}}; A_j^1(r,t), S(r,t)= S^1(r,t)\big| X_t =y \right)\\
	& \leq \mathbf{E}_{(x,i)}^h\left( \frac{v_{I_{t-r}}(r,y)}{h_{I_{t-r}}}; E_j^c, A_j^1(r,t), S(r,t)= S^1(r,t)\big| X_t =y \right)\\
	& \quad + \mathbf{E}_{(x,i)}^h\left(\exp\left\{\sum_{s\in D_J, j -j^{\delta}/2 < s\leq j} 1_{\{I_s = j_0\}}\log\eta \right\} \frac{v_{I_{t-r}}(r,y)}{h_{I_{t-r}}};  E_j , A_j^1(r,t), S(r,t)= S^1(r,t)\big| X_t =y \right)\\
	& \leq  e^{-j^{\delta}/4}\mathbf{E}_{(x,i)}^h\left( \frac{v_{I_{t-r}}(r,y)}{h_{I_{t-r}}}\right) \mathbf{P}_{(x,i)}^h\left( A_j^1(r,t), S(r,t)= S^1(r,t)\big| X_t =y \right) \\
	& \quad +  \mathbf{E}_{(x,i)}^h\left(\exp\left\{\sum_{s\in D_J, j -j^{\delta}/2 < s\leq j} 1_{\{I_s=j_0\}}\log\eta \right\} \frac{v_{I_{t-r}}(r,y)}{h_{I_{t-r}}}\right) \mathbf{P}_{(x,i)}^h\left( A_j^1(r,t), S(r,t)= S^1(r,t)\big| X_t =y \right).
\end{align}
By Lemma \ref{lemma8}, and the fact that $\inf_{j\in S}\inf_{r\geq 1, t>4r} \mathbf{P}_{(x,i)}^h \left(I_{t-r}= j\right) >c >0,$ we have
\begin{align}
	&\mathbf{E}_{(x,i)}^h\left(\exp\left\{\sum_{s\in D_J, j -j^{\delta}/2 < s\leq j} 1_{\{I_s = j_0\}} \log\eta \right\} \frac{v_{I_{t-r}}(r,y)}{h_{I_{t-r}}}\right) \\
	& \leq \sum_{j=1}^d \frac{v_j(r,y)}{h_j} \sup_{x\in \mathbb{R}, \ell \in S}\mathbf{E}_{(x,\ell)}^h\left(\exp\left\{\sum_{s\in D_J, s\leq j^\delta/2} 1_{\{I_s = j_0\}} \log\eta \right\} \right)
	\leq C_\eta \sum_{j=1}^d \frac{v_j(r,y)}{h_j} e^{-\varepsilon j^\delta/2} \\ &
	\leq \frac{C_\eta}{c} \mathbf{E}_{(x,i)}^h \left(\frac{v_{I_{t-r}}(r,y)}{h_{I_{t-r}}} \right) e^{-\varepsilon j^\delta/2}.
\end{align}
Combining the two displays above, we get
\begin{align}
	&\mathbf{E}_{(x,i)}^h\left(R_{t-r}((j-j^\delta/2, j];v_*)  \frac{v_{I_{t-r}}(r,y)}{h_{I_{t-r}}}; A_j^1(r,t), S(r,t)= S^1(r,t)\big| X_t =y \right) \\
	&\leq \mathbf{E}_{(x,i)}^h\left( \frac{v_{I_{t-r}}(r,y)}{h_{I_{t-r}}}\right) \mathbf{P}_{(x,i)}^h\left(  A_j^1(r,t), S(r,t)= S^1(r,t)\big| X_t =y \right) \left( e^{-j^{\delta}/4} + \frac{C_\eta}{c}e^{-\varepsilon j^\delta/2}\right)\\
	& \lesssim  e^{-j^\delta /C}\mathbf{E}_{(x,i)}^h\left( \frac{v_{I_{t-r}}(r,y)}{h_{I_{t-r}}}\right) \mathbf{P}_{(x,i)}^h\left(  A_j^1(r,t), S(r,t)= S^1(r,t)\big| X_t =y 	\right)
\end{align}
with
\[
C= \max\left\{ 4, \frac{2}{\varepsilon}\right\}.
\]
The second term can be treated similarly. Thus the assertion of the lemma is valid.

\hfill$\Box$

The next proposition is similar to \cite[Proposition 7.3, p. 108]{Bramson1}.
\begin{prop}\label{prop1}
	Let $\mathbf{v}$ be the solution of \eqref{F-KPP2} with initial value satisfying \eqref{Initial-Assumption}. Then for $r$ large enough, $t> 4r$, $y > -\log r$ and $x \geq m(t)$, it holds that
	\begin{align}
		\mathbf{E}_{(x,i)}^h\left(
		R((r, t-r];v_*)
		\frac{v_{I_{t-r}}(r,y)}{h_{I_{t-r}}}; (\mathcal{B}_{up})^c \big| X_t =y \right) \leq \frac{1}{r^2}\mathbf{E}_{(x,i)}^h \left(\frac{v_{I_{t-r}}(r,y)}{h_{I_{t-r}}} \right)  \mathbf{P}_{(x,i)}^h \left(\mathcal{B}_{up} \big| X_t =y\right).
	\end{align}
\end{prop}
\textbf{Proof: } Note that
\[
(\mathcal{B}_{up})^c \subset
\bigcup_{j=[r+r^\delta]}^{j_1} A_j(r,t) = \bigcup_{j=[r+r^\delta]}^{j_1} A_j^1(r,t) \bigcup \bigcup_{j=[r+r^\delta]}^{j_1} A_j^2(r,t).
\]
By Lemma \ref{lemma3}, and the independence of $X$ and $I$, for large $r$, when $y
>-\log r $ and $x \geq m(t)$,
\begin{align}
	&\mathbf{E}_{(x,i)}^h\left(R((r, t-r];v_*) \frac{v_{I_{t-r}}(r,y)}{h_{I_{t-r}}}; (\mathcal{B}_{up})^c \big| X_t =y \right) \\
	&\leq \sum_{j=[r+r^\delta]}^{j_1}\mathbf{E}_{(x,i)}^h\left(R((r, t-r];v_*) \frac{v_{I_{t-r}}(r,y)}{h_{I_{t-r}}}; A_j^1(r,t) \big| X_t =y \right) \\
	& \quad+\sum_{j=[r+r^\delta]}^{j_1}\mathbf{E}_{(x,i)}^h\left(\frac{v_{I_{t-r}}(r,y)}{h_{I_{t-r}}}; A_j^2(r,t) \big| X_t =y \right)\\
	& \lesssim \sum_{j=[r+r^\delta]}^{j_1} e^{-j^\delta/C}\mathbf{E}_{(x,i)}^h \left(\frac{v_{I_{t-r}}(r,y)}{h_{I_{t-r}}} \right)\mathbf{P}_{(x,i)}^h \left(A_j^1(r,t)\big| X_t =y \right)  \\ & \quad +\sum_{j=[r+r^\delta]}^{j_1}\mathbf{E}_{(x,i)}^h\left( \frac{v_{I_{t-r}}(r,y)}{h_{I_{t-r}}} \right)\mathbf{P}_{(x,i)}^h\left( A_j^2(r,t) \big| X_t =y \right).
\end{align}
Note that the estimates of the probabilities of $A_j^1(r,t)$ and $A_j^2(r,t)$ only relies on the path of Brownian bridge,
using Lemma \ref{lemma4} and the argument on \cite[p. 109]{Bramson1}, we get
\begin{align}
	&\sum_{j=[r+r^\delta]}^{j_1} e^{-j^\delta/C}\mathbf{P}_{(x,i)}^h \left(A_j^1(r,t)\big| X_t =y \right) \lesssim \sum_{j=[r+r^\delta]}^{j_1} e^{-j^\delta/C}\mathbf{P}_{(x,i)}^h \left(\mathcal{B}_{up}^{j+1}\big| X_t =y \right) \\
	& \lesssim \sum_{j=[r+r^\delta]}^{j_1} e^{-j^\delta/C}\frac{j+1}{r}\mathbf{P}_{(x,i)}^h \left(\mathcal{B}_{up}\big| X_t =y \right) \leq \frac{1}{r}\mathbf{P}_{(x,i)}^h \left(\mathcal{B}_{up}\big| X_t =y \right) \times  \sum_{j=[r+r^\delta]}^{\infty} (j+1) e^{-j^\delta/C},
\end{align}
and
\begin{align}
	&\sum_{j=[r+r^\delta]}^{j_1} \mathbf{P}_{(x,i)}^h\left( A_j^2(r,t) \big| X_t =y \right)\leq \sum_{j=[r+r^\delta]}^{j_1}e^{-j/2} \mathbf{P}_{(x,i)}^h\left( \mathcal{B}_{up}^{j+1} \big| X_t =y \right) \\
	&\lesssim\frac{1}{r} \mathbf{P}_{(x,i)}^h\left( \mathcal{B}_{up} \big| X_t =y \right) \times  \sum_{j=[r+r^\delta]}^{\infty}(j+1)e^{-j/2} .
\end{align}
Therefore, we conclude that
\begin{align}\label{step_14}
	&\mathbf{E}_{(x,i)}^h\left(R((r, t-r];v_*) \frac{v_{I_{t-r}}(r,y)}{h_{I_{t-r}}}; (\mathcal{B}_{up})^c \big| X_t =y \right)\\
	& \lesssim \mathbf{E}_{(x,i)}^h \left(\frac{v_{I_{t-r}}(r,y)}{h_{I_{t-r}}} \right) \times \frac{1}{r}\mathbf{P}_{(x,i)}^h \left(\mathcal{B}_{up}\big| X_t =y \right) \times \left(\sum_{j=[r+r^\delta]}^{\infty} (j+1) e^{-j^\delta/C} +  \sum_{j=[r+r^\delta]}^{\infty}(j+1)e^{-j/2}  \right).
\end{align}
As $r\to\infty$, the last term of \eqref{step_14} decays faster than $r^{-1}$,
thus the assertion of the proposition is valid.
\hfill$\Box$

We  now give an upper bound for $v_i(t,m(t)+x)$:

\begin{prop}\label{prop3}
	 Suppose that $\mathbf{v}$ satisfies \eqref{F-KPP2} and \eqref{Initial-Assumption}. Let $r$ be large enough, then for all $r< x \leq  \sqrt{t}$, it holds that
	\begin{align}
		v_i(t,m(t)+ x)&\leq C_{up} (r) e^{\lambda^* (t-r)}h_i\int_{-\log r}^\infty  \frac{e^{-\frac{(m(t)+x-y)^2}{2(t-r)}}}{\sqrt{2\pi (t-r)}} \times \mathbf{P}_{(m(t)+x,i)}^h\left(\mathcal{B}_{up} \Big|X_{t-r}=y  \right)\sum_{j=1}^d g_j v_j(r,y)	 \mathrm{d} y,
	\end{align}
         where $C_{up}(r)\downarrow 1$ as $r\to\infty$.
\end{prop}
\textbf{Proof: } By Lemma \ref{lemma2}, Proposition \ref{prop1},  and the independence of $X$ and $I$, we have
\begin{align}\label{step_13}
	&v_i (t,m(t)+x) \\
	& \leq C(r) e^{\lambda^* (t-r)}h_i\int_{-\log r} ^\infty \frac{e^{-\frac{(m(t)+x-y)^2}{2(t-r)}}}{\sqrt{2\pi (t-r)}}  \mathbf{E}_{(m(t)+x,i)}^h\left(R(t-r;v_*) \frac{v_{I_{t-r}}(r,y)}{h_{I_{t-r}}}\bigg| X_{t-r}=y \right)\mathrm{d} y \\
	& \leq   C(r) e^{\lambda^* (t-r)}h_i\int_{-\log r} ^\infty \frac{e^{-\frac{(m(t)+x-y)^2}{2(t-r)}}}{\sqrt{2\pi (t-r)}}  \mathbf{E}_{(m(t)+x,i)}^h\left( \frac{v_{I_{t-r}}(r,y)}{h_{I_{t-r}}};\mathcal{B}_{up}\bigg| X_{t-r}=y \right)\mathrm{d} y \\
	& \quad + C(r) e^{\lambda^* (t-r)}h_i\int_{-\log r} ^\infty \frac{e^{-\frac{(m(t)+x-y)^2}{2(t-r)}}}{\sqrt{2\pi (t-r)}}  \mathbf{E}_{(m(t)+x,i)}^h\left(R(t-r;v_*) \frac{v_{I_{t-r}}(r,y)}{h_{I_{t-r}}}; (\mathcal{B}_{up})^c \bigg| X_{t-r}=y \right)\mathrm{d} y \\
	& \leq  C(r) e^{\lambda^* (t-r)}h_i\int_{-\log r}^\infty \frac{e^{-\frac{(m(t)+x-y)^2}{2(t-r)}}}{\sqrt{2\pi (t-r)}} \times  \mathbf{P}_{(m(t)+x,i)}^h\left( \mathcal{B}_{up} \bigg| X_{t-r}=y \right) \mathbf{E}_{(m(t)+x,i)}^h\left( \frac{v_{I_{t-r}}(r,y)}{h_{I_{t-r}}} \right)\mathrm{d} y \\
	& \quad + C(r) e^{\lambda^* (t-r)}h_i\int_{-\log r}^\infty \frac{e^{-\frac{(m(t)+x-y)^2}{2(t-r)}}}{\sqrt{2\pi (t-r)}}  \left(\frac{1}{r^2}\mathbf{E}_{(m(t)+x,i)}^h \left(\frac{v_{I_{t-r}}(r,y)}{h_{I_{t-r}}} \right)  \mathbf{P}_{(m(t)+x,i)}^h \left(\mathcal{B}_{up} \big| X_{t-r} =y\right)\right)\mathrm{d} y\\
	& =  C(r)\left(1+\frac{1}{r^2}\right) e^{\lambda^* (t-r)}h_i\int_{-\log r}^\infty \frac{e^{-\frac{(m(t)+x-y)^2}{2(t-r)}}}{\sqrt{2\pi (t-r)}}   \mathbf{P}_{(m(t)+x,i)}^h\left( \mathcal{B}_{up} \bigg| X_{t-r}=y \right) \mathbf{E}_{(m(t)+x,i)}^h\left( \frac{v_{I_{t-r}}(r,y)}{h_{I_{t-r}}} \right)\mathrm{d} y.
\end{align}
Note that  $\lim_{r\to\infty} \mathbf{P}_{(x,i)}^h\left(I_r = j\right)= \mu_j=g_j h_j = \lim_{r\to\infty} \sup_{t>r} \mathbf{P}_{(x,i)}^h\left(I_t = j\right)$. Letting
\[
C_{up}(r):= C(r)\left(1+\frac{1}{r^2}\right)\sup_{i,j\in S}\sup_{t>r} \frac{\mathbf{P}_{(x,i)}^h\left(I_t = j\right)}{g_j h_j},
\]
we get the assertion of the proposition.
\hfill$\Box$

\subsection{Lower bound for $\mathbf{v}$}

Similar to \cite[(7.42), p. 111 and (7.9), p. 99]{Bramson1} or \cite[(2.10)]{RSZ},
we define
\begin{align}
	\overline{\mathcal{M}}_{r,t}^{x}(s):= 	\left\{\begin{array}{ll}\overline{L}_{r,t}(s) + \frac{s}{t}m(t) -\frac{t-s}{t}\log r,\quad &  0\leq s \leq t-2r,\\
		\displaystyle \frac{x}{2}+m(t),\quad &t-2r < s\leq t,
	\end{array}\right.
\end{align}
where $\overline{L}_{r,t}$ is defined in the beginning of Subsection \ref{ss:5.1}. Define
\[
    \mathcal{B}_{low}:=\left\{X_s> \overline{\mathcal{M}}_{r,t}^{x}(t-s),\ s\in [0,t-r] \right\}.
\]

\begin{prop}\label{prop2}
When $r$ is large enough, it holds that for all $t,x> 8r$ and all $i\in S$,
	\begin{align}
		v_i(t,m(t)+x)&\geq C_{low} (r) e^{\lambda^* (t-r)}h_i\int_\R  \frac{e^{-\frac{(m(t)+x-y)^2}{2(t-r)}}}{\sqrt{2\pi (t-r)}}  \mathbf{P}_{(m(t)+x,i)}^h\left(\mathcal{B}_{low} \Big|X_{t-r}=y  \right)\sum_{j=1}^d g_j v_j(r,y )
	 \mathrm{d} y,
	\end{align}
where $C_{low}(r)\uparrow 1$ as $r\to\infty$.
\end{prop}
\textbf{Proof: } When $s\in [0, 2r]$, on the set $\mathcal{B}_{low}$, for $t,x> 8r$  and $r$ large enough, we have
\[
X_s > \overline{\mathcal{M}}_{r,t}^{x}(t-s)\geq m(t)+4r \geq  m(t-s) + 4r + \sqrt{2\lambda^*}s +O(1)\geq m(t-s) + r+N_2,
\]
where $N_2$ is the constant in  \eqref{Initial-Assumption}. Note that for any $y > m(t-s)+ r+N_2$, by \eqref{Initial-Assumption}, \eqref{Prob-Representation} and  Proposition \ref{Upper-Bound-M-T},
\begin{align}
	&v_i(t-s, y)  = \E_{(y,i)} \left(1-\prod_{u \in Z(t-s)}\left(1-v_{I_u(t-s)}(0, X_u(t-s))\right)\right)\\
	& \leq  \E_{(y,i)}\left(1-\prod_{u \in Z(t-s)}\left(1-1_{(-\infty, N_2)}(X_u(t-s))\right)\right)\\
	&   = \P_{(0,i)} (M_{t-s}> y-N_2) \leq \P_{(0,i)} (M_{t-s} \geq m(t-s) +r)\lesssim  re^{-\sqrt{2\lambda^*}r}.
\end{align}
Using \eqref{Assumption1} and noting that $\log x \sim x-1$ as $x\to 1$, we get that
when $r$ is large enough, on $\mathcal{B}_{low}$,
\begin{align}\label{step_49}
	\sum_{s\in D_J, s\leq 2r} \log \left(\frac{ \varphi_{I_s}\left(\mathbf{v}(t-s, X_s)\right)}{\sum_{j=1}^d m_{I_s, j}v_j (t-s, X_s)}\right) \gtrsim -\sum_{s\in D_J, s\leq 2r} r^{\alpha_0} e^{-\alpha_0 \sqrt{2\lambda^*}r}
	=: \Gamma_1(0,2r).
\end{align}
Now we deal with the case $s\in [2r, t-r]$. Similar as above, when $r$ is large enough,
for all $s\in [2r, t-r]$ and $y> m (s+s^\delta \land (t-s)^\delta)$,
\begin{align}
	v_i(s,y)&\leq \P_{(0,i)}\left(M_s > y -N_2 \right)\leq \P_{(0,i)}\left(M_s-m(s) > m (s+s^\delta \land (t-s)^\delta) - m(s) -N_2 \right)\\
	& = \P_{(0,i)}\left(M_s-m(s) > \sqrt{2\lambda^*}\left(s^\delta \land (t-s)^\delta \right) + O(1) \right)\lesssim \left(s^\delta \land (t-s)^\delta\right) e^{-2\lambda^* \left(s^\delta \land (t-s)^\delta\right)}\\
	& \lesssim  e^{-\lambda^* \left(s^\delta \land (t-s)^\delta\right)}.
\end{align}
In this case, when $r$ is large enough
(see the display below \cite[(2.14)]{RSZ}),
\[
\overline{\mathcal{M}}_{r,t}^{x}(s)\geq m(s+s^\delta \land(t-s)^\delta) + s^\delta \land(t-s)^\delta\left(K- \frac{m(t)}{t}- \frac{\log r}{t}\right)> m (s+s^\delta \land(t-s)^\delta)
\]
since $K> \sqrt{2\lambda^*}$. Therefore, on $\mathcal{B}_{low}$, when $r$ is large enough,
\begin{align}\label{step_50}
	\sum_{s\in D_J, 2r < s\leq t-r} \log \left(\frac{ \varphi_{I_s}\left(\mathbf{v}(t-s, X_s)\right)}{\sum_{j=1}^d m_{I_s, j}v_j (t-s, X_s)}\right) & \gtrsim -  \sum_{s\in D_J, 2r < s\leq t-r} e^{- \lambda^* \left(s^\delta \land (t-s)^\delta\right)}
	\\ & = : \Gamma_2(2r, t-r).
\end{align}
Note that $\mathcal{B}_{low}$ is independent of $D_J$. By \eqref{step_49} and \eqref{step_50}, there exists a constant $c$ such that
\begin{align}\label{step_11}
	v_i(t,m(t)+x) &\geq e^{\lambda^* (t-r)}h_i\int_\R  \frac{e^{-\frac{(m(t)+x-y)^2}{2(t-r)}}}{\sqrt{2\pi (t-r)}}  \\
	  &\quad \times \mathbf{E}_{(m(t)+x,i)}^h\left(\exp\left\{-c\Gamma_1(0,2r) - c\Gamma_2(2r, t-r) \right\}\frac{v_{I_{t-r}}(r,y)}{h_{I_{t-r}}};\mathcal{B}_{low}  \bigg| X_{t-r} =y \right) \mathrm{d} y\\
	&= e^{\lambda^* (t-r)}h_i\int_\R  \frac{e^{-\frac{(m(t)+x-y)^2}{2(t-r)}}}{\sqrt{2\pi (t-r)}} \times \mathbf{P}_{(m(t)+x,i)}^h\left(\mathcal{B}_{low} \Big|X_{t-r}=y  \right) \\
	&\quad \times \mathbf{E}_{(m(t)+x,i)}^h\left(\exp\left\{-c\Gamma_1(0,2r) - c\Gamma_2(2r, t-r) \right\}\frac{v_{I_{t-r}}(r,y)}{h_{I_{t-r}}}  \right) \mathrm{d} y.
\end{align}
Since when $r$ is large enough, $\mathbf{P}_{(x,i)}^h \left( I_{t-r} = \ell \right) \geq c'>0$ for all $t >8r$ and $i,\ell \in S$,  we get that
\begin{align*}
	& \mathbf{E}_{(m(t)+x,i)}^h\left(\exp\left\{-c\Gamma_1(0,2r) - c\Gamma_2(2r, t-r) \right\}\frac{v_{I_{t-r}}(r,y)}{h_{I_{t-r}}}  \right) \mathrm{d} y\\
	& =\mathbf{E}_{(m(t)+x,i)}^h\left(\left(1-\exp\left\{-c\sum_{s\in D_J, s\leq 2r} r^{\alpha_0} e^{-\alpha_0\sqrt{2\lambda^*}r} - c\sum_{s\in D_J, 2r < s\leq t-r} e^{- \lambda^* (s^\delta \land (t-s)^\delta)} \right\}\right)\frac{v_{I_{t-r}}(r,y)}{h_{I_{t-r}}}  \right) \\
	&= :  \mathbf{E}_{(m(t)+x,i)}^h\left(M(r,t)\frac{v_{I_{t-r}}(r,y)}{h_{I_{t-r}}}  \right)
	\leq  \sum_{j=1}^d \frac{v_{j}(r,y)}{h_{j}} \mathbf{E}_{(m(t)+x,i)}^h\left(M(r,t)\right)\\
	& \leq \frac{1}{c'} \mathbf{E}_{(m(t)+x,i)}^h\left( \frac{v_{I_{t-r}}(r,y)}{h_{I_{t-r}}}\right)\mathbf{E}_{(m(t)+x,i)}^h\left(M(r,t)\right).
\end{align*}
Therefore,
\begin{align*}
	& \mathbf{E}_{(m(t)+x,i)}^h\left(\exp\left\{-c\sum_{s\in D_J, s\leq 2r} r^{\alpha_0}e^{-\alpha_0\sqrt{2\lambda^*}r} - c\sum_{s\in D_J, 2r < s\leq t-r} e^{- \lambda^* (s^\delta \land (t-s)^\delta)} \right\}\frac{v_{I_{t-r}}(r,y)}{h_{I_{t-r}}}  \right)\\
	& =  \mathbf{E}_{(m(t)+x,i)}^h\left(\left(1-M(r,t)\right)\frac{v_{I_{t-r}}(r,y)}{h_{I_{t-r}}}  \right)\geq \mathbf{E}_{(m(t)+x,i)}^h\left(\frac{v_{I_{t-r}}(r,y)}{h_{I_{t-r}}}\right)\times \left(1-\frac{1}{c'}\mathbf{E}_{(m(t)+x,i)}^h\left(M(r,t)\right)\right).
\end{align*}
By the Markov property, we see that
\begin{align}\label{step_12}
	& \mathbf{E}_{(m(t)+x,i)}^h\left(\exp\left\{-c\sum_{s\in D_J, s\leq 2r} r^{\alpha_0} e^{-\alpha_0\sqrt{2\lambda^*}r} - c \sum_{s\in D_J, 2r < s\leq t-r} e^{- \lambda^* (s^\delta \land (t-s)^\delta)} \right\} \right)\nonumber \\
	& \geq \mathbf{E}_{(m(t)+x,i)}^h\left(\exp\left\{-c\sum_{s\in D_J, s\leq 2r} r^{\alpha_0}e^{-\alpha_0\sqrt{2\lambda^*}r}  	\right\}\right)\nonumber\\
	& \quad\times \inf_{j\in S} \mathbf{E}_{(m(t)+x,j)}^h\left(\exp\left\{- c \sum_{s\in D_J, s\leq - 2r+ t/2} e^{- \lambda^*  (s+2r)^\delta } \right\}\right)\nonumber \\
	& \quad \times \inf_{\ell \in S} \mathbf{E}_{(m(t)+x,\ell)}^h\left(\exp\left\{- c\sum_{s\in D_J,  s\leq t/2 -r}
	e^{- \lambda^*  (t/2-s)^\delta} \right\}\right).
\end{align}
By \eqref{Coupling2},  the product of the first two terms on the right-hand side of \eqref{step_12} is bounded from below by
\begin{align*}
	\exp\left\{ - \overline{a}\int_0^{2r} \left(1-e^{-cr^{\alpha_0} e^{-\alpha_0\sqrt{2\lambda^*}r }} \right)\mathrm{d} s - \overline{a} \int_{2r}^\infty \left(1-e^{-c e^{-\lambda^*  s^\delta}}\right) \mathrm{d} s\right\}= : F_1 (r).
\end{align*}
For the last term of the right-hand side of \eqref{step_12}, let $\lceil x \rceil$ be the smaller integer larger than $x$, then by the Markov property and \eqref{Coupling2},
\begin{align}
	&\mathbf{E}_{(m(t)+x,\ell)}^h\left(\exp\left\{- c\sum_{s\in D_J,  s\leq t/2 -r} e^{- \lambda^*  (t/2-s)^\delta} \right\}\right) \\ &\geq  \mathbf{E}_{(m(t)+x,\ell)}^h\left(\exp\left\{- \sum_{k=1}^{\lceil t/2 -r \rceil}c\sum_{s\in D_J,  k-1< s\leq k \land (t/2 -r)} e^{- \lambda^*  (t/2-k)^\delta} \right\}\right)\\
	& \geq \prod_{k=1}^{\lceil t/2 -r \rceil} \inf_{\ell_k \in S} \mathbf{E}_{(m(t)+x,\ell_k)}^h\left(\exp\left\{- c\sum_{s\in D_J,  s\leq 1\land (t/2 -r-k+1)} e^{- \lambda^*  (t/2-k)^\delta} \right\}\right) \\
	& \geq \prod_{k=1}^{\lceil t/2 -r \rceil} \exp\left\{- \underline{a}\int_0^{1\land (t/2 -r-k+1)}\left(1- e^{-ce^{- \lambda^*  (t/2-k)^\delta}}\right) \mathrm{d}s\right\}\\
	&= \prod_{k=1}^{\lceil t/2 -r \rceil} \exp\left\{- \underline{a}\int_{k-1}^{k\land (t/2 -r)}\left(1- e^{-ce^{- \lambda^*  (t/2-k)^\delta}}\right) \mathrm{d}s\right\}\\
	&\geq  \exp\left\{- \underline{a}\int_{0}^{t/2 -r}\left(1- e^{-ce^{- \lambda^*  (t/2-s+1)^\delta}}\right) \mathrm{d}s\right\} \geq \exp\left\{- \underline{a}\int_{r+1}^{\infty}\left(1- e^{-ce^{- \lambda^*  s^\delta}}\right) \mathrm{d}s\right\}=: F_2(r).
\end{align}
By the definition of $M(r,t)$, we conclude that for large $r$ and $t> 8r$,
\[
\mathbf{E}_{(m(t)+x,i)}^h\left(M(r,t)\right) \leq 1- F_1(r) F_2(r).
\]
Since $ \lim_{t\to\infty}\mathbf{P}_{(x,i)}^h(I_t=j)
= \lim_{r\to\infty} \inf_{t>r}\mathbf{P}_{(x,i)}^h(I_t=j)=g_j h_j$, we have
\begin{align}
	& \mathbf{E}_{(m(t)+x,i)}^h\left(\exp\left\{-c\sum_{s\in D_J, s\leq 2r} r^{\alpha_0}e^{-\alpha_0\sqrt{2\lambda^*}r} - c\sum_{s\in D_J, 2r < s\leq t-r} e^{- \lambda^*  s^\delta \land (t-s)^\delta} \right\}\frac{v_{I_{t-r}}(r,y)}{h_{I_{t-r}}}  \right)\\
	&\geq \inf_{j\in S} \inf_{t>r} \frac{\mathbf{P}_{(x,i)}^h(I_t =j)}{g_jh_j} \left(1-\frac{1}{c'}\left(1-F_1(r) F_2(r)\right)\right)  \sum_{j=1}^d g_j v_j(r,y )\\
	&= : C_{low}(r)\sum_{j=1}^d g_j v_j(r,y ).
\end{align}
It is easy to see that $C_{low}(r)\uparrow 1$ as $r\to\infty$. The proof is now complete.
\hfill$\Box$

\subsection{Proof of Theorem \ref{thm1}}

\textbf{Proof of Theorem \ref{thm1}: }
Define
\begin{align}
	\Psi_{low}^i(r; t,x)&:=  e^{\lambda^* (t-r)}h_i\int_\R  \frac{e^{-\frac{(m(t)+x-y)^2}{2(t-r)}}}{\sqrt{2\pi (t-r)}} \times \mathbf{P}_{(m(t)+x,i)}^h\left(\mathcal{B}_{low} \Big|X_{t-r}=y  \right)\sum_{j=1}^dg_j v_j(r,y)\mathrm{d} y,\\
	\Psi_{up}^i(r; t,x) &:= e^{\lambda^* (t-r)}h_i\int_{\R}\frac{e^{-\frac{(m(t)+x-y)^2}{2(t-r)}}}{\sqrt{2\pi (t-r)}} \times \mathbf{P}_{(m(t)+x,i)}^h\left(\mathcal{B}_{up} \Big|X_{t-r}=y  \right)\sum_{j=1}^d g_j v_j(r,y )\mathrm{d} y.
\end{align}
By Propositions \ref{prop2} and \ref{prop3},  for all $8r < x \leq \sqrt{t}$,
\begin{align}
	C_{up}(r) \Psi_{up}^i(r; t,x) \geq v_i(t,x) \geq C_{low}(r) \Psi_{low}^i(r; t,x)
\end{align}
with $C_{low}(r) \uparrow 1$, $C_{up}(r)\downarrow 1$.
Note that the proof of  \cite[Proposition 8.3 (c)]{Bramson1}
only uses probabilities of Brownian bridge. Using the same argument, we get that for all $i\in S$,
\begin{align}
	1\leq \frac{\Psi_{up}^i(r;t,x)}{\Psi_{low}^i(r;t,x)}\leq \gamma(r)
\end{align}
with $\gamma(r)\downarrow 1$ as $r\to\infty.$
Define
\begin{align}\label{Def-Psi-i}
	\Psi_{mid}^i(r;t,x):= e^{\lambda^* (t-r)}h_i\int_\R  \frac{e^{-\frac{(m(t)+x-y)^2}{2(t-r)}}}{\sqrt{2\pi (t-r)}}  \mathbf{P}_{(m(t)+x,i)}^h\left( \mathcal{B}_{mid} \Big|X_{t-r}=y  \right)\sum_{j=1}^d g_j v_j(r,y ) \mathrm{d} y
\end{align}
with
$$
\mathcal{B}_{mid}:= \left\{ X_s > n_{r,t}(t-s) \mbox{ for all }  s\in [0, t-r] \right\}
$$
and
$$
n_{r,t}(s) := m(t)\frac{s-r}{t-r}+ \sqrt{2\lambda^*}r \frac{t-s}{t-r}, \quad s\in [r,t].
$$
Then by \cite[(8.61), p. 144]{Bramson1} or \cite[(2.26)]{RSZ},
\[
\underline{\mathcal{M}}_{r,t}' (t-s)\leq n_{r,t}(t-s)\leq \overline{\mathcal{M}}_{r,t}^x(t-s),\quad \mbox{for }s\in [0, t-r]
\]
when $r$ and $x$ are large enough. This yields that (see \cite[(8.62), p. 144]{Bramson1})
\begin{align}\label{step_23}
	\frac{1}{\gamma(r)}\leq \frac{v_i(t,m(t)+x)}{\Psi_{mid}^i(r;t,x)}\leq \gamma(r)
\end{align}
for $r \geq r_1,  8r < x\leq  \sqrt{t}$ with $r_1$ fixed.
Therefore, to find the limit of $v_i(t,m(t)+x)$ as $t\to\infty$, we first  get the limit of $\Psi_{mid}^i(r;t,x)$ as $t\to\infty$.

{\bf Step 1: }\quad In this step we study the limit of $\Psi_{mid}^i(r;t,x)$ as $t\to\infty$. Letting $ y_1 = y -\sqrt{2\lambda^*}r$ and  $x_1= x - \frac{3}{2\sqrt{2\lambda^*}}\log t$, using \cite[Lemma 2.2 (a), p. 15]{Bramson1}, similar to \cite[(8.63) and (8.64), pp. 144--145]{Bramson1},
we have that
\begin{align}\label{Psi-i}
	 \Psi_{mid}^i(r;t,x) &=  e^{\lambda^* (t-r)}h_i\int_{\sqrt{2\lambda^*}r}^\infty  \frac{e^{-\frac{(m(t)+x-y)^2}{2(t-r)}}}{\sqrt{2\pi (t-r)}} \times \left(1- e^{-2 x y_1 /(t-r)} \right)\sum_{j=1}^d g_j v_j(r,y) \mathrm{d} y\nonumber \\
	& = e^{\lambda^* (t-r)}h_i\int_{0}^\infty  \frac{e^{-\frac{(x_1-y_1 + \sqrt{2\lambda^*}(t-r))^2}{2(t-r)}}}{\sqrt{2\pi (t-r)}} \times \left(1- e^{-2 x y_1 /(t-r)} \right)\sum_{j=1}^d g_j v_j(r,y_1 +\sqrt{2\lambda^*}r) \mathrm{d} y_1\nonumber\\
	& = h_i\int_{0}^\infty  \frac{e^{-\frac{(x_1-y_1)^2}{2(t-r)}}}{\sqrt{2\pi (t-r)}} e^{-\sqrt{2\lambda^*}(x_1 -y_1)} \times \left(1- e^{-2 x y_1 /(t-r)} \right)\sum_{j=1}^d  g_j v_j(r,y_1 +\sqrt{2\lambda^*}r) \mathrm{d} y_1\nonumber\\
	& =  \sqrt{\frac{t^3}{2\pi (t-r)^3}} e^{-\sqrt{2\lambda^*}x}h_i\int_{0}^\infty e^{\sqrt{2\lambda^*}y_1}\nonumber\\
	& \quad \times \left( \sum_{j=1}^d g_j v_j(r,y_1 +\sqrt{2\lambda^*}r)\right) e^{-\frac{(x_1-y_1)^2}{2(t-r)}}  \times(t-r) \left(1- e^{-2 x y_1 /(t-r)} \right)
	\mathrm{d} y_1\\
	&=:C_v^i(r; t,x) x e^{-\sqrt{2\lambda^*}x}.\nonumber
\end{align}
Therefore, for $r\geq r_1$ and $8r < x \leq \sqrt{t}$,
\begin{align}\label{step_15}
	\frac{1}{\gamma(r)} C_v^i(r;t,x) x e^{-\sqrt{2\lambda^*}x} \leq v_i(t,m(t)+x)\leq \gamma(r) C_v^i(r;t,x) x e^{-\sqrt{2\lambda^*}x} .
\end{align}

Now we fix $r$ first and replace $x$ by $x(t)$ and suppose that $x(t)\to x$ as $t\to\infty$. Let
    $x_1(t)= x(t)-\frac{3}{2\sqrt{2\lambda^*}}\log t$.
Then we can easily see that for any fixed $y_1$,
\[
e^{-\frac{(x_1(t)-y_1)^2}{2(t-r)}} (t-r) \left(1- e^{-2 x(t) y_1 /(t-r)} \right)  \to 2x y_1,\quad \mbox{as }t \to \infty,
\]
and that $e^{-\frac{(x_1(t)-y_1)^2}{2(t-r)}} (t-r) \left(1- e^{-2 x(t) y_1 /(t-r)} \right)  \leq  2x (t)y_1\leq \left(2\sup_t x(t) \right) y_1.$
Note that, for any $j\in S$,
by \eqref{Feynman-Kac-2} and Markov's inequality,
\begin{align}\label{step_68}
	&\int_{0}^\infty y_1 e^{\sqrt{2\lambda^*}y_1}v_j(r,y_1 +\sqrt{2\lambda^*}r) \mathrm{d} y_1 \leq \int_{0}^\infty y_1e^{\sqrt{2\lambda^*}y_1}  \times e^{\lambda^* r}h_j \mathbf{E}_{(y_1 +\sqrt{2\lambda^*}r,j)}^h\left(\frac{v_{I_r}(0, X_r)}{h_{I_r}}\right) \mathrm{d} y_1\nonumber\\
	& \lesssim e^{\lambda^* r} \int_0^\infty  y_1 e^{\sqrt{2\lambda^*}y_1} \mathbf{P}_{y_1+\sqrt{2\lambda^*}r} \left(X_r \leq N_2\right)\mathrm{d} y_1 \nonumber \\
	& = e^{\lambda^* r} \int_0^\infty y_1 e^{\sqrt{2\lambda^*}y_1} \mathbf{P}_{0} \left(X_r \leq N_2-y_1-\sqrt{2\lambda^*}r\right)\mathrm{d} y_1 \nonumber \\
	& \leq  e^{\lambda^* r} \int_0^\infty y_1 e^{\sqrt{2\lambda^*}y_1} \times e^{-2\sqrt{2\lambda^*} (y_1+\sqrt{2\lambda^*}r -N_2)} \mathbf{E}_{0}\left(e^{-2\sqrt{2\lambda^*}B_r}\right)\mathrm{d} y_1	< \infty.
\end{align}
Therefore, using the dominated convergence theorem, letting $t\to\infty$ in \eqref{Psi-i}, we get that when $x(t)\to x$,
\begin{align}\label{Conver-Psi-i}
	 \lim_{t\to\infty} \Psi_{mid}^i(r;t,x(t)) &= \sqrt{\frac{2}{ \pi }} x e^{-\sqrt{2\lambda^*}x}h_i\int_{0}^\infty y_1 e^{\sqrt{2\lambda^*}y_1}\ \left(\sum_{j=1}^d  g_j v_j(r,y_1 +\sqrt{2\lambda^*}r)\right)  \mathrm{d} y_1\nonumber \\
	 & =: h_iC_v(r) x e^{-\sqrt{2\lambda^*}x},
\end{align}
or equivalently, by the definition of $C_v^i(r; t,x)$,
\begin{align}\label{Conver-C-V-r}
	\lim_{t\to\infty} C_v^i(r; t,x(t)) \to h_iC_v(r),
\end{align}
where
\begin{align}\label{step_27}
	C_v(r):= \sqrt{\frac{2}{ \pi }} \int_{0}^\infty y_1 e^{\sqrt{2\lambda^*}y_1}\ \left(\sum_{j=1}^d  g_j v_j(r,y_1 +\sqrt{2\lambda^*}r)\right)  \mathrm{d} y_1.
\end{align}

{\bf Step 2: }\quad In this step we
use the limit of $\Psi_{mid}^i(r;t,x)$ as $t\to\infty$ to get the limit of  $v_i(t,m(t)+x)$ as $t\to\infty$.
It is easy to see that for any $r >0$, $C_v(r)\in (0,\infty).$ Indeed, $C_v(r)< \infty$ follows from \eqref{step_68}. On the other hand, by \eqref{Lower-Bound-M-T2},
\begin{align}
	&v_j(r, y_1 + \sqrt{2\lambda^*}r) \geq \P_{(0,j)}\left(M_r^{i_0} \geq y_1 + \sqrt{2\lambda^*}r -N_1 \right)>0,
\end{align}
which implies that $C_v(r)>0$.

Therefore, for any $r\geq r_1$ and $x(t) \geq 8r$ with $x(t)\to x$,  by \eqref{step_15} and \eqref{Conver-C-V-r},
we get that as $t\to\infty$,
\begin{align}\label{step_16}
	&\frac{1}{\gamma(r)}\leq \frac{ \liminf_{t\to\infty}v_i\left(t,m(t)+x(t)\right)}{h_i C_v (r) x e^{-\sqrt{2\lambda^*}x}}\leq \frac{\limsup_{t\to\infty}v_i\left(t,m(t)+x(t)\right)}{h_i C_v (r) x e^{-\sqrt{2\lambda^*}x}} \leq \gamma(r).
\end{align}
 Now letting $x\to\infty$ in \eqref{step_16}, we get
\begin{align}\label{step_17}
&	0 < \frac{C_v(r)}{\gamma(r)}\leq \liminf_{x\to \infty} \frac{ \liminf_{t\to\infty}v_i\left(t,m(t)+x(t)\right)}{h_i x e^{-\sqrt{2\lambda^*}x}}\leq \limsup_{x\to \infty} \frac{ \liminf_{t\to\infty}v_i\left(t,m(t)+x(t)\right)}{h_i x e^{-\sqrt{2\lambda^*}x}}\nonumber\\
& \leq \liminf_{x\to \infty} \frac{ \limsup_{t\to\infty}v_i\left(t,m(t)+x(t)\right)}{h_i x e^{-\sqrt{2\lambda^*}x}}\nonumber\\&\leq \limsup_{x\to \infty} \frac{ \limsup_{t\to\infty}v_i\left(t,m(t)+x(t)\right)}{h_i x e^{-\sqrt{2\lambda^*}x}} \leq C_v(r) \gamma(r) < \infty.
\end{align}
Letting $r \to\infty$, using the facts that $\gamma(r)\to 1$ and that
 the 4 quantities between $C_v(r)/\gamma(r) $ and $C_v(r)\gamma(r)$ in \eqref{step_17}
are independent of $r$,
we get
\begin{align}\label{step_19}
	&\lim_{r\to\infty} C_v(r) = \lim_{x\to \infty} \frac{ \liminf_{t\to\infty}v_i\left(t,m(t)+x(t)\right)}{h_i x e^{-\sqrt{2\lambda^*}x}}\nonumber \\& = \lim_{x\to \infty} \frac{ \limsup_{t\to\infty}v_i\left(t,m(t)+x(t)\right)}{h_i x e^{-\sqrt{2\lambda^*}x}} = : C_v(\infty)  \in (0, \infty),
\end{align}
which implies that \eqref{Cv(infty)} holds.
Now let $r=[x]/8$ and
\begin{align}
	Y_1(x)&:= \left(\frac{C_v(r)}{ C_v(\infty)} \vee \frac{ C_v(\infty)}{C_v(r)} \right)\gamma (r),\\
	Y_2(t; n)&:= 2 \left(\sup_{x>0} \left(xe^{-\sqrt{2\lambda^*}x}\right) \right)\times \sup_{i\in S}\sup \left\{ \left| C_v^i(r;t,x) - h_i C_v(r) \right|:\ x \in [n,n+1)	\right\}.
\end{align}
Let $ n_0\geq 8r_1$ be large enough such that $\gamma([x]/8)\leq 2$ for all $x\geq n_0$, then by \eqref{step_15},  for any $n \geq n_0$ and $x \in [n, n+1),$ we have
\begin{align}
	& v_i(t, m(t)+x)\leq  \gamma(r) C_v^i(r;t,x) x e^{-\sqrt{2\lambda^*}x}\\
	& = \gamma(r) \left(C_v^i(r;t,x) -h_i C_v(r) \right)x e^{-\sqrt{2\lambda^*}x}+ h_i \gamma(r) C_v(r) x e^{-\sqrt{2\lambda^*}x}\\
	& \leq h_i C_v(\infty) x e^{-\sqrt{2\lambda^*}x} Y_1(x)+Y_2(t;n)
\end{align}
and similarly,
\begin{align}
	v_i(t, m(t)+x)\geq h_i C_v(\infty) x e^{-\sqrt{2\lambda^*}x} \frac{1}{Y_1(x)}-Y_2(t;n).
\end{align}
Note that uniformly for all $i\in S$, it holds that
\begin{align}\label{step_18}
	Y_1(x)\to 1,\quad \mbox{as } x\to\infty,\quad \quad
	 Y_2(t;n)\to0 ,\quad \mbox{as } t\to\infty.
\end{align}
In conclusion, for all $n \geq n_0$, $x\in [n,n+1)$, $t \geq (n+1)^2$ and all $i\in S$,
\begin{align}\label{step_69}
	h_i C_v(\infty)x e^{-\sqrt{2\lambda^*}x} \frac{1}{Y_1(x)}-Y_2(t;n)\leq 	v_i(t, m(t)+x) \leq  h_i C_v(\infty) x e^{-\sqrt{2\lambda^*}x} Y_1(x)+Y_2(t;n)
\end{align}
with $Y_1(x), Y_2(t;n)$ satisfying \eqref{step_18}.

Fix $s>0 ,x  \in \R$, let
\[
x(t,s):= m(t+s)-m(t)+x -\sqrt{2\lambda^*}s= x  -\frac{3}{2\sqrt{2\lambda^*}}\log \left(1+\frac{s}{t}\right).
\]
Then $m(t+s)+x= m(t)+x(t,s)+\sqrt{2\lambda^*}s$ and
\begin{align}
	&1-v_i(t+s, m(t+s)+x)= \E_{(m(t+s)+x,i)}\left(1-\prod_{u\in Z(s)} \left(1- v_{I_u(s)}\left(t,X_u(s)\right)\right)\right)\\& = \E_{(0,i)} \left(1-\prod_{u\in Z(s)} \left(1- v_{I_u(s)}\left(t,X_u(s)+ \sqrt{2\lambda^*}s+m(t)+ x(t,s) \right)\right) \right).
\end{align}
 For any $1> \delta >0$, let $\varepsilon>0$ be sufficient small such that for all $x\in (0, \varepsilon)$,
 \[
 e^{-(1+\delta)x}\leq  1-x \leq e^{-x}.
 \]
 For this $ \varepsilon>0$, let $n_0$ be sufficient large such that
 when $x\geq n_0$, we have
 \[
 h_i C_v(\infty) x e^{-\sqrt{2\lambda^*}x}  < \frac{\varepsilon}{2},\quad \mbox{and}\quad \frac{1}{Y(x)},\  Y(x) \in [1-\delta, 1+\delta],
 \]
 which implies that
 \[
  h_i C_v(\infty) x e^{-\sqrt{2\lambda^*}x} \max\left\{\frac{1}{Y(x)}, Y(x)\right\}< \varepsilon.
 \]
 Thus, for any $\delta>0$, there exists $n_0$ such that when $x \geq n_0$,
 \begin{align}
 	1-  h_i C_v(\infty) x e^{-\sqrt{2\lambda^*}x} \frac{1}{Y_1 (x)}&\leq \exp\left\{ -  h_i C_v(\infty) x e^{-\sqrt{2\lambda^*}x} \frac{1}{Y_1 (x)} \right\}\\
 & \leq \exp\left\{ -  (1-\delta)h_i C_v(\infty) x e^{-\sqrt{2\lambda^*}x} \right\}, \\
 	1-  h_i C_v(\infty) x e^{-\sqrt{2\lambda^*}x} Y_1 (x)&\geq \exp\left\{- (1+\delta)  h_i C_v(\infty) x e^{-\sqrt{2\lambda^*}x} Y_1 (x)\right\}\\ &\geq \exp\left\{- (1+\delta)^2  h_i C_v(\infty) x e^{-\sqrt{2\lambda^*}x} \right\}.
 \end{align}
Since
$M_s^-+\sqrt{2\lambda^*}s:=\min_{u\in Z(s)} X_u(s) +\sqrt{2\lambda^*}s \to\infty$
(see \cite[Theorem 4]{RY}), there exists $s(\omega)$ such that for $t\geq s\geq s(\omega)$,
\begin{align*}
	\Delta_u(s,t):& =X_u(s)+\sqrt{2\lambda^*}s + x(t,s)\\
	& \geq X_u(s)+\sqrt{2\lambda^*}s + x-\frac{3}{2\sqrt{2\lambda^*}}\log 2 \geq  n_0,\quad \forall u\in Z(s).
\end{align*}
It follows from \eqref{step_69} that when $t\geq s\geq s(\omega)$,
\begin{align}
	&1- v_{I_u(s)}\left(t,\Delta_u(s,t) +m(t)\right)\geq 1- h_{I_u(s)}C_v(\infty)\Delta_u(s,t) e^{-\sqrt{2\lambda^*}\Delta_u(s,t)}{Y_1(\Delta_u(s,t))} -Y_2(t; [\Delta_u(s,t)])\\
	& \geq \exp\left\{ - (1+\delta)^2 h_{I_u(s)}C_v(\infty)\Delta_u(s,t) e^{-\sqrt{2\lambda^*}\Delta_u(s,t)}\right\}- Y_2(t; [\Delta_u(s,t)]),
\end{align}
and
\begin{align}
	&1- v_{I_u(s)}\left(t,\Delta_u(s,t) +m(t)\right)\leq 1- h_{I_u(s)}C_v(\infty)\Delta_u(s,t) e^{-\sqrt{2\lambda^*}\Delta_u(s,t)}Y_1(\Delta_u(s,t)) +Y_2(t; [\Delta_u(s,t)])\\
	& \leq \exp\left\{ -(1-\delta) h_{I_u(s)}C_v(\infty)\Delta_u(s,t) e^{-\sqrt{2\lambda^*}\Delta_u(s,t)}\right\}+ Y_2(t; [\Delta_u(s,t)]).
\end{align}
Therefore, on the event that $\left\{M_s^- +\sqrt{2\lambda^*}s + x(t,s)\geq n_0 \right\}$, it holds that
\begin{align}\label{step_20}
	&\prod_{u\in Z(s)} \left(\exp\left\{ - (1-\delta) h_{I_u(s)}C_v(\infty)\Delta_u(s,t) e^{-\sqrt{2\lambda^*}\Delta_u(s,t)}\right\}- Y_2(t; [\Delta_u(s,t)])\right)\nonumber \\
	&\geq \prod_{u\in Z(s)} \left(1- v_{I_u(s)}\left(t,X_u(s)+ \sqrt{2\lambda^*}s+m(t)+ x(t,s) \right)\right)\nonumber\\
	& \geq \prod_{u\in Z(s)} \left(\exp\left\{ - (1+\delta)^2 h_{I_u(s)}C_v(\infty)\Delta_u(s,t) e^{-\sqrt{2\lambda^*}\Delta_u(s,t)}\right\}- Y_2(t; [\Delta_u(s,t)])\right).
\end{align}
Since $x(s,t)\to x$ as $t\to\infty$, we have $\Delta_u(s,t)\to X_u(s)+\sqrt{2\lambda^*}s +x$ as $t\to\infty$.
Letting $t\to\infty$ in \eqref{step_20}, we get from \eqref{step_18} that
\begin{align}
	&\limsup_{t\to\infty} \left(1- v_i(t, m(t)+x)\right)= \limsup_{t\to\infty} \left(1- v_i(t+s, m(t+s)+x)\right)\\
	& \leq \P_{(0,i)}\left(M_s^- +\sqrt{2\lambda^*}s + x< n_0 \right)\\
	& \quad + \E_{(0,i)}\left(\exp\left\{ - (1-\delta) C_v(\infty)\left( x W_{\sqrt{2\lambda^*}}(s)+ M_{\sqrt{2\lambda^*}}(s)\right)e^{-\sqrt{2\lambda^*}x}\right\} 1_{\left\{M_s^- +\sqrt{2\lambda^*}s + x\geq n_0 \right\}}\right),
\end{align}
and
\begin{align}
	& \liminf_{t\to\infty} \left(1- v_i(t, m(t)+x)\right)\\
	&\geq \E_{(0,i)}\left(\exp\left\{ - (1+\delta)^2 C_v(\infty) \left(x W_{\sqrt{2\lambda^*}}(s)+ M_{\sqrt{2\lambda^*}}(s)\right)e^{-\sqrt{2\lambda^*}x}\right\} 1_{\left\{M_s^- +\sqrt{2\lambda^*}s + x\geq n_0 \right\}}\right),
\end{align}
where
$\{W_{\sqrt{2\lambda^*}}(s), s\ge 0\}$  is the additive martingale  defined by \eqref{Additive-martingale},  and $\{M_{\sqrt{2\lambda^*}}(s), s\ge 0\}$ is the derivative martingale defined by \eqref{Derivative-martingale}.
By \eqref{addMlimit} and \eqref{derivMlimit},
letting $s\to\infty$ and noting that
$\P_{(0,i)}\left(M_s^- +\sqrt{2\lambda^*}s + x\geq n_0 \right)\to 1$ for every fixed $x$ and $n_0$, we get that
\begin{align}
	&\P_{(0,i)}\left(\exp\left\{ - (1-\delta) C_v(\infty) M_{\sqrt{2\lambda^*}}(\infty)\right\} \right)\geq \limsup_{t\to\infty} \left(1- v_i(t, m(t)+x)\right)\\
	&\geq \liminf_{t\to\infty} \left(1- v_i(t, m(t)+x)\right)\\
	& \geq \P_{(0,i)}\left(\exp\left\{ - (1+\delta)^2 C_v(\infty) M_{\sqrt{2\lambda^*}}(\infty)\right\} \right).
\end{align}
Letting $\delta\to 0$, we get the desired convergence.
\hfill$\Box$

\section{Extremal Process for multitype branching Brownian motion}

In this section, we  study the asymptotic behavior of the extremal
 process of multitype branching Brownian motion
and prove Theorems \ref{Conver-Gap-Process} and  \ref{thm2}.

\begin{prop}\label{prop4}
	For any $\phi \in \mathcal{C}_c^+ (\R \times S)$ and $x\in \R$,
	\begin{align}
		\lim_{t\to\infty} \E_{(0,i)}\left(\exp\left\{- \int \phi(y+x,j) \mathcal{E}_t(\mathrm{d} y \mathrm{d}j) \right\} \right) = \E_{(0,i)} \left(\exp\left\{- C(\phi)M_{\sqrt{2\lambda^*}}(\infty) e^{-\sqrt{2\lambda^*}x}\right\} 		\right),
	\end{align}
where
	\begin{align}\label{C-phi}
		C(\phi) :=\lim_{r\to\infty} \sqrt{\frac{2}{ \pi }} \int_{0}^\infty y e^{\sqrt{2\lambda^*}y}\ \left(\sum_{j=1}^d  g_j v_j(r,y +\sqrt{2\lambda^*}r)\right)  \mathrm{d} y		 \in (0,\infty)
	\end{align}
	with $\mathbf{v}$ a solution of \eqref{F-KPP2} with initial value $v_j(0,y)= 1- e^{-\phi(-y,j)}$.
\end{prop}
\textbf{Proof: } For $L \in \R$, define
\begin{align}\label{step_30}
	v_j (0, y; L):= 1- \left(e^{-\phi(-y,j)}1_{\left\{ -y \leq L \right\}}\right),
\end{align}
then
\begin{align}\label{v-t-x-L}
	&v_i(t, x;L)= 1- \E_{(x,i)} \left( \prod_{u \in Z(t)}\left(1- v_{I_u(t)}(0, X_u(t); L) \right) \right)\nonumber \\
	&= 1-  \E_{(0,i)} \left(\exp\left\{-\sum_{u\in Z(t)} \phi\left(X_u(t)-x, I_u(t)\right) \right\} 1_{\left\{M_t \leq x+ L \right\}} \right).
\end{align}
For any fixed $L$,  $v_j(0, y; L)$ satisfies \eqref{Initial-Assumption}. Therefore, by Theorem \ref{thm1},
\begin{align}
		\lim_{t\to\infty} \left(1- v_i \left(t, m(t)+x ; L\right) \right) =\P_{(0,i)}\left(\exp\left\{ -  C(\phi; L) M_{\sqrt{2\lambda^*}}(\infty)e^{-\sqrt{2\lambda^*}x}\right\}  \right),
\end{align}
with $C(\phi; L)$ defined by
\begin{align}\label{C-phi-L}
	C(\phi; L) = \lim_{r\to\infty} \sqrt{\frac{2}{ \pi }} \int_{0}^\infty y e^{\sqrt{2\lambda^*}y}\ \left(\sum_{j=1}^d g_j v_j(r,y +\sqrt{2\lambda^*}r; L)\right)  \mathrm{d} y.
\end{align}
Since
\begin{align}\label{step_31}
	& 0\leq  v_j (t,x;L) - v_j(t,x)\\
	&  =\E_{(0,j)} \left(\exp\left\{-\sum_{u\in Z(t)} \phi\left(X_u(t)-x, I_u(t)\right) \right\} 1_{\left\{M_t > x+ L \right\}} \right)\leq \P_{(0,j)}\left(M_t >x +L \right),
\end{align}
we get that
\begin{align*}
H_1(r;L) - H_2(r;L)
&:=	\sqrt{\frac{2}{ \pi }} \int_{0}^\infty y e^{\sqrt{2\lambda^*}y}\ \left(\sum_{j=1}^d g_j v_j(r,y +\sqrt{2\lambda^*}r;L)\right)  \mathrm{d} y \\
	& \quad - \sqrt{\frac{2}{ \pi}} \int_{0}^\infty y e^{\sqrt{2\lambda^*}y}\ \left(\sum_{j=1}^d g_j \P_{(0,j)}\left(M_r > y+\sqrt{2\lambda^*}r +L \right)\right)  \mathrm{d} y
\\	& \leq \sqrt{\frac{2}{\pi}} \int_{0}^\infty y e^{\sqrt{2\lambda^*}y}\ \left(\sum_{j=1}^d g_j v_j(r,y +\sqrt{2\lambda^*}r)\right)  \mathrm{d} y = : H(r)\\
	& \leq \sqrt{\frac{2}{ \pi }} \int_{0}^\infty y e^{\sqrt{2\lambda^*}y}\ \left(\sum_{j=1}^d g_j v_j(r,y +\sqrt{2\lambda^*}r;L)\right)  \mathrm{d} y
=H_1(r;L).
\end{align*}
Therefore,
\begin{equation}\label{ineq-H}
H_1(r;L) - H_2(r;L)\le H(r)\le H_1(r;L).
\end{equation}
Note that
\[
H_2(r;L)\leq e^{-\sqrt{2\lambda^*}L} \sqrt{\frac{2}{ \pi }} \int_{0}^\infty y e^{\sqrt{2\lambda^*}y}\ \left(\sum_{j=1}^d g_j \P_{(0,j)}\left(M_r > y+\sqrt{2\lambda^*}r  \right)\right)  \mathrm{d} y\stackrel{r\to\infty}{\longrightarrow} e^{-\sqrt{2\lambda^*}L}C_\infty
\]
with $C_\infty$ given in Corollary \ref{cor4}. Thus
\begin{align}
\lim_{L\to +\infty} \limsup_{r\to\infty}	H_2(r;L)=0.
\end{align}
Also note that $\lim_{r\to\infty}H_1(r;L)=C(\phi; L).$
Since $C(\phi;L)$ is positive and decreasing in $L$,
letting $r\to\infty$ and then $L\to \infty$ in \eqref{ineq-H},
we have
\begin{align}
	C(\phi):= \lim_{L\to\infty} C(\phi; L) \leq \liminf_{r\to\infty} H(r)\leq \limsup_{r\to\infty} H(r) \leq \lim_{L\to\infty} C(\phi; L),
\end{align}
which implies that $\lim_{r\to\infty} H(r)= C(\phi).$

Next, for any $\phi\in \mathcal{C}_c^+(\R)$ with $\phi \neq 0$,
there exist $\ell_0\in S$,  $a_{\ell_0}< b_{\ell_0}$ and $c_0>0$ such that $\phi(y,\ell_0)\geq c_0$ for all $y\in [a_{\ell_0},b_{\ell_0}]$.
Thus,
\begin{align}
	\E_{(0,i)}\left(\exp\left\{- \int \phi(y,j) \mathcal{E}_t(\mathrm{d} y \mathrm{d}j) \right\} \right)
	&\leq  e^{-c_0}\P_{(0,i)}\left( M_t^{\ell_0}-m(t)\in [a_{\ell_0}, b_{\ell_0}) \right)\\
	&\quad + \left(1-\P_{(0,i)}\left( M_t^{\ell_0}-m(t)\in [a_{\ell_0}, b_{\ell_0}) \right)\right).
\end{align}
It follows immediately from  Corollary \ref{cor3} that
\begin{align}\label{step_46}
	\lim_{t\to\infty}\P_{(0,i)}\left( M_t^{\ell_0}-m(t)\in [a_{\ell_0}, b_{\ell_0}) \right)>0.
\end{align}
Thus $C(\phi)>0$. Hence we have shown that $C(\phi)>0$ when $\phi\in \mathcal{C}_c^+(\R \times S)$ and $\phi \neq 0$.

For any $x\in \R$, when $L$ is large enough so that $x+L\geq 1$, by Proposition \ref{Upper-Bound-M-T}, there exists a constant $C_0$  such that
\begin{align}
	& 1-v_i\left(t,m(t)+x; L\right)\leq
	1-v_i\left(t,m(t)+x\right) \\
	&\leq 1-v_i\left(t,m(t)+x;L\right) + \P_{(0,i)}\left(M_t > m(t)+x +L\right)\\
	&\leq 1-v_i\left(t,m(t)+x;L\right) + C_0 (x+L)e^{-\sqrt{2\lambda^*}(x+L)}.
\end{align}
Letting $t \to\infty$, we get
\begin{align}
	&\P_{(0,i)}\left(\exp\left\{ -  C(\phi; L) M_{\sqrt{2\lambda^*}}(\infty)e^{-\sqrt{2\lambda^*}x}\right\}  \right) \leq \liminf_{t\to\infty} \left(1-v_i\left(t,m(t)+x\right)\right)\\
	&\leq \limsup_{t\to\infty} \left(1-v_i \left(t,m(t)+x\right)\right) \\
	& \leq \P_{(0,i)}\left(\exp\left\{ -  C(\phi; L) M_{\sqrt{2\lambda^*}}(\infty)e^{-\sqrt{2\lambda^*}x}\right\}  \right)+C_0 (x+L)e^{-\sqrt{2\lambda^*}(x+L)}.
\end{align}
Next, letting $L\to\infty$, we get the desired result.
\hfill$\Box$

\begin{cor}\label{cor2}
	The point process $\mathcal{E}_t$ converges
	in distribution to a random measure $\mathcal{E}_\infty$,
	where the Laplace transform of $\mathcal{E}_\infty$ is given by
	\[
	\E_{(0,i)}\left(\exp\left\{- \int \phi(y+x,j) \mathcal{E}_\infty(\mathrm{d} y \mathrm{d}j) \right\} \right) = \E_{(0,i)} \left(\exp\left\{- C(\phi)M_{\sqrt{2\lambda^*}}(\infty) e^{-\sqrt{2\lambda^*}x}\right\} \right)
	\]
	with $C(\phi)$ given in \eqref{C-phi}.
\end{cor}
\textbf{Proof: } Without loss of generality, we assume $x=0$, otherwise we may consider $\widehat{\phi}(\cdot, j)= \phi(x+\cdot,j)$. It suffices to prove the tightness for $\mathcal{E}_t$, which is equivalent to the tightness for $\int \phi(y,j)\mathcal{E}_t(\mathrm{d}y \mathrm{d}j)$.
By Proposition \ref{prop4}, it suffices to show that $\lim_{\theta\downarrow 0} C(\theta \phi)= 0.$
Choose $m_\phi$ so that $\phi(y,j)=0$ for  all $y< m_\phi$ and $j\in S$. Let $|\Vert \phi \Vert_\infty:= \sup_{x\in\R, j\in S} |\phi(x,j)|$, then
\begin{align}
	&\E_{(0,i)}\left(\exp\left\{- \int\theta \phi(y,j) \mathcal{E}_t(\mathrm{d} y \mathrm{d}j) \right\} \right)\geq
	\E_{(0,i)}\left(\exp\left\{-\theta \Vert \phi \Vert_\infty \mathcal{E}_t\left((m_\phi,\infty)\times S \right) \right\}\right)\\
	&\geq e^{-\theta \Vert \phi \Vert_\infty N}
	\P_{(0,i)}\left(\mathcal{E}_t\left((m_\phi,\infty)\times S \right)\leq N \right).
\end{align}
First letting $t\to\infty$, next $\theta \to 0$ and then $N\to\infty$, we only need to prove that for all $i \in S$,
\begin{align}\label{step_35}
	\lim_{N\to\infty} \limsup_{t\to\infty} \P_{(0,i)}\left(\mathcal{E}_t\left((m_\phi,\infty)\times S \right)> N \right) =0.
\end{align}
Suppose that under $\{ \mathcal{E}_t\left((m_\phi,\infty)\times S \right)> N\}$,  $u_1,...,u_N \in \left\{u\in Z(t), X_u(t)-m(t)> m_\phi \right\}$, then
\begin{align}
	&\P_{(0,i)}\left(\mathcal{E}_t\left((m_\phi,\infty)\times S \right)> N, M_{t+1}-m(t+1)\leq n \right) \\
	& \leq \P_{(0,i)}\left(\mathcal{E}_t\left((m_\phi,\infty)\times S \right)> N, \max_{1\leq n\leq N} M_1^{u_n} + X_{u_n}(t) -m(t+1)\leq n \right) \\
	& \leq \P_{(0,i)}\left(\mathcal{E}_t\left((m_\phi,\infty)\times S \right)> N, \max_{1\leq n\leq N} M_1^{u_n} + m(t)+m_\phi -m(t+1)\leq n \right)\\
	& \leq \left(\sup_{j\in S} \P_{(0,j)}\left(M_1 + m(t)+m_\phi -m(t+1)\leq n\right)\right)^N.
\end{align}
By Proposition \ref{Upper-Bound-M-T}, we have
\begin{align}\label{step_34}
	&\P_{(0,i)}\left(\mathcal{E}_t\left((m_\phi,\infty)\times S \right)> N \right)
\nonumber\\
	&\leq \P_{(0,i)}\left(M_{t+1}-m(t+1)> n \right)+ \left(\sup_{j\in S} \P_{(0,j)}\left(M_1 + m(t)+m_\phi -m(t+1)\leq n\right)\right)^N
\nonumber\\
	& \lesssim ne^{-\sqrt{2\lambda^*}n}+ \left( \mathbf{P}_0\left(B_1 + m(t)+m_\phi -m(t+1)\leq n\right)\right)^N,\quad n,t, N\geq 1.
\end{align}
For every $n\geq 1,$ letting $t\to\infty$ first and then $N\to\infty$ in \eqref{step_34}, we get that
\[
\lim_{N\to\infty} \limsup_{t\to\infty} \P_{(0,i)}\left(\mathcal{E}_t\left((m_\phi,\infty)\times S \right)> N \right)\lesssim ne^{-\sqrt{2\lambda^*}n}.
\]
Letting $n\to\infty$, we get \eqref{step_35} and thus $\mathcal{E}_t$ converges in distribution to a random point process $\mathcal{E}_\infty$.

\hfill$\Box$

Recall the definition of $\Psi_{mid}^i(r;t,x)$ defined in \eqref{Def-Psi-i}. By \eqref{Psi-i}, we have
\begin{align}
	&\Psi_*^i(r;t,x) := \Psi_{mid}^i\left(r;t,x+\frac{3}{2\sqrt{2\lambda^*}}\log t \right) =
	\sqrt{\frac{t^3}{2\pi (t-r)^3}} e^{-\sqrt{2\lambda^*}\left(x + \frac{3}{2\sqrt{2\lambda^*}}\log t\right)}h_i\int_{0}^\infty e^{\sqrt{2\lambda^*}y}\\
	& \quad \times \left( \sum_{j=1}^d g_j v_j(r,y +\sqrt{2\lambda^*}r)\right) e^{-\frac{(x-y)^2}{2(t-r)}}  \times(t-r) \left(1- e^{-2 \left(x +\frac{3}{2\sqrt{2\lambda^*}}\log t\right)y /(t-r)} \right) \mathrm{d} y\\
	& = \frac{1}{\sqrt{2\pi (t-r)}} e^{-\sqrt{2\lambda^*}x}h_i\int_{0}^\infty e^{\sqrt{2\lambda^*}y}\\
	& \quad \times \left( \sum_{j=1}^d g_j v_j(r,y +\sqrt{2\lambda^*}r)\right) e^{-\frac{(x-y)^2}{2(t-r)}}  \left(1- e^{-2 \left(x +\frac{3}{2\sqrt{2\lambda^*}}\log t\right)y /(t-r)} \right) \mathrm{d} y.
\end{align}
It follows from \eqref{step_23} that
\begin{align}\label{step_28}
\frac{1}{\gamma(r)}\leq \frac{v_i(t,\sqrt{2\lambda^*}t+x)}{\Psi_*^i(r;t,x)}\leq \gamma(r)
\end{align}
holds for $r \geq r_1,  8r-\frac{3}{2\sqrt{2\lambda^*}}\log t < x\leq \sqrt{t}- \frac{3}{2\sqrt{2\lambda^*}}\log t$ with  $r_1$ fixed.

\begin{lemma}\label{lemma6}
	Let $\mathbf{v}$ solve \eqref{F-KPP2}
with initial value satisfying \eqref{Initial-Assumption}.
Then for any fixed $x\in \R$,
	\begin{align}
		\lim_{t\to\infty} \frac{t^{3/2}}{\frac{3}{2\sqrt{2\lambda^*}}\log t}\Psi_*^i(r;t,x) = h_i C_v(r)e^{-\sqrt{2\lambda^*}x},
	\end{align}
where  $C_v(r)$ is given by \eqref{step_27}.
\end{lemma}
\textbf{Proof: } The proof is very similar to that of \cite[Lemma 4.5]{ABK} and we omit the details.

\hfill$\Box$

Let $\mathbf{v}$ solve \eqref{F-KPP2} with initial value satisfying \eqref{Initial-Assumption}. By Lemma \ref{lemma6} and \eqref{step_28},  we have for every $x\in \R$,
	\begin{align}\label{step_29}
	\lim_{t\to \infty} \frac{t^{3/2}}{\frac{3}{2\sqrt{2\lambda^*}}\log t}	v_i(t,\sqrt{2\lambda^*}t +x) = h_i C_v(\infty) e^{-\sqrt{2\lambda^*}x},
	\end{align}
where  $C_v(\infty)=\lim_{r\to\infty} C_v(r)$, given by \eqref{Cv(infty)}.
Hence for every $x\in \R$,
\begin{align}\label{step_42}
	\lim_{t\to \infty} \frac{t^{3/2}}{\frac{3}{2\sqrt{2\lambda^*}}\log t}	\P_{(0,i)}\left(M_t> \sqrt{2\lambda^*}t +x\right)= h_i C_\infty e^{-\sqrt{2\lambda^*}x},
\end{align}
with $C_\infty$ defined in \eqref{C-Infinity}.
Now we extend \eqref{step_29} to the case $v_i(0,y)=1-e^{-\phi(-y,i)}$:
\begin{lemma}\label{lemma7}
	For any $\phi \in \mathcal{C}_c^+ (\R \times S)$,
let $\mathbf{v}$ solve \eqref{F-KPP2}
with $v_i(0,y)=1-e^{-\phi(-y,i)}$, then
	\begin{align}
		\lim_{t\to \infty} \frac{t^{3/2}}{\frac{3}{2\sqrt{2\lambda^*}}\log t}	v_i(t,\sqrt{2\lambda^*}t +x) = h_i C(\phi )e^{-\sqrt{2\lambda^*}x}
	\end{align}
	with $C(\phi)$ defined in \eqref{C-phi}.
\end{lemma}
\textbf{Proof: } Let $v_i(t,x;L)$ solves \eqref{F-KPP2} with initial value \eqref{step_30}.
By \eqref{step_29},
\begin{align}
	\lim_{t\to \infty} \frac{t^{3/2}}{\frac{3}{2\sqrt{2\lambda^*}}\log t}	v_i(t,\sqrt{2\lambda^*}t +x;L) = h_i C(\phi;L) e^{-\sqrt{2\lambda^*}x}
\end{align}
with $C(\phi; L)$ defined in \eqref{C-phi-L}.  By \eqref{step_31} and Proposition \ref{Upper-Bound-M-T}, we have
\begin{align}
	&0\leq  v_i (t,\sqrt{2\lambda^*}t +x;L) - v_i(t,\sqrt{2\lambda^*}t +x)\\
	&\leq \P_{(0,i)}\left(M_t >\sqrt{2\lambda^*}t+x +L \right)\lesssim \left(\frac{3}{2\sqrt{2\lambda^*}}\log t+x+L\right)e^{-\sqrt{2\lambda^*}\left(\frac{3}{2\sqrt{2\lambda^*}}\log t+x+L\right)}.
\end{align}
Note that $C(\phi)= \lim_{L\to\infty} C(\phi; L)$.
Letting $t\to\infty$ first and then $L\to +\infty$, we arrive at the desired conclusion.

\hfill$\Box$

Define
\[
\overline{\mathcal{E}}_t:= \sum_{u\in Z(t)} \delta_{(X_u(t)-\sqrt{2\lambda^*}t, I_u(t))},\quad \overline{\mathcal{E}}_t -z := \sum_{u\in Z(t)} \delta_{(X_u(t)-\sqrt{2\lambda^*}t -z, I_u(t))},
\]

\begin{prop}\label{prop5}
	For any $z\in \R$ and $i\in S$, under $\mathbb{P}_{(0,i)}\left(\cdot \big| M_t>\sqrt{2\lambda^*}t +z\right),$
	 $\left(\overline{\mathcal{E}}_t-z, M_t-\sqrt{2\lambda^*}t-z\right)$ converges in distribution to a limit $(\overline{\mathcal{E}}_\infty, Y)$
	 independent of $z$ and $i$,
	 where $\overline{\mathcal{E}}_\infty$ is a point process, $Y$ is an exponential random variable with parameter $\sqrt{2\lambda^*}$ and
	\begin{align}\label{Joint-law}
		\E_{(0,i)}\left(\exp\left\{-\int \phi(y,j)\overline{\mathcal{E}}_\infty (\mathrm{d} y\mathrm{d} j) \right\}; Y>x \right)=
\frac{\widetilde C(\phi,x)}{C_\infty} - \frac{C(\phi)}{C_\infty},
	\end{align}
where $C(\phi)$ is given by \eqref{C-phi}, $C_\infty$ is given by \eqref{Cv(infty)}, and
\begin{equation}\label{tildeC(phi, x)}
\widetilde{C}(\phi, x):= \lim_{r\to\infty} \sqrt{\frac{2}{ \pi }} \int_{0}^\infty y e^{\sqrt{2\lambda^*}y}\ \left(\sum_{j=1}^d  g_j (v_1)_j(r,y +\sqrt{2\lambda^*}r)\right)  \mathrm{d} y
\end{equation}
with  $\mathbf{v}_1$ being a solution of  \eqref{F-KPP2} with initial value
\begin{align}\label{step_36}
	(v_1)_i (0, y)= 1- e^{-\phi(-y,i)} 1_{\{-y \leq x \}},\quad i\in S.
\end{align}
\end{prop}
\textbf{Proof: } By \eqref{step_42}, for any $x>0$,
\begin{align}\label{Conditioned-tightness1}
	&\lim_{t\to\infty} \P_{(0,i)}\left(M_t-\sqrt{2\lambda^*}t-z >x \big| M_t>\sqrt{2\lambda^*}t +z \right) \nonumber\\
	&=\lim_{t\to\infty}\frac{\frac{t^{3/2}}{\frac{3}{2\sqrt{2\lambda^*}}\log t}	\P_{(0,i)}\left(M_t > \sqrt{2\lambda^*}t +x+z\right)}{\frac{t^{3/2}}{\frac{3}{2\sqrt{2\lambda^*}}\log t}	\P_{(0,i)}\left(M_t > \sqrt{2\lambda^*}t +z\right)}=e^{-\sqrt{2\lambda^*}x}.
\end{align}
Then we have  under $\mathbb{P}_{(0,i)}\left(\cdot \big| M_t>\sqrt{2\lambda^*}t +z\right),$
	 $ M_t-\sqrt{2\lambda^*}t-z$ converges in distribution to $Y$,  an exponential random variable with parameter $\sqrt{2\lambda^*}$.
For any $\phi\in \mathcal{C}_c^+\left(\R \times S \right)$ and $x>0$,
\begin{align}\label{step_32}
	& \E_{(0,i)}\left(\exp\left\{ -\int \phi (y,j)\left( \overline{\mathcal{E}}_t -z\right)(\mathrm{d}y \mathrm{d}j) \right\}; M_t> \sqrt{2\lambda^*}t +z+x \big|  M_t> \sqrt{2\lambda^*}t +z \right) \\
	& = \frac{1}{\P_{(0,i)}\left( M_t> \sqrt{2\lambda^*}t +z \right)} \E_{(0,i)}\left(\prod_{u\in Z(t)} e^{-\phi(X_u(t)-\sqrt{2\lambda^*}t -z, I_u(t))}; M_t> \sqrt{2\lambda^*}t +z+x  \right)\\
	& = \frac{1}{\P_{(0,i)}\left( M_t> \sqrt{2\lambda^*}t +z \right)} \left(1-\E_{(0,i)}\left(\prod_{u\in Z(t)} e^{-\phi(X_u(t)-\sqrt{2\lambda^*}t -z, I_u(t))}; M_t\leq \sqrt{2\lambda^*}t +z+x  \right)\right) \\
	& \quad -  \frac{1}{\P_{(0,i)}\left( M_t> \sqrt{2\lambda^*}t +z \right)} \E_{(0,i)}\left(1-\prod_{u\in Z(t)} e^{-\phi(X_u(t)-\sqrt{2\lambda^*}t -z, I_u(t))}  \right) \\
	&=: \frac{1}{\P_{(0,i)}\left( M_t> \sqrt{2\lambda^*}t +z \right)} (v_1)_i(t,\sqrt{2\lambda^*}t+z) - \frac{1}{\P_{(0,i)}\left( M_t> \sqrt{2\lambda^*}t +z \right) }(v_2)_i(t,\sqrt{2\lambda^*}t+z)
\end{align}
where $\mathbf{v}_1$ and $\mathbf{v}_2$ solve \eqref{F-KPP2} with
\begin{align}
	\label{Initial-v1-v2}
	(v_1)_i (0, y)= 1- e^{-\phi(-y,i)} 1_{\{-y \leq x \}},\quad (v_2)_i (0,y) = 1-e^{-\phi(-y,i) },\quad i\in S
\end{align}
according to \eqref{v-t-x-L}.
Using Lemma \ref{lemma7} and \eqref{step_42}, it is easy to see that
\begin{align}\label{Conditioned-Tightness3}
	&\lim_{t\to\infty}\frac{1}{\P_{(0,i)}\left( M_t> \sqrt{2\lambda^*}t +z \right)} (v_1)_i(t,\sqrt{2\lambda^*}t+z) - \lim_{t\to\infty}\frac{1}{\P_{(0,i)}\left( M_t> \sqrt{2\lambda^*}t +z \right) }(v_2)_i(t,\sqrt{2\lambda^*}t+z) \nonumber \\
	& = \frac{\widetilde{C}(\phi,x)}{C_\infty} - \frac{C(\phi)}{C_\infty},
\end{align}
with $\widetilde{C}(\phi, x)$ being defined by \eqref{tildeC(phi, x)},
and the right-hand side of \eqref{Conditioned-Tightness3} is independent of $z\in \R$ and $i\in S$.
Let $x=0$ in \eqref{Conditioned-Tightness3}, then
\begin{align}\label{step_33}
	\lim_{t\to\infty} \E_{(0,i)}\left(\exp\left\{ -\int \phi (y,j) \left(\overline{\mathcal{E}}_t -z\right)(\mathrm{d}y \mathrm{d}j) \right\}\bigg|  M_t> \sqrt{2\lambda^*}t +z \right) = \frac{\widetilde{C}(\phi,0)}{C_\infty} - \frac{C(\phi)}{C_\infty}.
\end{align}
 Note that $\left(\overline{\mathcal{E}}_t -z\right)$ under $\P_{(0,i)}\left(\cdot| M_t> \sqrt{2\lambda^*}t +z\right)$ is still a  point process.
We now prove the convergence of $\left(\overline{\mathcal{E}}_t -z\right)$ in distribution
 under $\P_{(0,i)}\left(\cdot| M_t> \sqrt{2\lambda^*}t +z\right)$. By \eqref{step_33}, it suffices to prove that
\begin{align}\label{Conditioned-Tigheness2}
	\lim_{\theta\downarrow0} \left(\frac{\widetilde{C}(\theta\phi,0)}{C_\infty} - \frac{C(\theta\phi)}{C_\infty} \right)=1.
\end{align}
By Corollary \ref{cor2}, we have $\lim_{\theta\downarrow 0} C(\theta\phi)=0$.
Note that the initial value of $\mathbf{v}_1$  in \eqref{Initial-v1-v2} with $x=0$ satisfies condition \eqref{Initial-Assumption}, it follows from Theorem \ref{thm1} that
\begin{align}
	&\lim_{t\to\infty} \E_{(0,i)}\left(\prod_{u\in Z(t)} e^{-\theta \phi(X_u(t)-m(t) -z, I_u(t))}; M_t\leq m(t) +z \right)\\
	&= \E_{(0,i)}\left(\exp\left\{ - \widetilde{ C}(\theta \phi, 0) M_{\sqrt{2\lambda^*}}(\infty) e^{-\sqrt{2\lambda^*}z}\right\}  \right)
\end{align}
and by Corollary \ref{cor4},
\begin{align}
	\lim_{t\to\infty} \P_{(0,i)}\left( M_t\leq m(t) +z \right)= \E_{(0,i)}\left(\exp\left\{ -  C_\infty  M_{\sqrt{2\lambda^*}}(\infty) e^{-\sqrt{2\lambda^*}z}\right\}  \right).
\end{align}
Also note that
\begin{align}
&\left|	\P_{(0,i)}\left( M_t\leq m(t) +z \right) - \E_{(0,i)}\left(\prod_{u\in Z(t)} e^{-\theta \phi(X_u(t)-m(t) -z, I_u(t))}; M_t\leq m(t) +z \right)\right| \\
	& \leq 1- \E_{(0,i)}\left(\exp\left\{- \theta \int \phi(y-z,j) \mathcal{E}_t(\mathrm{d} y \mathrm{d}j) \right\} \right).
\end{align}
Letting $t\to\infty$, we get that
\begin{align}
	&\left|\E_{(0,i)}\left(\exp\left\{ -  C_\infty  M_{\sqrt{2\lambda^*}}(\infty) e^{-\sqrt{2\lambda^*}z}\right\}  \right) - \E_{(0,i)}\left(\exp\left\{ -  \widetilde{C}(\theta \phi, 0) M_{\sqrt{2\lambda^*}}(\infty) e^{-\sqrt{2\lambda^*}z}\right\}  \right)\right| \\
	&\leq 1- \E_{(0,i)} \left(\exp\left\{- C(\theta \phi)M_{\sqrt{2\lambda^*}}(\infty) e^{\sqrt{2\lambda^*}z}\right\} \right).
\end{align}
Let $\theta \downarrow 0$, we get that $\lim_{\theta\downarrow 0} \widetilde{C}(\theta\phi,0) = C_\infty$, which implies \eqref{Conditioned-Tigheness2}.
Combining \eqref{Conditioned-tightness1},
 \eqref{Conditioned-Tightness3}, \eqref{Conditioned-Tigheness2}
and the fact
that the process $(X_t, Y_t)$ is tight if $X_t$ and $Y_t$ are both tight, which follows from the inequality
\[
\inf_{t>0} \mathbb{P}\left(|X_t|\leq K, |Y_t|\leq K \right)\geq \inf_{t>0} \mathbb{P} \left(|X_t| \leq K\right)+ \inf_{t>0} \mathbb{P} \left(|Y_t| \leq K\right)-1,
\]
we get that under $\P_{(0,i)},$
\begin{align}
	\left(\overline{\mathcal{E}}_t-z, M_t -\sqrt{2\lambda^*}t-z \right)|_{M_t>\sqrt{2\lambda^*}t +z}
\end{align}
converges joint in distribution to $(\overline{\mathcal{E}}_\infty, Y)$, where the joint law is given in \eqref{Joint-law} and is independent of $z\in \R$ and $i\in S$.

\hfill$\Box$

 \textbf{Proof of Theorem \ref{Conver-Gap-Process}: }
Define
$\mathcal{D}:=\overline{\mathcal{E}}_\infty -Y$.
By Proposition \ref{prop5} and \cite[Lemma 4.13]{ABK}, also note that $\mathcal{D}_t = \left(\overline{\mathcal{E}}_t-z \right)-  \left(M_t-\sqrt{2\lambda^*}t-z\right)$, we get that
under $\mathbb{P}_{(0,i)}\left(\cdot \big| M_t>\sqrt{2\lambda^*}+z\right)$,
$\mathcal{D}_t$ converges in distribution to $\mathcal{D}$. Also, for all $x>0$,
\begin{align}
	&\E_{(0,i)}\left(\exp\left\{-\int \phi(y,j)\mathcal{D}_t(\mathrm{d}y \mathrm{d}j) \right\}; M_t> \sqrt{2\lambda^*}t +z+x \big|  M_t> \sqrt{2\lambda^*}t +z\right)\\
	& = \E_{(0,i)}\left(\exp\left\{-\int \phi(y-M_t+\sqrt{2\lambda^*}t +z+x,j)\overline{\mathcal{E}}_t(\mathrm{d}y \mathrm{d}j) \right\}\big|  M_t> \sqrt{2\lambda^*}t +z+x\right)\\
	& \quad \times \P_{(0,i)}\left( M_t> \sqrt{2\lambda^*}t +z+x \big|  M_t> \sqrt{2\lambda^*}t +z\right)\\
	& \to \E_{(0,i)}\left(\exp\left\{-\int \phi(y-Y,j)\overline{\mathcal{E}}_\infty(\mathrm{d}y \mathrm{d}j) \right\}\right)\P_{(0,i)}\left(Y>x\right)\\
	& = \E_{(0,i)}\left(\exp\left\{-\int \phi(y,j)\mathcal{D}(\mathrm{d}y \mathrm{d}j) \right\}\right)\P_{(0,i)}\left(Y>x\right).
\end{align}
The desired result follows.

\hfill$\Box$

\textbf{Proof of Theorem \ref{thm2}: } By Proposition \ref{prop4} and Corollary \ref{cor2}. We only need to show that for any $\phi \in \mathcal{C}_c^+(\R \times S)$,
\begin{align}
	\E_{(0,i)} \left(\exp\left\{- C(\phi)M_{\sqrt{2\lambda^*}}(\infty)\right\} \right)= \E_{(0,i)}\left(\exp\left\{-\sum_{k,n} \phi\left(p_k+\Delta_n^{(k)}, q_n^{(k)} \right) \right\}\right).
\end{align}
Note that by Campbell's formula,
\begin{align}
	&\E_{(0,i)}\left(\exp\left\{-\sum_{k,n} \phi\left(p_k+\Delta_n^{(k)}, q_n^{(k)} \right) \right\}\right)\\
	& = \E_{(0,i)}\left(\exp\left\{-\sum_{k}\int  \phi\left(p_k+y, j \right)\mathcal{D}^{(k)}\left(\mathrm{d} y \mathrm{d} j\right) \right\}\right)\\
	& = \E_{(0,i)}\left(\prod_k \E_{(0,i)} \left(\exp\left\{-\int  \phi\left(z+y, j \right)\mathcal{D}\left(\mathrm{d} y \mathrm{d} j\right) \right\}\right)\big|_{z=p_k}\right)\\
	& =\E_{(0,i)} \left(\exp\left\{-\int_{\R} \left(1- \E_{(0,i)} \left(\exp\left\{-\int  \phi\left(z+y, j \right)\mathcal{D}\left(\mathrm{d} y \mathrm{d} j\right) \right\}\right)\right) C_\infty M_{\sqrt{2\lambda^*}}(\infty)\sqrt{2\lambda^*}e^{-\sqrt{2\lambda^*}z}\mathrm{d}z \right\}\right).
\end{align}
It suffices to show that for every $\phi \in \mathcal{C}_c^+(\R \times S)$,
\begin{align}\label{toprove-C}
	C(\phi) = C_\infty \int_{\R}  \left(1- \E_{(0,i)} \left(\exp\left\{-\int  \phi\left(z+y, j \right)\mathcal{D}\left(\mathrm{d} y \mathrm{d} j\right) \right\}\right)\right) \sqrt{2\lambda^*}e^{-\sqrt{2\lambda^*}z}\mathrm{d}z.
\end{align}
Suppose that $\phi(y,j)=0$ for all  $y\le m_\phi$ and $j\in S$.
Recalling that $Y$ is an exponential random variable with parameter $\sqrt{2\lambda^*}$ and $\overline{\mathcal{E}}_\infty = \mathcal{D}+Y$, we get that
\begin{align}\label{right1}
	&
\int_{\R}  \left(1- \E_{(0,i)} \left(\exp\left\{-\int  \phi\left(z+y, j \right)\mathcal{D}\left(\mathrm{d} y \mathrm{d} j\right) \right\}\right)\right) \sqrt{2\lambda^*}e^{-\sqrt{2\lambda^*}z}\mathrm{d}z \nonumber\\
	 & =
  \int_{m_\phi}^\infty  \left(1- \E_{(0,i)} \left(\exp\left\{-\int  \phi\left(z+y, j \right)\mathcal{D}\left(\mathrm{d} y \mathrm{d} j\right) \right\}\right)\right) \sqrt{2\lambda^*}e^{-\sqrt{2\lambda^*}z}\mathrm{d}z\nonumber\\
	 & =
e^{-\sqrt{2\lambda^*}m_\phi} \int_{0}^\infty   \left(1- \E_{(0,i)} \left(\exp\left\{-\int  \phi\left(z+y+m_\phi, j \right)\mathcal{D}\left(\mathrm{d} y \mathrm{d} j\right) \right\}\right)\right) \sqrt{2\lambda^*}e^{-\sqrt{2\lambda^*}z}\mathrm{d}z\nonumber\\
	 & =
e^{-\sqrt{2\lambda^*}m_\phi} \left(1- \E_{(0,i)} \left(\exp\left\{-\int \phi\left(Y+y+m_\phi, j \right)\mathcal{D}\left(\mathrm{d} y \mathrm{d} j\right) \right\}\right)\right) \nonumber\\
	 & =
 e^{-\sqrt{2\lambda^*}m_\phi} \left(1- \E_{(0,i)} \left(\exp\left\{-\int  \phi\left(y+m_\phi, j \right)\overline{\mathcal{E}}_\infty\left(\mathrm{d} y \mathrm{d} j\right) \right\}\right)\right).
\end{align}
Applying Proposition \ref{prop5} with $z=m_\phi$, we get
\begin{align}\label{right2}
	 &e^{-\sqrt{2\lambda^*}m_\phi} \left(1- \E_{(0,i)} \left(\exp\left\{-\int  \phi\left(y+m_\phi, j \right)\overline{\mathcal{E}}_\infty\left(\mathrm{d} y \mathrm{d} j\right) \right\}\right)\right) \nonumber\\
	 &= e^{-\sqrt{2\lambda^*}m_\phi}\lim_{t\to\infty} \left(1- \E_{(0,i)} \left(\exp\left\{-\int  \phi\left(y+m_\phi, j \right)\left(\overline{\mathcal{E}}_t-m_\phi\right)\left(\mathrm{d} y \mathrm{d} j\right) \right\}\Big| M_t>\sqrt{2\lambda^*}t +m_\phi \right)\right)\nonumber\\
	 & = 	e^{-\sqrt{2\lambda^*}m_\phi}\lim_{t\to\infty} \frac{ \E_{(0,i)} \left(1- \exp\left\{-\int  \phi\left(y, j \right) \overline{\mathcal{E}}_t \left(\mathrm{d} y \mathrm{d} j\right)\right\}; M_t>\sqrt{2\lambda^*}t +m_\phi \right)}{ \P_{(0,i)} \left(M_t>\sqrt{2\lambda^*}t +m_\phi \right)} 	\nonumber\\
	 & = 			e^{-\sqrt{2\lambda^*}m_\phi}\lim_{t\to\infty} \frac{ \E_{(0,i)} \left(1- \exp\left\{-\int  \phi\left(y, j \right) \overline{\mathcal{E}}_t \left(\mathrm{d} y \mathrm{d} j\right)\right\}\right)}{ \P_{(0,i)} \left(M_t>\sqrt{2\lambda^*}t +m_\phi \right)}.
\end{align}
By \eqref{step_42},
\begin{align}
	e^{-\sqrt{2\lambda^*}m_\phi}  \lim_{t\to\infty} \frac{\P_{(0,i )}\left(M_t > \sqrt{2\lambda^* }t \right) }{\P_{(0,i)} \left(M_t>\sqrt{2\lambda^*}t +m_\phi \right)} =1.
\end{align}
Therefore, by the probabilistic representation of $(v_2)_i(t,\sqrt{2\lambda^*}t)$ given by \eqref{Prob-Representation}, we continue \eqref{right2} to obtain
\begin{align}\label{right3}
	&e^{-\sqrt{2\lambda^*}m_\phi} \left(1- \E_{(0,i)} \left(\exp\left\{-\int  \phi\left(y+m_\phi, j \right)\overline{\mathcal{E}}_\infty\left(\mathrm{d} y \mathrm{d} j\right) \right\}\right)\right) \nonumber\\
	& = \lim_{t\to\infty} \frac{ \E_{(0,i)} \left(1- \exp\left\{-\int  \phi\left(y, j \right) \overline{\mathcal{E}}_t \left(\mathrm{d} y \mathrm{d} j\right)\right\}\right)}{ \P_{(0,i)} \left(M_t>\sqrt{2\lambda^*}t \right)} \nonumber\\
  &=  \lim_{t\to\infty}\frac{1}{\P_{(0,i)}\left( M_t> \sqrt{2\lambda^*}t  \right) }(v_2)_i(t,\sqrt{2\lambda^*}t)  = \frac{C(\phi)}{C_\infty},
\end{align}
where
$\mathbf{v}_2$ solves \eqref{F-KPP2} whose initial value is defined in
\eqref{Initial-v1-v2}
 and  in the last equality above we used  \eqref{Conditioned-Tightness3}.
Combining \eqref{right1} and \eqref{right3}, we get  \eqref{toprove-C}. The proof is now complete.

\hfill$\Box$

\noindent

\begin{singlespace}

\end{singlespace}

\vskip 0.2truein
\vskip 0.2truein

\noindent{\bf Haojie Hou:}  School of Mathematical Sciences, Peking
University,   Beijing, 100871, P.R. China. Email: {\texttt
houhaojie@pku.edu.cn}

\smallskip

\noindent{\bf Yan-Xia Ren:} LMAM School of Mathematical Sciences \& Center for
Statistical Science, Peking
University,  Beijing, 100871, P.R. China. Email: {\texttt
yxren@math.pku.edu.cn}

\smallskip
\noindent {\bf Renming Song:} Department of Mathematics,
University of Illinois at Urbana-Champaign,
Urbana, IL 61801, U.S.A.
Email: {\texttt rsong@illinois.edu}

\end{document}